\documentclass[12pt]{article}
\usepackage{theorem,amsfonts,amssymb,amsmath,latexsym}
\usepackage[mathscr]{eucal}
\usepackage[all]{xy}

\newcommand{\n}{\noindent}
\newcommand{\vn}{\vspace{5mm} \noindent}
\def\RR{{\mathbb R}}
\def\CC{{\mathbb C}}
\def\QQ{{\mathbb Q}}

\def\ZZ{{\mathbb Z}}

\def\GG{{\mathbb G}}
\def\AA{{\mathbb A}}
\def\DD{{\mathbb D}}
\def\FF{{\mathbb F}}

\def\GG{{\mathbb G}}
\def\BT1{{$\mbox{\rm BT}_{1}$}}
\def\DM1{{$\mbox{\rm DM}_{1}$}}
\def\Om{{\Omega}}

\newcommand{\Image}{\mbox{\rm Im}}
\newcommand{\dm}{\mbox{\rm dim}}

\newcommand{\inv}{\mbox{\rm inv}}
\newcommand{\degree}{\mbox{\rm deg}}
\newcommand{\sdm}{\mbox{\rm sdim}}

\def\diam{{\Diamond}}

\def\tr{{\triangle}}
\def\diam{{\Diamond}}

\newcommand{\va}{\varphi}

\def\cF{{\cal F}}
\def\cV{{\cal V}}

\def\cX{{\cal X}}
\def\cN{{\cal N}}
\def\cA{{\cal A}}
\def\cW{{\cal W}}
\def\cC{{\cal C}}

\def\cH{{\cal H}}
\def\cO{{\cal O}}
\def\cG{{\cal G}}
\def\cD{{\cal D}}

\newcommand{\B}{{\hfill$\Box$}}

\makeatletter
\def\setseccntfmt{\renewcommand{\@seccntformat}[1]{\S
    \csname the##1\endcsname.\hspace{1ex}}}
\def\setsubseccntfmt{\renewcommand{\@seccntformat}[1]{%
    (\csname the##1\endcsname)\hspace{0.5ex}}}
\renewcommand{\section}{\setseccntfmt\@startsection
  {section}{1}{0mm}{-\baselineskip}{0.5\baselineskip}{\sf\bfseries\Large}}
\def\presubsection{\setsubseccntfmt\@startsection
  {subsection}{2}{0mm}{-\baselineskip}{-0.5ex}{\bfseries\upshape}}
\def\tmpa{}\def\tmpb{}
\newcommand{\addspaceifnonempty}[1]{\def\tmpa{}\def\tmpb{#1}%
  \ifx\tmpa\tmpb{}\else{\hspace{0.5ex}#1\hspace{1ex}}\fi}
\renewcommand{\subsection}[1][]{\presubsection{\addspaceifnonempty{#1}}}
\def\presubsubsection{\setsubseccntfmt\@startsection
  {subsubsection}{3}{0mm}{-\baselineskip}{-0.5ex}{\bfseries\upshape}}
\renewcommand{\subsubsection}[1][]{\presubsubsection{\addspaceifnonempty{#1}}}
\gdef\th@nonumplain{\normalfont\itshape
  \def\@begintheorem##1##2{\item[\hskip\labelsep\theorem@headerfont ##1]}%
  \def\@opargbegintheorem##1##2##3{%
    \item[\hskip\labelsep \theorem@headerfont ##1\ (##3)]}}
\gdef\th@change{
  \def\@begintheorem##1##2{\item[\hskip\labelsep
    {\bfseries\upshape(##2)\hskip 1ex}\theorem@headerfont ##1]}%
  \def\@opargbegintheorem##1##2##3{\item[\hskip\labelsep
    {\bfseries\upshape(##2)\hskip 1ex}\theorem@headerfont ##1\ (##3)]}}
\makeatother   \def\nthm{\newtheorem}
\theoremheaderfont{\normalfont\bfseries\scshape}
{\theoremstyle{change}{\theorembodyfont{\normalfont\itshape}
\nthm{conj}[subsection]{Conjecture}  \nthm{*conj}[subsubsection]{Conjecture}
\nthm{thm}[subsection]{Theorem}      \nthm{*thm}[subsubsection]{Theorem}
\nthm{prop}[subsection]{Proposition} \nthm{*prop}[subsubsection]{Proposition}
\nthm{lemma}[subsection]{Lemma}      \nthm{*lemma}[subsubsection]{Lemma}
\nthm{cor}[subsection]{Corollary}    \nthm{*cor}[subsubsection]{Corollary}
}
{\theorembodyfont{\normalfont\rmfamily}
\nthm{rem}[subsection]{Remark.}       \nthm{*rem}[subsubsection]{Remark.}
\nthm{defn}[subsection]{Definition.}  \nthm{*defn}[subsubsection]{Definition.}
\nthm{ques}[subsection]{Question.}    \nthm{*ques}[subsubsection]{Question.}
}
}
{\theoremstyle{nonumplain}
{\theorembodyfont{\normalfont\itshape}
\nthm{thm*}[subparagraph]{Theorem}   \nthm{prop*}[subparagraph]{Proposition}
\nthm{lemma*}[subparagraph]{Lemma}   \nthm{cor*}[subparagraph]{Corollary}
\nthm{conj*}[subparagraph]{Conjecture}}{\theorembodyfont{\normalfont\rmfamily}
 \nthm{rem*}[subparagraph]{Remark}   \nthm{defn*}[subparagraph]{Definition}}}

{\theorembodyfont{\normalfont\rmfamily}
\theoremheaderfont{\normalfont\scshape}

}
\def\makeop#1{\expandafter\def\csname#1\endcsname
  {\mathop{\mathrm{#1}}\nolimits}\ignorespaces}
\makeop{Hom}   \makeop{End}   \makeop{Aut}   \makeop{Isom}  \makeop{Pic} 
\makeop{Gal}   \makeop{ord}   \makeop{Char}  \makeop{Div}   \makeop{Lie} 
\makeop{PGL}   \makeop{Corr}  \makeop{PSL}   \makeop{sgn}   \makeop{Spf}
\makeop{Spec}  \makeop{Tr}    \makeop{Nr}    \makeop{Fr}    \makeop{disc}
\makeop{Proj}  \makeop{supp}  \makeop{ker}   \makeop{im}    \makeop{dom}
\makeop{coker} \makeop{Stab}  \makeop{SO}    \makeop{SL}    \makeop{SL}
\makeop{Cl}    \makeop{cond}  \makeop{Br}    \makeop{inv}   \makeop{rank}
\makeop{id}    \makeop{Fil}   \makeop{Frac}  \makeop{GL}    \makeop{SU}
\makeop{Trd}   \makeop{Sp}    \makeop{Tr}    \makeop{Trd}
\makeop{CSp}   \makeop{diag}
\makeop{Res}   \makeop{ind}   \makeop{depth} \makeop{Tr}    \makeop{st}
\makeop{Ad}    \makeop{Int}   \makeop{tr}    \makeop{Sym}   \makeop{can}
\makeop{length}\makeop{SO}    \makeop{torsion}
\makeop{Ext}   \makeop{Tor}   \makeop{Def}
\makeop{fr}

\def\makebb#1{\expandafter\def
  \csname bb#1\endcsname{{\mathbb{#1}}}\ignorespaces}
\def\makebf#1{\expandafter\def\csname bf#1\endcsname{{\bf #1}}\ignorespaces}
\def\makegr#1{\expandafter\def
  \csname gr#1\endcsname{{\mathfrak{#1}}}\ignorespaces}
\def\makescr#1{\expandafter\def
  \csname scr#1\endcsname{{\mathscr{#1}}}\ignorespaces}
\def\makecal#1{\expandafter\def\csname cal#1\endcsname{{\cal #1}}\ignorespaces}
\def\doLetters#1{#1A #1B #1C #1D #1E #1F #1G #1H #1I #1J #1K #1L #1M
                 #1N #1O #1P #1Q #1R #1S #1T #1U #1V #1W #1X #1Y #1Z}
\def\doletters#1{#1a #1b #1c #1d #1e #1f #1g #1h #1i #1j #1k #1l #1m
                 #1n #1o #1p #1q #1r #1s #1t #1u #1v #1w #1x #1y #1z}
\doLetters\makebb   \doLetters\makecal  \doLetters\makebf \doLetters\makescr
\doletters\makebf   \doLetters\makegr   \doletters\makegr

\newcommand{\rmitem}[1][{}]{\item[{\rm #1}]}
\newcommand{\tworows}[2]{\genfrac{}{}{0pt}{}{#1}{#2}}

\def\Zp{{\mathbb{Z}}_p}
\def\Qp{{\mathbb{Q}}_p}
\def\Zl{{\mathbb{Z}}_{\ell}}
\def\Ql{{\mathbb{Q}}_{\ell}}

\def\Fp{{\mathbb{F}_p}}

\def\Fpbar{{\mathbb{F}}}

\def\Cart#1{{\rm Cart}(#1)}
\def\Cartp#1{{\rm Cart}_p(#1)}

\def\Ker{{\rm Ker}}

\def\rk{{\rm rk}}
\def\ringO{{\mathscr{O}}}

\def\ZZ{\bbZ}
 
\def\Spec{\mbox{\rm Spec}}
\def\Spec{\mbox{\rm Spec}}
\def\cF{{\cal F}}

\def\Bl{\fbox{\bf BB}}
\def\Th{\fbox{\bf Th}}
\def\Ex{\fbox{\bf Extra}}

\def\proof{\medbreak\noindent{\scshape Proof.}\enspace}%
     
\def\qed{\qedmark\medbreak} 
\def\qedmark{{\enspace\raisebox{3pt}%
{\framebox[5pt][l]{\raisebox{0pt}[0.5pt][.5pt]{}}}}}

\begin{document}


\begin{center}{\Large\scshape Moduli of abelian varieties and 
$p$-divisible groups:\\
{\normalsize Density of Hecke orbits, and a conjecture by Grothendieck}}\\
\bigbreak
{\large Ching-Li Chai\footnote{Partially supported by a grant DMS04-00482
from the National Science Foundation}
\&
Frans Oort
 }\\
\end{center}

\hfill{\normalsize\bf Conference on Arithmetic Geometry,}

\hfill{\normalsize\bf G\"ottingen 17 July  - 11 August 2006.}

\bigbreak\bigbreak

\vn
In the week 7 -- 11 August 2006 we gave a course, and here are notes for 
that course. Our main topic is: {\it  geometry and arithmetic  
of $\cA_g \otimes \FF_p$,  the moduli space of polarized abelian 
varieties of dimension $g$  in positive characteristic}. 
We illustrate properties, and some of the available techniques by 
treating two topics:

\begin{center}
{\bf Density of ordinary Hecke orbits}
\end{center}
and
\begin{center}
{\bf A conjecture by Grothendieck on deformations of 
$p$-divisible groups.}
\end{center}
\n
{\bf Contents:}

\vn
{\bf 1.} {\it Introduction:}
{\it Hecke orbits, and the Grothendieck conjecture.}
\\
{\bf  2.}  
{\it  Serre-Tate theory.}
\\
{\bf 3.}  
{\it    The Tate-conjecture: $\ell$-adic and 
$p$-adic.}
\\
{\bf 4.} 
{\it   Dieudonn\'e modules and Cartier modules.}
\\
{\bf 5.} {\it Cayley-Hamilton:}
{\it    --   a conjecture by Manin and the weak 
Grothendieck conjecture.}
\\
{\bf 6.}  
{\it    Hilbert modular varieties.  }
\\
{\bf  7.} 
{\it   Deformations of $p$-divisible groups to 
$a \leq 1$.}
\\
{\bf 8.}  
{\it   Proof of the Grothendieck conjecture.}
\\
{\bf 9.} 
{\it    Proof of the density of ordinary 
Hecke orbits.}
\\
{\bf 10.}
{\it Notations and some results used.}
\\{\bf 11.}
{\it A remark and a question.}


\vn
We  present proofs of two recent results. 
The main point is that 
the {\it methods} used for these proofs are interesting. 
The emphasis will be to discuss various 
techniques available.

\vn
In characteristic zero we have strong tools at our disposal: 
besides algebraic-geometric theories we can use  analytic  
and topological methods. It seems that we are at a loss in 
positive characteristic. However the opposite is true. 
Phenomena, only occurring in positive characteristic  
provide us with strong tools to study moduli spaces. 
And, as it turns out again and again, several results in 
characteristic zero 
can be derived using reduction modulo $p$. It is about these tools 
in positive characteristic that will be the focus of our talks.

\vn
Here is a list of some of the central topics:

Serre-Tate theory.

Abelian varieties over finite fields.

Monodromy: $\ell$-adic and $p$-adic, geometric an arithmetic.

Dieudonn\'e modules and Newton polygons.

Theory of Dieudonn\'e modules, Cartier modules and  displays.

Cayley-Hamilton and deformations of $p$-divisible groups.

Hilbert modular varieties.

Purity of the Newton polygon stratification in families of 
$p$-divisible groups.

\vn
The strategy  is that we have chosen certain central topics, 
and for those we took ample time for explanation and for proofs. 
Besides that we need certain results which we label as ``Black Box''. 
These are results which we need for our proofs, which are either 
fundamental theoretical results (but it would take too much time to 
explain their proofs), or it concerns lemmas which are computational, 
important for the proof, but not very interesting to explain in a course. 
We hope that we explain well enough what every relevant statement is. 
 We write:

\vn
\Bl \ \  A Black Box, please accept that this result is true.
\smallskip

\n
\Th \ \  This is one of the central results, and we will explain.
\smallskip

\n
\Ex \ \  This is a result, which is interesting, but was not 
 discussed in the course.

\vn
Notation to be used will be explained in Section \ref{Notat}. 
In order to be somewhat complete 
we will gather related interesting other results and questions and  
conjectures in Section \ref{ques}. 
Part of our general convention is that \(K\) denotes a field of
characteristic \(p>0\), and \(k\) denotes an algebraically closed
field containing \(K\).

\section{Introduction: Hecke orbits, and the Grothendieck 
conjecture}\label{1}
\n
In this section we discuss the two theorems we are going to consider.

\vn
\underline{\bf Hecke orbits.}
\subsection{} An abelian variety $A$ of dimension $g$ over a field $K$ is 
called {\it ordinary} if 
$$
\#\left(A[p](k)\right) = p^g.
$$
More  generally, 
\begin{center}
the number $f$ such that $\#\left(A[p](k)\right) = p^f$ is called the 
$p$-rank of $A$, 
\end{center}
and $f=g$ is the case of ordinary abelian varieties.
\smallbreak

We say an elliptic curve $E$ is {\it supersingular} if it is not 
ordinary. 
Equivalently: $E$ is supersingular if $E[p](k) = 0$. 
This terminology stems 
from Deuring; explanation: an elliptic curve in characteristic zero 
is said to determine a {\it singular} $j$-value if the 
endomorphism ring over an algebraically closed field 
is larger than $\ZZ$, 
in fact of rank $2$ over $\ZZ$; a {\it supersingular} elliptic curve 
$E$ over $k$ has $\rk_{\ZZ}(\End(E)) = 4$. Note that an elliptic curve 
is non-singular; a more correct terminology would be: 
an elliptic curve with a supersingular $j$-invariant, 
and for a an elliptic curve over $\CC$ with endomorphism ring 
not equal to $\ZZ$: an elliptic curve with a singular $j$-invariant.
\smallbreak

We say an   abelian variety $A$ of dimension $g$ over a field $K$ is   
{\it supersingular} if there exists an isogeny 
$A \otimes_K k  \sim E^g$, 
where $E$ is a supersingular elliptic curve. 
An equivalent condition will be 
given later, and more explanation will follow as soon as we have 
Dieudonn\'e modules and the theory of Newton polygons at our disposal, 
see Section \ref{4}.
\smallbreak

Note that $f(A) = 0$ does not imply the abelian variety is supersingular 
in case $\dim(A) \geq 3$.
\smallbreak

Let \(A\) and \(B\) be abelian varieties over a field \(K\).
A \(\bbQ\)-isogeny, also called a \emph{quasi-isogeny}, 
from \(A\) to \(B\),
is an element of \(\Hom(A,B)\otimes_{\bbZ}{\bbQ}\)
which has an inverse in
\(\Hom(B,A)\otimes_{\bbZ}{\bbQ}\); any element
\(\phi\in \Hom(A,B)\otimes_{\bbQ} \)
can be realized by a diagram
\(
A \xleftarrow{\alpha}C \xrightarrow{\beta}B
\)
where \(\alpha,\beta\) are isogenies of abelian varieties.
Below is a more general definition.

\subsection{\bf Definition.} Let $\Gamma \subset \QQ$ be a subring. 
A \(\Gamma\)-isogeny from \(A\) to \(B\) is an element 
\(\psi\in \Hom(A,B)\otimes_{\bbZ}\Gamma\) such that 
there exists an element \(\psi'\in \Hom(B,A)\otimes_{\bbZ}\Gamma\)
with the property that \(\psi'\circ \psi\in \Gamma^{\times}\cdot {\rm Id}_A\)
and \(\psi\circ\psi'\in \Gamma^{\times}\cdot {\rm Id}_B\).

\noindent{\bf Remark}.\ \enspace
(i) A \(\Gamma\) isogeny \(\psi\) as above can be realized by a diagram
\(\,
A \xleftarrow{\alpha}C \xrightarrow{\beta}B\,
\), where
\(\alpha\) and \(\beta\) are isogenies such that there
exists an integer \(N\in\Gamma^{\times}\) such that
\(N\cdot \Ker(\alpha)=N\cdot \Ker(\beta)=0\).

(ii) When \(\Gamma=\bbZ_{(p)}\) (resp. \(\bbZ[1/\ell]\)),
we say that \(\psi\) is a \emph{prime-to-\(p\) isogeny}
(resp.\ an \emph{\(\ell\)-power isogeny}).
Here \(\bbZ_{(p)}=\bbQ\cap \bbZ_p\) is the localization of \(\bbZ\) at
the prime ideal \((p)=p\bbZ\).

\subsection{\bf Definition.} \label{def_ho}
Let  $[(A,\lambda)] = x \in \cA_g$ be the moduli 
point of a polarized abelian variety over a field $K$. 
We say $[(B,\mu)] = y$ 
is in the {\it Hecke orbit} of $x$  if there exists a field $\Om$, and 
\begin{center}
a \(\bbQ\)-isogeny $\va: A_{\Om}\to B_{\Om}$   
such that 
$\va^{\ast}(\mu) = \lambda$.
\end{center} 
{\bf Notation:} $y \in \cH(x)$. The set  $\cH(x)$ is called the 
{\it Hecke orbit} of $x$.

\vn
{\bf Hecke-prime-to-$p$-orbits.} If in the previous definition moreover 
$\va$ is a \(\bbZ_{(p)}\)-isogeny,
we say $[(B,\mu)] = y$ is in the  Hecke-prime-to-$p$-orbit of $x$. \\
{\bf Notation:} $y \in \cH^{(p)}(x)$.

\vn
{\bf Hecke-$\ell$-orbits.} Fix a prime number $\ell$ different from $p$. 
We say $[(B,\mu)] = y$ is in the  Hecke-$\ell$-orbit of $x$ if 
in the previous definition moreover \(\va\) is a \(\bbZ[1/\ell]\)-isogeny.
\\
{\bf Notation:} $y \in \cH_{\ell}(x)$.

\vn
{\bf Remark.} We have given the definition of the  so-called 
$\Sp_{2g}$-Hecke-orbit. On can also define the (slightly bigger) 
$\CSp_{2g}(\ZZ)$-Hecke-orbits by the usual Hecke correspondences, 
see  \cite{F.C}, VII.3, also see \ref{Isog} below,
$$
\cH^{\rm Sp}(x) = \cH(x) \subset \cH^{\rm CSp}(x). 
$$

\vn
{\bf Remark.} Note that $y \in \cH(x)$ is equivalent by requiring 
the existence of a diagram
$$
(B,\mu)  \stackrel{\psi}{\longleftarrow}  
(C,\zeta) \stackrel{\va}{\longrightarrow} (A,\lambda).
$$
where  $[(B,\mu)] = y$ and   
 $[(A,\lambda)] = x$.

\vn
{\bf Remark.} Suppose $y \in \cH^{(p)}(x)$ and suppose 
$\deg(\lambda) = \deg(\mu)$; then $\deg(\va)$ and 
$\deg(\psi)$ are not divisible by $p$, 
which explains the terminology ``prime-to-$p$''.

\vn
{\bf Remark.} The diagrams which define $\cH(x)$ as above give representable 
correspondences between components of the moduli scheme; 
these correspondences could be denoted by $\Sp$-Isog, 
whereas the correspondences considered in  \cite{F.C}, VII.3 
could be denoted by $\CSp$-Isog.

\subsection{\bf Remark/Exercise.} (Characteristic zero.) 
The Hecke orbit of 
a point in the moduli space $\cA_g \otimes \CC$ in 
{\it characteristic zero} is dense in that moduli space 
(dense in the classical topology, dense in the Zariski topology).

\subsection{\bf Hecke orbits of elliptic curves.} Consider the 
moduli point 
$[E]= j(E) = x \in \cA_{1,1} \cong \AA^1$ of an elliptic curve in 
characteristic $p$. Here $\cA_{1,1}$ stands for  $\cA_{1,1} \otimes \FF_p$. 
Note that every elliptic curve has a 
unique principal polarization.\\
(1) {\bf Remark.} {\it If $E$ is supersingular 
$\cH(x) \cap \cA_{1,1}$ 
is a finite set; we conclude that 
$\cH(x)$ is nowhere dense in $\cA_1$.} 

Indeed, the supersingular locus in $\cA_{1,1}$ is closed, 
there do exist ordinary elliptic curves, hence that locus is finite; 
Deuring and Igusa computed the exact number of geometric points 
in this locus. \\
(2) {\bf Remark/Exercise.} {\it If $E$ is ordinary, its 
Hecke-$\ell$-orbit is dense in 
$\cA_{1,1}$.} There are several ways of proving this. Easy and direct 
considerations show that in this case $\cH_{\ell}(x) \cap \cA_{1,1}$ 
is not finite, note that every component of $\cA_1$ has dimension one; 
conclude $\cH(x)$ is dense in $\cA_1$.\\
{\bf Remark.} More generally in fact, as we see in \cite{Chai-Density}, 
Proposition 1 on page 448:   $\cH_{\ell}(x) \cap \cA_{g,1}$ 
{\it is finite if and only if $[(A,\lambda)] = x \in \cA_g$ where 
$A$ is supersingular.}\\
{\bf Remark.} For elliptic curves we have defined 
(supersingular)$ \Leftrightarrow$ (non-ordinary). 
For $g=2$ on can see that 
(supersingular) $\Leftrightarrow$ ($f=0$). 
However, see Section \ref{4}, we can define supersingular 
as those abelian varieties where the Newton polygon has all slopes equal 
to $1/2$; for  $g>2$ there do exist abelian varieties of $p$-rank zero 
which are not supersingular.  

\subsection{}\label{Isog}  {\bf A bigger Hecke orbit.} We define the
notion of \({\rm CSp}\)-Hecke orbits.  Two \(K\)-points 
\([(A,\lambda)]\), \([(B,\mu)]\) of \(\cA_{g,1}\) are in the
same \({\rm CSp}\)-Hecke orbit 
(resp.\ prime-to-\(p\) \(\rm CSp\)-Hecke orbit,
resp.\ \(\ell\)-power ${\rm CSp}$-Hecke orbit)  
if there exists an isogeny 
\(\va: A \otimes k \to B\otimes k\) and a positive integer \(n\)
(resp.\ a positive integer \(n\) which is relatively prime to \(p\),
resp.\ a positive integer which is a power of \(\ell\))
such that \(\va^*(\mu)=n\cdot \lambda\).
Such Hecke correspondences are representable by morphisms
${\rm Isog}_g \subset \cA_g \times \cA_g$ on $\cA_g$, 
also see \cite{F.C}, VII.3.

\vn
The set of all such  $(B,\mu)$ for a fixed $x := [(A,\lambda)]$ 
is called the $\CSp$-Hecke orbit  
(resp.\ \(\CSp(\bbA_f^{(p)})\)-Hecke orbit 
resp.\  \(\CSp(\Ql)\)-Hecke orbit) of $x$;
notation $\cH^{\rm Sp}(x)$
(resp.\ \(\cH^{(p)}_{\rm Sp}(x)\),
resp.\ \(\cH^{\CSp}_{\ell}(x)\).)
Note that $\cH^{\CSp}(x) \supset \cH(x)$. This slightly bigger 
Hecke orbit will play no role in this paper. 
However it is nice to see the relation between the Hecke orbit defined 
previously in \ref{def_ho}, which could be called the 
$\Sp$-Hecke orbits and \(\Sp\)-Hecke correspondences,
with the \(\CSp\)-Hecke orbits and \(\CSp\)-Hecke correspondences.

\begin{thm}\label{dens} \Th \ 
{\rm ({\bf Density of ordinary Hecke orbits.})} Let $[(A,\lambda)] = x$ 
be the moduli point of a polarized {\rm ordinary} abelian variety. 
Let $\ell$ be a prime number different from $p$. 
The Hecke-$\ell$-orbit $\cH_{\ell}(x)$ is dense in $\cA_{g,1}$:
$$ \left(\cH_{\ell}(x) \cap \cA_{g,1}\right)^{\rm Zar} = \cA_{g,1}.
$$
From this we conclude: $\cH(x)$ is dense in $\cA_g$.
\end{thm}

\n
See Theorem \ref{dens_siegel}.
This theorem was proved by Ching-Li Chai in 1995, 
see \cite{Chai-Density}, 
Theorem 2 on page 477. Although \(\CSp\)-Hecke orbits were used in 
\cite{Chai-Density}, the same argument works for 
\(\Sp\)-Hecke orbits as well.
We  present a proof of this theorem; 
we  follow  \cite{Chai-Density}  partly, but also present new insight 
which was necessary for solving the general Hecke orbit problem. 
This final strategy will provide us with a proof which seems easier than 
the one given previously. More information on the general Hecke orbit 
problem can be obtained from \cite{Chai-HOsurvey} as long as  
\cite{Chai.FO-HO} is not yet available.

\subsection{\bf Exercise.} (Any characteristic.) Let $k$ be any 
algebraically closed field (of any characteristic). 
Let $E$ be an elliptic curve 
over $k$ such that $\End(E) = \ZZ$. Let $\ell$ be a prime number 
different from the characteristic of $k$. 
Let $E'$ be an elliptic curve such that there exists an isomorphism 
$E'/\underline{(\ZZ/\ell)}_k \cong E$. 
Let $\lambda$ be the principal polarization on $E$, 
let $\mu$ be the pull back of $\lambda$ to $E'$,
 hence $\mu$ has degree $\ell^2$, and 
let $\mu' = \mu/\ell^2$, hence $\mu'$ is a principal polarization 
on $E'$. 
{\it Remark that $[(E',\mu')] \in \cH(x)$. 
Show that   $[(E',\mu')] \not\in \cH^{\rm Sp}(x)$.}

\subsection{\bf Exercise.} Let $E$ be an elliptic curve 
in characteristic $p$ 
which is not supersingular (hence ordinary); 
let $\mu$ be any polarization 
on $E$, and $x:= [(E,\mu)]$. Show $\cH^{\rm Sp}(x)$ is dense in $\cA_1$. 

\begin{thm}\label{dualth}  
{\rm (Duality theorem for abelian schemes, 
see \cite{FO-CGS}, Theorem 19.1)} 
Let $\va: B \to A$ be an isogeny of abelian schemes. 
We obtain an exact sequence
$$
0 \quad\to\quad \Ker(\va)^D \quad\longrightarrow\quad A^t 
\quad\stackrel{\va^t}{\longrightarrow}\quad B^t \quad\to\quad 0. 
$$
\end{thm} 

\subsection{\bf An example.} Write 
$\mu_s = \Ker(\times s: \GG_m \to \GG_m)$ 
for every $s \in \ZZ_{>0}$. It is not difficult to see that 
$(\mu_p)^D = \underline{\ZZ/p}$, and in fact, 
$(\mu_{p^b})^D = \underline{\ZZ/{p^b}}$.\\
{\bf Conclusion.} {\it For an ordinary abelian variety $A$ 
over $k$, we have}
$$A[p] \quad\cong\quad \left((\mu_p)^g\right) \times 
\left((\underline{\ZZ/p})^g\right). $$
In fact, by definition we have that $A[p](k) \cong (\ZZ/p)^g$.  
This implies that  $(\underline{\ZZ/p^b})^g \subset A$. 
By the duality theorem we have $A[p]^D \subset A^t$. 
Hence  $(\mu_{p^b})^g \subset A^t$. 
As $A$ admits a polarization we have an isogeny $A \sim A^t$. 
We conclude that $(\mu_{p})^g \subset A$. Hence the result.
\qed   

\vn
Note that for an ordinary abelian variety $A$  over an arbitrary field 
$K$ the Galois group $\Gal(K^{\rm sep}/K)$ acts on $A[p]^{\rm loc}$ 
and on $A[p]^{\rm et} =  A[p]/A[p]^{\rm loc}$, and these actions 
need not be trivial. Moreover if $K$ is not perfect, the extension 
$0 \to A[p]^{\rm loc} \to A[p] \to  A[p]^{\rm et} \to 0$ 
need not be split; this will be studied extensively in Section \ref{2}.

\subsection{\bf Reminder.} Let $N$ be a finite group scheme over a field 
$K$. Suppose that the rank of $N$ is is prime to the characteristic 
of $k$. Then $N$ is\'etale over $K$; e.g. see \cite{FO-reduced}.

\subsection{} (1) Write ${\rm Isog} =  \Sp$-${\rm Isog}$. 
Consider a component $I$  of ${\rm Isog}_g$ defined by  
diagrams as in \ref{Isog} with $\deg(\psi) = b$ and $\deg(\va) = c$.
If $b$ is not divisible by $p$, the first projection 
$\cA_g \leftarrow I$ is\'etale; if $c$ is not divisible by $p$, 
 the second projection $I \to \cA_g$ is\'etale. \\
(2) Consider ${\rm Isog}_g^{\rm ord} \subset {\rm Isog}_g$, 
the largest subscheme (it is locally closed) lying over the ordinary 
locus (either in the first projection, or in the second projection, 
that is the same).\\
{\bf Exercise.} {\it The two projections 
$(\cA_g)^{\rm ord} \leftarrow {\rm Isog}_g^{\rm ord} 
\to (\cA_g)^{\rm ord}$ 
are both surjective,\'etale, finite and flat.
}\\
(3) \Ex  \ \ {\it The projections   $(\cA_g) \leftarrow 
{\rm Isog}_g \to (\cA_g)$ 
are both surjective and proper on every irreducible component of 
${\rm Isog}_g$}; this follows from  \cite{F.C}, VII.4. 
The previous exercise (2) is easy; fact (3) is more difficult; 
it uses the computation in \cite{Norman.FO}.

\subsection{}\Bl \ \ In \cite{Norman.FO} it has been proved: 
$(\cA_g)^{\rm ord}$ {\it is dense in} $\cA_g$.

\subsection{} {\it We see that for an {\rm ordinary}  
$[(A,\lambda)] = x$ we have:}
$$
\left(\cH_{\ell}(x) \cap \cA_{g,1}\right)^{\rm Zar} = \cA_{g,1} 
\quad\Longrightarrow\quad \left(\cH(x)\right)^{\rm Zar} = \cA_g.
$$
Work over $k$. In fact, consider an irreducible component 
$T$ of $\cA_g$. As proven in \cite{Norman.FO} there is an ordinary point 
$y = [(B,\mu)]  \in T$. By \cite{Mumford-AV}, Corollary 1 on page 234, 
we see that there is an isogeny  $(B,\mu) \to (A,\lambda)$, 
where $\lambda$ is a principal polarization. By (2) in the previous 
subsection we see that density of  $\cH_{\ell}(x) \cap \cA_{g,1}$ in 
$\cA_{g,1}$ implies density of  $\cH_{\ell}(x) \cap T$ in $T$.
\qed

\smallbreak\n
Therefore, from now on we shall be mainly interested in Hecke orbits 
in the principally polarized case.

\subsection{\bf Theorem}\label{HdenseNP} \ \Ex \  
(Ching-Li Chai and Frans Oort). 
{\it For any $[(A,\mu)] = x \in \cA_g \otimes \FF_p$ with 
$\xi = \cN(A)$,  
the Hecke orbit $\cH(x)$ is dense in the Newton polygon locus 
$\cW_{\xi}(\cA_g \otimes \FF_p)$.}\\
A proof will be presented in \cite{Chai.FO-HO}. \\
Note that in case $f(A) \leq g-2$ the $\ell$-Hecke orbit is not dense in 
$\cW_{\xi}(\cA_g \otimes \FF_p)$. In \cite{FO-Fol} we find 
a precise conjectural description of the Zariski closure of 
$\cH_{\ell}(x)$; that conjecture has been proved now, 
and it implies \ref{HdenseNP}.

\vn
\subsection{\bf Lemma}\label{ssfinite} \Bl \ \  (Chai). 
{\it Let $[(A,\lambda)] = x \in \cA_{g,1}$. 
Suppose that $A$ is supersingular} (i.e. over an algebraically closed field 
$A$ is isogenous with a product of supersingular elliptic curves, 
equivalently: all slopes in the Newton polygon $\cN(A)$ are equal to $1/2$). 
{\it Then}
$$\cH^{(p)}(x) \cap \cA_{g,1} \quad \mbox{\it is finite}.$$
See \cite{Chai-Density}, Proposition 1 on page 448. \\
Note that $\cH(x)$ equals {\it the whole 
supersingular Newton polygon stratum: 
the prime-to-$p$ Hecke orbit is small, but the Hecke orbit including 
$p$-power quasi-isogenies is large}.

\vn\vn
\underline{\bf A conjecture by Grothendieck.}

\subsection{\bf Definition,  $p$-divisible groups.} Suppose given 
$h \in \ZZ_{>0}$. Suppose given a base scheme $S$. 
Suppose given for every $i \in \ZZ_{>0}$ a finite, flat group scheme 
$G_i \to S$ of rank $p^{ih}$, and inclusions 
$G_i \subset G_{i+1}$ for every $i$  such that $G_{i+1}[p^i] = G_i$. 
The inductive system $X=\{G_i \mid i\} \to S$ is called a 
$p$-divisible group of height $h$ over S.\\
The notion of a Barsotti-Tate group, or BT-group, is the same as that of 
a $p$-divisible group; both terminologies will be used. 
For more information see \cite{Ill}, Section 1.

\vn
{\bf Remark.} Note that for every $j$ and every $s\geq 0$ the map 
$\times p^s$ induces a surjection:
$$
G_{j+s} \twoheadrightarrow G_j = G_{j+s}[p^j] \subset G_{j+s}.
$$
For every $i$ and $s\geq 0$ we have an exact sequence of 
finite flat group schemes
$$
0 \to G_i \longrightarrow G_{i+s} \longrightarrow G_s \to 0.
$$
{\bf Example.} Let $A \to S$ be an abelian scheme. For every $i$ we write 
$G_i = A[p^i]$. The inductive system  $G_i \subset G_{i+s} \subset A$  
defines a $p$-divisible group of height $2g$. We shall denote this by 
$X = A[p^{\infty}]$ (although of course ``$p^{\infty}$'' 
strictly speaking is not defined). 

\vn
For $p$-divisible groups, inductive systems, we define homomorphisms by 
$$
\Hom(\{G_i\},\{H_j\}) \quad=\quad {\rm lim}_{\leftarrow i} \ 
\ {\rm lim}_{j \rightarrow} \ \Hom(G_i,H_j).
$$
Note that a homomorphism $A \to B$ of abelian schemes defines a morphism 
$A[p^{\infty}] \to B[p^{\infty}]$ of $p$-divisible groups.

\subsection{\bf Discussion.} Over any base scheme $S$ (in any 
characteristic) 
for an abelian scheme $A \to S$ and for a prime number $\ell$ invertible 
on $S$ one can define $T_{\ell}(A/S)$ as follows. 
For $i \in \ZZ_{>0}$ one chooses $N_i := A[\ell^i]$, 
and we  define $\times \ell: N_{i+1} \to N_i$. 
This gives a {\it projective system}, and we write 
$$
T_{\ell}(A/S) \quad=\quad \{A[\ell^i] \mid i\in \ZZ_{>0}\} 
\quad=\quad  {\rm lim}_{\leftarrow i} \ A[\ell^i].
$$
This is called the $\ell$-Tate group of $A/S$. Any geometric fiber is 
$T_{\ell}(A/S)_s \cong (\underline{\ZZ_{\ell}}_s)^{2g}$. 
If $S$ is the spectrum of a field $K$, the Tate-$\ell$ can be considered 
as the structure of a $\Gal(K^{\rm sep}/K)$-module   
on the group $\ZZ_{\ell}^{2g}$, see \ref{et}. 
{\it One should like to have an analogous concept for this 
notion in case $p$ 
is not invertible on $S$.} This is precisely the role of 
$A[p^{\infty}]$ 
defined above. Historically a Tate-$\ell$-group is defined as a 
projective system, and the $p$-divisible group as an inductive system; 
it turns out that these are the best ways of handling these concepts 
(but the way in which direction to choose the limit is not very 
important). 
We see that the $p$-divisible group of an abelian variety 
should be considered as the natural substitute for the Tate-$\ell$-group. 

\vn
In order to carry this analogy further we investigate aspects of  
$T_{\ell}(A)$ and wonder whether these can be carried over to  
$A[p^{\infty}]$. The first is a twist of a pro-group scheme 
defined over $\Spec(\ZZ)$. What can be said in analogy about $A[p^{\infty}]$? 
We will see that {\it up to isogeny}  $A[p^{\infty}]$ is a twist of 
an ind-group scheme over $\FF_p$; however ``twist'' here 
should be understood not only in the sense of separable Galois theory, 
but also using inseparable aspects: the main idea of Serre-Tate 
parameters, to be discussed in Section \ref{2}.

\subsection{\bf The Serre dual of a $p$-divisible group.} 
Consider a $p$-divisible group $X=\{G_i\}$. The exact sequence 
$G_{j+s}/G_j = G_s$ by Cartier duality, see \cite{FO-CGS}, I.2, 
defines an exact sequence
$$
0 \to G_s^D \longrightarrow G_{j+s}^D \longrightarrow G_j^D \to 0.
$$
These are  used, in particular the inclusions 
$$G_s^D \hookrightarrow G_{s+1}^D \quad=\quad  \left(G_{s+1} 
\twoheadrightarrow G_s \right)^D,$$
to define the $p$-divisible group $X^t =\{G_s^D\}$, 
called {\it the Serre dual} of $X$. Using \ref{dualth} we conclude: 
$(A[p^\infty])^t = A^t[p^{\infty}]$ (which is less trivial than 
notation suggests...).

\vn
In order to being able to handle the isogeny class of $A[p^{\infty}]$
 we need the notion of Newton polygons.

\subsection{\bf Newton polygons.} Suppose given integers 
$h, d \in \ZZ_{\geq 0}$;  
here $h$ = ``height'', $d$ = ``dimension''. In case of abelian 
varieties 
we will choose $h = 2g$, and $d = g$. A Newton polygon $\gamma$ 
(related to $h$ and $d$) is a polygon $\gamma \subset \QQ \times \QQ$ 
(or, if you wish in $\RR \times \RR$), such that:

\begin{picture}(250,90)(0,0)
\put(290,0){\line(6,1){60}}
\put(350,10){\line(2,1){60}}
\put(410,40){\line(1,1){30}}

\put(290,0){\line(1,0){150}}
\put(290,0){\line(0,1){80}}

\put(290,0){\circle*{2}}
\put(350,10){\circle*{2}}
\put(410,40){\circle*{2}}
\put(440,70){\circle*{2}}

\put(430,5){$h$}
\put(295,70){$d$}
\put(380,35){$\zeta$}

\put(0,60){$\bullet$ \quad $\gamma$ starts at $(0,0)$ and  
ends at $(h,d)$;
}
\put(0,40){$\bullet$ \quad $\gamma$ is lower convex;}
\put(0,20){$\bullet$ \quad any slope $\beta$ of $\gamma$ has the property 
$0 \leq \beta \leq 1$;}
\put(0,0){$\bullet$ \quad the breakpoints of $\gamma$ are 
in $\ZZ \times \ZZ$; hence $\beta \in \QQ$.}

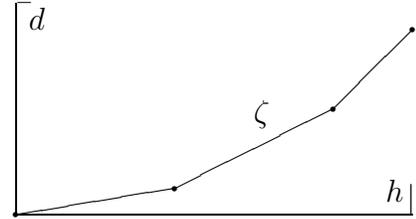
\put(290,80){$_{-}$}
\put(440,2.5){\!\!\! $\mid$}


\end{picture}

\vn
Note that a Newton polygon determines (and is determined by) 
\begin{center}
$\beta_1, \cdots , \beta_h \in \QQ$ with 
$0 \leq \beta_1 \leq \cdots \leq \beta_h \leq 1 \quad\leftrightarrow
\quad \zeta$.
\end{center}

\vn
Sometimes we will give a Newton polygon by data $\sum_i \ (m_i,n_i)$; 
here $m_i, n_i \in \ZZ_{\geq 0}$, with $\gcd(m_i, n_i) = 1$, 
and $m_i/(m_i+n_i) \leq m_j/(m_j+n_j)$  for $i \leq j$, 
and $h = \sum_i \ (m_i+n_i)$, \  $d= \sum_i m_i$.  From these data we 
construct the related Newton polygon by choosing the slopes 
$m_i/(m_i+n_i)$ with multiplicities  $h_i = m_i+n_i$. 
Conversely clearly any Newton polygon 
can be encoded in a unique way in such a form.

\vn {\bf Remark. The Newton polygon of a polynomial.} 
Let $g \in \QQ_p[T]$ be a 
monic polynomial of degree $h$. We are interested in the 
$p$-adic values of its zeroes (in an algebraic closure of  $\QQ_p$). 
These can be computed by the Newton polygon of this polynomial. 
Write $g = \sum_j \gamma_j T^{h-j}$. 
Plot the pairs $(j, v_p(\gamma_j))$ for $0 \leq j \leq h$. 
Consider the lower convex hull of $\{(j, v_p(\gamma_j)) \mid j\}$. 
This is a Newton polygon according to the definition above. 
The slopes of the sides of this polygon are precisely the 
$p$-adic values of the zeroes of $g$, ordered in non-decreasing order. 
(Suggestion: prove this as an exercise.) 

\vn
Later we will see: {\it a $p$-divisible group $X$ over a field 
of characteristic $p$ 
determines uniquely a Newton polygon.} 
In Section \ref{4} a correct and precise 
definition will be given; moreover we will see (Dieudonn\'e-Manin) 
that the isogeny class of a  $p$-divisible group over $k$ uniquely 
determines (and is uniquely determined by) its Newton polygon.

\vn  
(Incorrect.) 
Here we indicate what the Newton polygon of a $p$-divisible group is 
(in a slightly incorrect way ...). 
Consider ``the Frobenius endomorphism`` of $X$. 
This has a ``characteristic polynomial''. 
This polynomial determines a Newton polygon, which we write 
as $\cN(X)$, the Newton polygon of $X$. For an abelian variety $A$ 
we write $\cN(A)$ instead of $\cN(A[p^{\infty}])$.

Well, this ``definition'' is correct over $\FF_p$ as ground field. 
However over any other field $F: X \to X^{(p)}$ is not an endomorphism, 
and the above ``construction'' fails. 
Over a finite field there is a method 
which repairs this, see  \ref{pi}. 
However we need the Newton polygon of an abelian variety 
over an arbitrary field. Please accept for the time being the 
``explanation'' given above:  $\cN(X)$ is the 
``Newton polygon of the Frobenius on $X$'', 
which will be made precise later, see Section \ref{4}.\\
Some examples. The Newton polygon of $X = \GG_m[p^{\infty}]$ 
has one slope equal to $1$. 
In fact, on $\GG_m$ the Frobenius endomorphism is $\times p$.

The Newton polygon of $X = \underline{\QQ_p/\ZZ_p}$ 
 has one slope equal to $0$. 

The Newton polygon of an ordinary elliptic curve has two slopes, 
equal to $0$ and to $1$. 

The Newton polygon of a supersingular elliptic curve has two slopes, 
both equal to $1/2$.

\subsection{\bf Newton polygons go up under specialization.} 
In 1970 Grothendieck observed that`` Newton polygons go up''  
under specialization. In order to study this and related questions 
we introduce the notation of a {\it partial ordering} 
between Newton polygons.

\begin{picture}(150,50)(0,0)
\put(210,0){\line(6,1){60}}
\put(270,10){\line(1,1){30}}
\put(210,0){\line(4,1){40}}
\put(250,10){\line(5,3){50}}
\put(210,0){\line(1,0){120}}
\put(210,0){\line(0,1){40}}
\put(250,20){$\zeta_2$}
\put(280,10){$\zeta_1$}
\put(-15,30){We write $\zeta_1 \succ \zeta_2$ if  $\zeta_1$ is ``below'' 
$\zeta_2$,}

\put(-15,15){i.e. if no point of $\zeta_1$ is strictly above $\zeta_2$.}

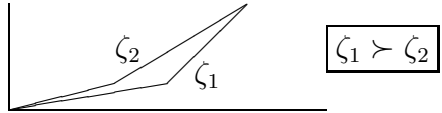
\put(330,20){\fbox{$\zeta_1 \succ \zeta_2$}}

\end{picture}

\vn
Note that we use this notation only if Newton polygons with the same 
endpoints are considered.

This notation may seem unnatural. However if  $\zeta_1$ is strictly below 
$\zeta_2$ the stratum defined by $\zeta_1$ is larger than the 
stratum defined by $\zeta_2$; this explains the choice for this notation.
 
\subsection{}\label{sym} Later we will show that 
isogenous $p$-divisible groups have the same Newton polygon. 
Using the construction defining a Newton polygon, see Section \ref{4}, 
and using \ref{dualth}, we will see that if $\cN(X)$ is given by 
$\{\beta_i \mid 1 \leq i \leq h\}$ then $\cN(X^t)$ is given by 
$\{1-\beta_h, \cdots , 1-\beta_1\}$.

\vn
A Newton polygon $\xi$, given by the slopes 
$\beta_1 \leq \cdots \leq \beta_h$ is called {\it symmetric} if 
$\beta_i = 1 - \beta_{h+1-i}$  for all $i$. 
We see that $X \sim X^t$ implies that $\cN(X)$ is symmetric; 
in particular for an abelian variety $A$ we see that $\cN(A)$ 
is symmetric. 
This was proved over finite fields by Manin, see \cite{Manin}, page 70; 
for any base field we can use the duality theorem over any base, 
see \cite{FO-CGS}, Th. 19.1, also see \ref{dualth}.

\subsection{}
If $S$ is a base scheme, $\cX \to S$ is a $p$-divisible 
group over $S$ and $\zeta$ is a Newton polygon  we write 
$$\cW_{\zeta}(S) := \{s \in S \mid \cN(\cX_s) \prec \zeta\} \subset S$$
and
$$\cW_{\zeta}^0(S) := \{s \in S \mid \cN(\cX_s) = \zeta\} \subset S.$$

\subsection{\bf Theorem} \Bl \ \  (Grothendieck and Katz; 
see \cite{Katz-Slope}, 2.3.2). 
$$\cW_{\zeta}(S) \subset S \quad\mbox{\it is a closed set.}
$$

\vn
Working over $S = \Spec(K)$, where $K$ is a perfect field, 
$\cW_{\zeta}(S)$ and $\cW_{\zeta}^0(S)$ will be given the 
induced reduced scheme structure.

As the set of Newton polygons of a given height is finite 
we conclude: 
$$
\cW_{\zeta}^0(S) \subset S \quad \mbox{\it  is a locally closed set}.
$$

\subsection{\bf Notation.} Let $\xi$ be a symmetric Newton polygon. 
We write $W_{\xi} = \cW_{\xi}(\cA_{g,1} \otimes \FF_p)$.

\subsection{} We have seen that ``Newton polygons go up under 
specialization''. Does a kind of converse hold? 
In 1970 Grothendieck conjectured the converse.
In \cite{G}, the appendix, we find a letter of Grothendieck to Barsotti, 
and on page 150 we read: {\it ``$\cdots$ The wishful 
conjecture I have 
in mind now is the following: the necessary conditions 
$\cdots$ that $G'$ 
be a specialization of $G$ are also sufficient. 
In other words, starting 
with a BT group $G_0 = G'$, taking its formal 
modular deformation $\cdots$ 
we want to know if every sequence of rational 
numbers satisfying $\cdots$ 
these numbers occur as the sequence of slopes 
of a fiber of $G$ as some 
point of $S$."}

\subsection{}\label{GC}{\bf Theorem} \Th \ {\rm ({\bf The 
Grothendieck Conjecture})} 
(conjectured by Grothendieck, Montreal 1970). {\it Let $K$ be a 
field of characteristic $p$, and let $X_0$ be a 
$p$-divisible group over $K$. We write 
${\cal N}({\cal X}_0) =: \beta$ for its Newton 
polygon. Suppose given a Newton polygon 
$\gamma$ ``below" $\beta$, i.e.  $\beta \prec \gamma$. 
There exists a deformation $X_{\eta}$ of $X_0$ 
such that ${\cal N}({\cal X}_{\eta}) = \gamma$.}\\
See \S 9.  This was proved by Frans Oort in 2001. 
For a proof see  \cite{AJdJ.FO}, 
\cite{FO-CH}, \cite{FO-NP}.

\vn
We say ``$X_{\eta}$ is a deformation  of $X_0$'' if there exists 
an integral scheme $S$ over $K$,  with generic point $\eta \in S$ 
and $0 \in S(K)$, and a $p$-divisible group $\cX \to S$ such that 
$\cX_0 = X_0$ and   $\cX_{\eta} = X_{\eta}$.

\vn
A (quasi-) polarized version will be given later.

\vn
In this paper we record  a proof of this theorem, 
and we will see that this is an important tool in understanding 
Newton polygon strata in $\cA_g$.

\vn
Why is the proof of this theorem difficult? A direct approach seems 
obvious: write down deformations of $X_0$,  compute Newton polygons of 
all fibers, and inspect whether all relevant Newton polygons 
appear in this way. 
However, computing the Newton polygon of a $p$-divisible group in general 
is difficult (but see Section \ref{5} how to circumvent this 
in an important special case). Moreover, abstract deformation theory 
is easy, but in general Newton polygon strata are ``very singular''; 
in Section \ref{7} we  describe how to ``move out'' of 
a singular point to a non-singular point of a Newton polygon stratum. 
Then, at non-singular points the deformation theory can be described 
more easily, see Section \ref{5}. By a combination of these two methods 
we achieve a proof of the Grothendieck conjecture. 
Later we will formulate and prove the analogous ``polarized case'' 
of the Grothendieck conjecture, 
see Section \ref{9}. 

\vn
We see: a direct approach did not work, but the detour via 
``deformation to $a \leq 1$''  plus the  results via 
Cayley-Hamilton gave the essential 
ingredients for a proof. Note the analogy of this method with 
the approach to liftability of abelian varieties to characteristic zero, 
as proposed by Mumford, and carried out in \cite{Norman.FO}.


\section{Serre-Tate theory}\label{2}

\noindent
In this section we explain the deformation theory of 
abelian varieties and Barsotti-Tate groups.
The content can be divided into several parts:
\smallbreak
\noindent (1) In \ref{subsec_fmdef} we give the formal definitions of
deformation functors for abelian varieties and Barsotti-Tate
groups.
\smallbreak

\noindent (2) In contrast to the deformation theory for general algebraic
varieties, the deformation theory for abelian varieties
and Barsotti-Tate groups can be efficiently dealt with
by linear algebra, as follows from the crystalline deformation
theory of Grothendieck-Messing. 
It says that, over an extension by an ideal with a divided
power structure, deforming abelian varieties or Barsotti-Tate
groups is the same as lifting the Hodge filtration.
See Thm.\ \ref{GM} for the precise statement, and 
Thm.\ref{cry-dual} for the behavior of the theory under duality.
The smoothness of the moduli space \(\calA_{g,1,n}\) follows
quickly from this theory.
\smallbreak

\noindent (3) The Serre-Tate theorem: deforming an abelian
variety is the same as deforming its Barsotti-Tate group.
See Thm.\ \ref{thm:ST} for a precise statement.  
A consequence is that the deformation space
of a polarized abelian variety admits an natural action by
a large \(p\)-adic group, see \ref{cor:action1}, \ref{cor:action2}.
In general this action is poorly understood.
\smallbreak

\noindent (4) There is one case when the action on the deformation
space mentioned in (3) above is linearized and well-understood.
This is the case when the abelian variety is ordinary.
The theory of Serre-Tate coordinates says that the deformation space
of an ordinary abelian variety has a natural structure as a 
formal torus.  See Thm.\ \ref{STcoord} for the statement.
In this case the action on the local moduli space mentioned in (3)
above preserve the group structure and gives a linear representation
on the character group of the Serre-Tate formal torus.
This phenomenon has important consequences later. 
A local rigidity result Thm.\ \ref{localrig}  
is important for the Hecke orbit problem in that it provides
an effective linearization of the Hecke orbit problem.
Also, computing the deformation using the Serre-Tate coordinates is 
often easy.  The reader is encouraged to try Exer.\ \ref{compdef}
as an example of the last sentence.

\noindent 
Barsotti-Tate groups: \cite{Messing-thesis}, \cite{Ill}.\\
Crystalline deformation theory: \cite{Messing-thesis}, \cite{BBM}.\\
Serre-Tate Theorem: \cite{Messing-thesis}, \cite{Katz-STthm}.\\
Serre-Tate coordinates: \cite{Katz-ST}.\\

\subsection{\scshape Deformation of abelian varieties and BT-groups: 
definitions}\label{subsec_fmdef}
\medbreak

\noindent{\bf Definition.}  Let \(K\) be a perfect field of
characteristic \(p\). Denote by \(W(K)\) the ring of \(p\)-adic 
Witt vectors with coordinates in \(K\).\\
(i) Denote by \(\,{\rm{Art}}_{W(K)}\,\) the category 
of Artinian local algebras over \(W(K)\).
 An objects of \({\rm{Art}}_{W(K)}\)
is a pair \((R,j)\), where 
\(R\) is an Artinian local algebra and \(\epsilon:W(K)\to R\)
is an local homomorphism of local rings. 
A morphism in \({\rm{Art}}_{W(K)}\) from 
\((R_1,j_1)\) to \((R_2,j_2)\) is a homomorphism
\(h:R_1\to R_2\) between Artinian local rings such that
\(h\circ j_1=j_2\).
\smallbreak

\noindent (ii) Denote by \(\,{\rm{Art}}_K\,\) 
the category of Artinian local 
\(K\)-algebras.  An object in \(\,{\rm{Art}}_K\,\) is a pair
\((R,j)\), where \(R\) is an Artinian local algebra and
\(\epsilon:K\to R\,\) is a ring homomorphism.
A morphism in \(\,{\rm{Art}}/K\,\) from 
\((R_1,j_1)\) to \((R_2,j_2)\) is a homomorphism
\(h:R_1\to R_2\) between Artinian local rings such that
\(h\circ j_1=j_2\).
Notice that \(\,{\rm{Art}}_K\,\) is a fully faithful subcategory 
of \(\,{\rm{Art}}_{W(K)}\,\).
\medbreak

\noindent{\bf Definition.} Denote by 
\[\rm{SETS}\] the
category whose objects are sets and whose morphisms are
isomorphisms of sets. 
\medbreak

\noindent{\bf Definition.} Let \(A_0\) be an abelian variety
over a perfect field \(K\supset \Fp\).
The deformation functor of \(A_0\)
is a functor 
\[\,{\rm{Def}}(A_0/W(K)): 
{\rm{Art}}_{W(K)}\to {\rm{SETS}}\,\]
defined as follows.
For every object \((R,\epsilon)\) of \({\rm{Art}}_{W(K)}\),
\(\,{\rm{Def}}(A_0/W(K))(R,\epsilon)\,\) is
the set of isomorphism classes of pairs
\((A\to \Spec(R), \epsilon)\), where
\(A\to \Spec(R)\) is an abelian scheme, and 
\[\epsilon:A\times_{\Spec(R)}\Spec(R/\grm_R)\xrightarrow{\sim}
A_0\times_{\Spec(K)}\Spec(R/\grm_R)\]
is an isomorphism of abelian varieties over \(R/\grm_R\).
Denote by \(\,{\rm{Def}}(A_0/K)\) the restriction of 
\(\,{\rm{Def}}(A_0/W(K))\,\) to the faithful subcategory
\(\,{\rm{Art}}_K\,\) of \(\,{\rm{Art}}_{W(K)}\,\).
\bigbreak

\noindent{\bf Definition.} Let \(A_0\) be an abelian variety
over a perfect field \(K\supset \Fp\), and let \(\lambda_0\) be
a polarization on \(A_0\).
The deformation functor of \((A_0,\lambda_0)\)
is a functor 
\[\,{\rm{Def}}(A_0/W(K)): 
{\rm{Art}}_{W(K)}\to {\rm{SETS}}\,\]
defined as follows.
For every object \((R,\epsilon)\) of \(\,{\rm{Art}}_{W(K)}\,\),
\(\,{\rm{Def}}(A_0/W(K))(R,\epsilon)\,\)
is the set of isomorphism classes of pairs
\((A,\lambda)\to \Spec(R), \epsilon)\), where
\((A,\lambda)\to \Spec(R)\) is a polarized abelian scheme, and 
\[\epsilon:(A,\lambda)\times_{\Spec(R)}\Spec(R/\grm_R)\xrightarrow{\sim}
(A_0,\lambda_0)\times_{\Spec(K)}\Spec(R/\grm_R)\]
is an isomorphism of polarized abelian varieties over \(R/\grm_R\).
Denote by \(\,{\rm Def}((A_0,\lambda_0)/K)\) the restriction of 
\(\,{{\rm Def}}(A_0/W(K))\,\) to the faithful subcategory
\(\,{\rm Art}_K\,\) of \(\,{\rm Art}_{W(K)}\,\).
\bigbreak

\noindent{\bf Definition.} Let \(X_0\) be a Barsotti-Tate group
over a perfect field \(K\supset \Fp\).
Let 
\[\,{\rm{Def}}(X_0/W(K)): 
{\rm{Art}}_{W(k)}\to {\rm{SETS}}\,\]
be the deformation functor of \(X_0\), defined as follows.
For every object \((R,\epsilon)\) of \(\,{\rm{Art}}_{W(K)}\,\),
\(\,{\rm{Def}}(A_0/W(K))(R,\epsilon)\,\) is
the set of isomorphism classes of pairs
\((X\to \Spec(R), \epsilon)\), where
\(X\to \Spec(R)\) is a BT-group over \(\Spec(R)\),
and 
\[\epsilon:X\times_{\Spec(R)}\Spec(R/\grm_R)\xrightarrow{\sim}
X_0\times_{\Spec(K)}\Spec(R/\grm_R)\]
is an isomorphism of BT-groups over \(R/\grm_R\).
Denote by \(\,{\rm{Def}}(X_0/K)\) the restriction of the functor
\(\,{\rm{Def}}(X_0/W(K))\,\) to the faithful subcategory
\(\,{\rm{Art}}_K\,\) of \(\,{\rm{Art}}_{W(K)}\,\).
\smallbreak

Suppose that \(\lambda_0\) is a polarization on \(X_0\);
the deformation functor of \((X_0,\lambda_0)\) is the functor
\(\,{\rm{Def}}((X_0,\lambda_0)/W(K))\,\) on
\(\,{\rm{Art}}_{W(K)}\,\)
which sends an object \((R,\epsilon)\) of \(\,{\rm{Art}}_{W(K)}\,\)
to the set of isomorphism classes of pairs
\(((X,\lambda)\to \Spec(R), \epsilon)\), where
\((X,\lambda)\to \Spec(R)\) is a polarized BT-group over \(\Spec(R)\),
and 
\[\epsilon:(X,\lambda)\times_{\Spec(R)}\Spec(R/\grm_R)\xrightarrow{\sim}
(X_0,\lambda)\times_{\Spec(k)}\Spec(R/\grm_R)\]
is an isomorphism of polarized BT-groups over \(R/\grm_R\).
Denote by \(\,{\rm{Def}}((X_0,\lambda_0)/K)\,\) the
restriction of \(\,{\rm{Def}}((X_0,\lambda_0)/W(K))\,\)
to \(\,{\rm{Art}}_{K}\,\).
\bigbreak

\begin{defn} Let \(R\) be a commutative ring, and
let \(I\subset R\) be an ideal of \(I\).
A \emph{divided power structure} (a DP structure for short) on \(I\)
is a collection of 
maps \(\,\gamma_i:I\to R\), \(i\in \bbN\),  such that
\begin{itemize}
\item \(\gamma_0(x)=1\) \ \ \(\forall x\in I\),
\item \(\gamma_1(x)=x\) \ \ \(\forall x\in I\),
\item \(\gamma_i(x)\in I\) \ \ 
\(\forall x\in I\), \(\forall i\geq 1\),
\item \(\gamma_j(x+y)=\sum_{0\leq i\leq j} \gamma_i(x)\gamma_{j-i}(y)\)
\ \ \(\forall x, y\in I\), \(\forall j\geq 0\),

\item \(\gamma_i(ax)=a^i\) \ \ \(\forall a\in R\), \(\forall x\in I\),
\(\forall i\geq 1\),

\item \(\gamma_i(x)\gamma_j(y)=\frac{(i+j)!}{i!j!} \gamma_{i+j}(x)\)
\ \ \(\forall i,j\geq 0\), \(\forall x\in I\),

\item \(\gamma_i(\gamma_j(x))=\frac{(ij)!}{i! (j!)^i} \gamma_{ij}(x)\)
\ \ \(\forall i, j\geq 1\), \(\forall x\in I\).
\end{itemize}
A divided power structure \((R, I, (\gamma_i)_{i\in \bbN})\) as above
is \emph{locally nilpotent} if there exist \(n_0\in \bbN\) such that
\(\gamma_i(x)=0\) for all \(i\geq n_0\) and all \(x\in I\).
A \emph{locally nilpotent} DP \emph{extension} 
of a commutative ring \(R_0\)
is a locally nilpotent DP structure \((R, I, (\gamma_i)_{i\in \bbN})\)
together with an isomorphism \(R/I\xrightarrow{\sim} R_0\).
\end{defn}
\bigbreak

\begin{rem}\enspace\label{trivDP}
Let \(R\) be a commutative ring with \(1\), and let \(I\) be an
ideal of \(R\) such that \(I^2=(0)\). 
Define a DP structure on \(I\) by requiring that
\(\,\gamma_i(x)=0\) for all \(i\geq 2\) and all \(x\in I\).
This DP structure on a square-zero ideal \(I\) will be called the
\emph{trivial DP structure} on \(I\).
An extension of a ring \(R_0\) by a square-zero ideal \(I\) is
forms a standard ``input data'' in deformation theory.
So we can feed such input data into the crystalline deformation
theory summarized in Thm.\ \ref{GM} below to translate the deformation
of an abelian scheme \(A\to \Spec(R_0)\)
over a square-zero extension \(R\twoheadrightarrow R_0\)
into a question about lifting Hodge filtrations.
\end{rem}
\medbreak

\begin{thm}\Bl\  {\rm (Grothendieck-Messing).} \label{GM}
Let \(X_0\to \Spec(R_0)\) be an Barsotti-Tate group
over a commutative ring \(R_0\).
\begin{itemize}
\rmitem[(i)] To every locally nilpotent {\rm DP} extension 
\((R, I, (\gamma_i)_{i\in \bbN})\) of \(R_0\) there is a
functorially attached locally free \(R\)-module \(\bbD(X_0)_R\)
of rank \({\rm ht}(X_0)\).  The functor \(\bbD({X_0})\)
is called the covariant Dieudonn\'e crystal attached to \(X_0\).

\rmitem[(ii)] Let \((R, I, (\gamma_i)_{i\in \bbN})\) be a locally nilpotent
{\rm DP} extension of \(R_0\).  Suppose that \(X\to \Spec(R)\) is
an Barsotti-Tate group extending \(X_0\to \Spec(R_0)\).
Then there is a functorial short exact sequence
\[
0\to {\rm Lie}(X^t/R)^{\vee} \to
\bbD(X_0)_R
\to {\rm Lie}(X/R)\to 0\,.
\]
Here \({\rm Lie}(X/R)\) is the tangent space of the
BT-group \(X\to\Spec(R)\),
which is a projective \(R\)-module of rank \(\dim(X/R)\), and
\({\rm Lie}(X^t/R)^{\vee}\) is the \(R\)-dual of the 
tangent space of the Serre dual \(X^t\to\Spec(R)\) of \(X\to \Spec(R)\).

\rmitem[(iii)] Let \((R, I, (\gamma_i)_{i\in \bbN})\) be a locally nilpotent
{\rm DP} extension of \(R_0\).
Suppose that \(A\to \Spec(R)\) is an abelian scheme such that
there exist an isomorphism
\[\beta:A[p^{\infty}]\times_{\Spec(R)}\Spec(R_0)\xrightarrow{\sim} X_0\] 
of BT-groups over \(R_0\).
Then there exists a natural isomorphism
\[
\bbD(X_0)_R\to {\rm H}_1^{\rm DR}(A/R)\,,
\]
where \({\rm H}_1^{\rm DR}(A/R)\) is the first De Rham homology of
\(A\to \Spec(R)\).
Moreover the above isomorphism identifies the short exact sequence
\[0\to {\rm Lie}(A[p^{\infty}]^t/R)^{\vee} \to
\bbD(X_0)_R
\to {\rm Lie}(A[p^{\infty}]/R)\to 0\]
described in {\rm (ii)} 
with the Hodge filtration
\[
0\to {\rm Lie}(A^t/R)^{\vee} \to
{\rm H}_1^{\rm DR}(A/R)
\to {\rm Lie}(A/R)\to 0
\]
on \({\rm H}_1^{\rm DR}(A/R)\).

\rmitem[(iv)] Let \((R, I, (\gamma_i)_{i\in \bbN})\) be a
locally  nilpotent
{\rm DP} extension of \(R_0\).
Denote by \(\,\mathfrak{E}\,\) the category whose objects are
short exact sequences 
\(
0\to F\to \bbD(X_0)_R\to Q\to 0
\)
such that \(F\) and \(Q\) are projective \(R\)-modules,
plus an isomorphism
from the short exact sequence
\(\,(0\to F\to \bbD(X_0)_R\to Q\to 0)\otimes_R R_0\,\)
of projective \(R_0\)-modules
to the short exact sequence
\[0\to {\rm Lie}(X_0^t)^{\vee}\to \bbD(X_0)_{R_0}\to 
{\rm Lie}(X)\to 0
\]
attached to the BT-group \(X_0\to \Spec(R_0)\)
as a special case of {\rm (ii)} above.
The morphisms in  \(\,\mathfrak{E}\,\) are maps between
diagrams.
Then the functor from the category of BT-groups over
\(R\) lifting \(X_0\) to the category \(\,\mathfrak{E}\,\)
described in {\rm (ii)}  is an equivalence of categories.
\end{itemize}
\end{thm}
\bigbreak

\begin{cor} \label{cor:defbt}
Let \(X_0\) be a BT-group
over a perfect field \(K\supset \Fp\).
Let \(d=\dim(X_0)\), \(c=\dim(X_0^t)\).
The deformation functor \({\rm{Def}}(X_0/W(K))\) of \(X_0\)
is representable by a smooth formal scheme over \(W(K)\)
of dimensional \(cd\). 
In other words,  \({\rm{Def}}(X_0/W(K))\) is non-canonically isomorphic
to \(\,\Spf\left(W(K)[[x_1,\ldots,x_{cd}]]\right)\).
\end{cor}

\proof
Apply Thm.\ \ref{GM} to the trivial DP structure on pairs \((R,I)\)
with \(I^2=(0)\), we see that \(\,{\rm{Def}}(X_0/W(K))\) is
formally smooth over \(W(K)\).  Apply Thm.\ \ref{GM} again to the pair
\((K[t]/(t^2), tK[t]/(t^2)\), we see that the dimension of the
tangent space of \(\,{\rm{Def}}(X_0/K)\,\) is equal to \(cd\).
\hfill\qed

\noindent{\bf Remark.} From the definition of \({\rm{Def}}(A_0/W(K))\),
there is a natural action of \({\rm Aut}(A_0)\) on the 
smooth formal scheme \({\rm{Def}}(A_0/W(K))\cong
\Spf(W(K)[[x_1,\ldots,x_{g^2}]])\).
\bigbreak

\subsection
We set up notation for the Serre-Tate theorem \ref{thm:ST}.
Let \(p\) be a prime number.
Let \(S\) be a scheme such that \(p\) is locally nilpotent in 
\(\ringO_S\).
Let \(I\subset \ringO_S\) be a coherent sheaf of ideals such that
\(I\) is locally nilpotent.
Let \(S_0=\underline{\rm Spec}(\ringO_S/I)\).
Denote by \(\,{\rm{AV}}_S\,\) the category of abelian schemes over 
\(S\).  
Denote by \(\,{\rm{AVBT}}_{S_0,S}\,\) 
the category whose objects
are triples \((A_0\to S_0,X\to S,\epsilon)\), where
\(A_0\to S_0\) is an abelian scheme over \(S_0\), 
\(X\to S\) is a Barsotti-Tate group over \(S\),
and \(\epsilon: X\times_S S_0\to A_0[p^{\infty}]\) is an
isomorphism of BT-groups. A morphism from 
\((A_0\to S_0,X\to S,\epsilon)\) to \((A_0'\to S_0,X'\to S,\epsilon')\)
is a pair
\((h,f)\), where \(h_0:A_0\to A_0'\) is a homomorphism
of abelian schemes over \(S_0\), \(f:X\to X'\) is a homomorphism of
BT-groups over \(S\), such that
\(h[p^{\infty}]\circ \epsilon
= \epsilon' \circ(f\times_S S_0)\).
Let \[\,\FF: {\rm{AV}}_S \to {\rm{AVBT}}_{S_0,S}\,\] 
be the functor which sends an abelian scheme 
\(A\to S\) to the triple \((A\times_S S_0, A[p^{\infty}], {\rm can})\)
where \({\rm can}\) is the canonical isomorphism
\(A[p^{\infty}]\times_S S_0\xrightarrow{\sim}(A\times_S S_0)[p^{\infty}]\).
\smallbreak

\begin{thm}\Bl\ {\rm (Serre-Tate).} \label{thm:ST}
Notation and assumptions as in the above paragraph. 
The functor \(\FF\) is an equivalence
of categories.
\end{thm}
\bigbreak

\noindent{\bf Remark.}\enspace See \cite{WoodsHole}. 
A proof of Thm.\ \ref{thm:ST} first
appeared in print in \cite{Messing-thesis}.  See also \cite{Katz-STthm}.
\medbreak

\begin{cor} \label{cor:ST1}
Let \(A_0\) be an variety over a perfect field \(K\).
Let \[\GG:{\rm{Def}}(A_0/W(K))\to 
{\rm{Def}}(A_0[p^{\infty}]/W(K))\]
be the functor which sends any object
\[
\left(A\to \Spec(R), \epsilon:A\times_{\Spec(R)}\Spec(R/\grm_R)
\xrightarrow{\sim} A_0\times_{\Spec(k)}\Spec(R/\grm_R)
\right)
\]
in \({\rm{Def}}(A_0/W(K))\) to the object
\[
\left(A[p^{\infty}]\to \Spec(R), 
\epsilon[p^{\infty}]:A[p^{\infty}]\times_{\Spec(R)}\Spec(R/\grm_R)
\xrightarrow{\sim} A_0[p^{\infty}]\times_{\Spec(K)}\Spec(R/\grm_R)
\right)
\]
in \({\rm{Def}}(A_0[p^{\infty}]/W(K))\).
The functor \(\GG\) is an equivalence of categories.
\end{cor}
\medbreak

\noindent{\bf Remark.} In words, Cor.\ \ref{cor:ST1} says that deforming
an abelian variety is the same as deforming its 
\(p\)-divisible group.
\bigbreak

\begin{cor}\label{cor:g2}
Let \(A_0\) be a \(g\)-dimensional abelian variety
over a perfect field \(K\supset \Fp\).
The deformation functor \({\rm{Def}}(A_0/W(K))\) of \(A_0\)
is representable by a smooth formal scheme over \(W(K)\)
of relative dimension \(g^2\). 
\end{cor}

\proof We have  
\(\,{\rm{Def}}(A_0/W(K))\cong {\rm{Def}}(A_0[p^{\infty}]/W(K))\,\) 
by Thm.\ \ref{thm:ST}.
Cor.\ \ref{cor:g2} follows from Cor.\ \ref{cor:defbt}.
\hfill\qed
\medbreak

\subsection
Let \(R_0\) be a commutative ring.
Let \(A_0\to \Spec(R_0)\) be an abelian scheme.
Let \(\bbD(A_0):=\bbD(A_0[p^{\infty}])\) be the covariant 
Dieudonn\'e crystal attached to \(A_0\).
Let \(\bbD(A_0^t)\) be the covariant Dieudonn\'e crystal attached
to the dual abelian scheme \(A^t\).
Let \(\bbD(A_0)^{\vee}\) be the dual of \(\bbD(A_0)\), i.e.
\[
\bbD(A_0)^{\vee}_{R}={\rm Hom}_R(\bbD(A_0)_R,R)
\]
for any DP extension \((R, I, (\gamma_i)_{i\in\bbN})\) 
of \(R_0=R/I\).

\begin{thm}\label{cry-dual}
We have functorial isomorphisms 
\[
\phi_{A_0}: \bbD(A_0)^{\vee}\xrightarrow{\sim} \bbD(A_0^t)
\]
for abelian varieties \(A_0\) over \(K\)
with the following properties.
\begin{itemize}
\rmitem[(1)] The composition
\[
\bbD(A_0^t)^{\vee}\xrightarrow[\sim]{\phi_{A_0}^{\vee}} 
(\bbD(A_0)^{\vee})^{\vee}
=\bbD(A_0)\xrightarrow[\sim]{j_{A_0}}
\bbD((A_0^t)^t)
\]
is equal to 
\[
-\phi_{A_0^t}: \bbD(A_0^t)^{\vee}\xrightarrow{\sim}
\bbD((A_0^t)^t)\,,
\]
where the isomorphism \(\,\bbD_{A_0}\xrightarrow[\sim]{j_{A_0}}
\bbD((A_0^t)^t)\,\) is induced by the
canonical isomorphism 
\[\,A_0\xrightarrow{\sim}(A_0^t)^t\,.\]

\rmitem[(2)] For any locally nilpotent
DP extension \((R, I, (\gamma_i)_{i\in\bbN})\) 
of \(R_0=R/I\), and any lifting \(A\to \Spec(R)\) of \(A_0\to\Spec(R_0)\)
to \(R\), 
the following diagram 
\[
\xymatrix{
0\ar[r] &\Lie(A/R)^{\vee}\ar[d]^{\cong}\ar[r]
&\bbD(A_0)_R^{\vee}\ar[d]^{\phi_{A_0}}\ar[r]
&(\Lie(A^t/R)^{\vee})^{\vee}\ar[d]^{=}\ar[r]&0
\\
0\ar[r]&\Lie((A^t)^t/R)^{\vee}\ar[r]
&\bbD(A_0^t)_R\ar[r]
&\Lie(A^t/R)\ar[r]&0
}
\]
commutes. Here the bottom horizontal exact sequence is as in \ref{GM},
the top horizontal sequence is the dual of the short exact 
sequence in \ref{GM}, 
and the left vertical isomorphism is induced by the canonical
isomorphism \(\,A\xrightarrow{\sim} (A^t)^t\).
\end{itemize}
\end{thm}
\medbreak

\begin{cor}
Let \((A_0,\lambda_0)\) be a \(g\)-dimensional principally polarized 
abelian variety over a perfect field \(K\supset\Fp\).
The deformation functor \({\rm{Def}}((A_0,\lambda)/W(K))\) 
of \(A_0\) is representable by a smooth formal scheme over \(W(K)\)
of dimensional \(g(g+1)/2\).
This statement can be reformulated as follows.
Let \(\eta_0\) be a \(K\)-rational symplectic level-\(n\) structure on \(A_0\),
\(n\geq 3\), \((n,p)=1\), and let \(x_0=[(A_0,\lambda_0,\eta_0)]\in
\calA_{g,1,n}(K)\).
The formal completion \(\calA_{g,1,n}^{/x_0}\) of 
the moduli space \(\calA_{g,1,n}^{/x_0}\to \Spec(W(K))\)
is non-canonically isomorphic to 
\(\,\Spf\left(W(K)[[x_1,\ldots,x_{g(g+1)/2}]]\right)\).
\end{cor}
\medbreak

The proof follows quickly from Thm.\ \ref{GM} and Thm.\ \ref{cry-dual},
and is left as an exercise.
\medbreak

\begin{cor}
Let \((A_0,\lambda_0)\) be a polarized abelian variety over
a perfect field \(K\supset \Fpbar\); let \({\rm deg}(\lambda_0)=d^2\).
\begin{itemize}
\rmitem[(i)] The natural map \(\,{\rm{Def}}((A_0,\lambda_0)/W(K))
\to {\rm{Def}}(A_0/W(K))\,\) is represented by a closed
embedding of formal schemes.

\rmitem[(ii)] Let \(n\) be a positive integer, \(n\geq 3\), \((n,pd)=1\).
Let \(\eta_0\) be a \(K\)-rational
symplectic level-\(n\) structure on \((A_0,\lambda_0)\).
Let \(x_0=[(A_0,\lambda_0,\eta_0)]\in\calA_{g,d,n}(K)\).
The formal completion \(\calA_{g,d,n}^{/x_0}\) of 
the moduli space \(\calA_{g,d,n}\to \Spec(W(K))\) at the closed
point \(x_0\) is isomorphic to \(\,{\rm{Def}}((A_0,\lambda_0)/W(K))\).
\end{itemize}
\end{cor}
\medbreak

\begin{cor} \label{cor:action1}
Let \(A_0\) be an variety over a perfect field \(K\supset \Fpbar\). 
There is a natural action of
the profinite group \({\rm Aut}(A_0)\) on the
smooth formal scheme \(\,{\rm{Def}}(A_0/W(K))\).
\end{cor}
\medbreak

\begin{cor} \label{cor:action2}
Let \((A_0,\lambda_0)\) be a polarized abelian variety 
over a perfect field \(K\supset \Fpbar\,\).
Let \({\rm Aut}((A_0,\lambda_0)[p^{\infty}])\)
be the closed subgroup of \({\rm Aut}(A_0[p^{\infty}])\)
consisting of all automorphisms of \({\rm Aut}(A[p^{\infty}])\)
compatible with the quasi-polarization \(\lambda_0[p^{\infty}]\).
The natural action in {\rm Cor.\ \ref{cor:action1}}
induces a natural action of 
on the closed formal subscheme
\(\,{\rm{Def}}(A_0,\lambda_0)\,\)
of \(\,{\rm{Def}}(A_0)\).
\end{cor}
\medbreak

\subsection{\scshape\'Etale and multiplicative BT-groups: Notation}
\smallbreak

\noindent{\bf Remark.}\enspace  Let \(E\to S\) be an 
\'etale Barsotti-Tate group,
where \(S\) is a scheme.  The \(p\)-adic Tate module of \(E\), 
defined by 
\[{\rm T}_p(E):=\varprojlim_{n} E[p^n]\,,\]
is representable by a smooth \(\Zp\)-sheaf on \(S_{\rm et}\)
whose rank is equal to \({\rm ht}(E/S)\).
When \(S\) is the spectrum of a field \(K\), \({\rm T}_p(E)\)
``is'' a free \(\Zp\)-module with an action by \({\rm Gal}(K^{\rm sep}/K)\); 
see \ref{et}.\\
{\bf Warning.} The notation ${\rm T}_p(-)$ could lead to confusion 
when used for arbitrary $p$-divisible groups: 
the \'etale part need not be of constant height for a $p$-divisible 
group over a base.
\bigbreak

\noindent{\bf Remark.}\enspace Attached to any multiplicative 
Barsotti-Tate group \(T\to S\) is its character group
\(X^*(T):=\underline{\rm Hom}(T,\bbG_m[p^{\infty}])\)
and cocharacter group
\(X_*(T):=\underline{\rm Hom}(\bbG_m[p^{\infty}],T)\).
The character group of \(T\) can be identified with 
the \(p\)-adic Tate module of the Serre-dual \(T^t\) of \(T\),
and \(T^t\) is an \'etale Barsotti-Tate group over \(S\).
Both \(X^*(T)\) and \(X_*(T)\) are smooth \(\Zp\)-sheaves of
rank \(\dim(T/S)\) on \(S_{\rm et}\), and they are naturally dual to
each other. 
\medbreak

\begin{defn} Let \(S\) be either a scheme such that
\(p\) is locally nilpotent in \(\ringO_S\), or an adic
formal scheme such that \(p\) is locally topologically nilpotent in 
\(\ringO_S\). An Barsotti-Tate group 
\(X\to S\) is \emph{ordinary} if \(X\) sits in the middle
of a short exact sequence
\[
0\to T\to X\to E\to 0
\]
where \(T\) (resp.\ \(E\)) is a multiplicative (resp.\ \'etale) 
Barsotti-Tate group. 
Such an exact sequence is unique up to unique isomorphisms.
\end{defn}
\medbreak

\noindent{\bf Remark.}\enspace 
If \(S=\Spec(K)\), \(K\) is a perfect field
of characteristic \(p\), and \(X\) is an ordinary Barsotti-Tate group
over \(K\),
 there exists a unique splitting
of the short exact sequence \(0\to T\to X\to E\to 0\) over \(K\).
\bigbreak

\begin{prop}\Bl \label{prop:ordinary}
Suppose that \(S\) is a scheme over \(W(K)\)
and \(p\) is locally nilpotent in \(\ringO_S\). 
Let \(S_0=\underline{\rm Spec}(\ringO_S/p\ringO_S)\), the
closed subscheme of \(S\) defined by the ideal \(p\ringO_S\)
of the structure sheaf \(\ringO_S\).
If \(X\to S\) is a BT-group such that \(X\times_S S_0\) is
ordinary, then \(X\to S\) is ordinary.
\end{prop}
\bigbreak

\subsection 
We set up notation for Thm.\ \ref{STcoord} on
the theory of Serre-Tate local coordinates.
Let \(K\supset \Fp\) be a perfect field
and let \(X_0\) be an ordinary Barsotti-Tate
group over \(K\). This means that there is a natural split
short exact sequence
\[
0\to T_0\to X_0\to E_0\to 0\,
\]
where \(T_0\) (resp.\ \(E_0\)) is a multiplicative (resp.\ \'etale)
Barsotti-Tate group over \(K\).
Let \(T_i\to \Spec(W(K)/p^iW(K))\) (resp.\ \(E_i\to \Spec(W(K)/p^iW(K))\))
be the multiplicative (resp.\ \'etale) BT-group over \(\Spec(W(K)/p^iW(K))\)
which lifts \(T_0\) (resp.\ \(E_0\)) for each \(i\geq 1\).
Both \(T_i\) and \(E_i\) are unique up to unique isomorphism.
Taking the limit of \(T_i[p^n]\) (resp.\ \(E_i[p^i]\)) as 
\(i\to \infty\), we get a multiplicative (resp.\ \'etale) 
BT$_n$-group \(\tilde{T}\to \Spec(W(K))\) 
(resp.\ \(\tilde{E}\to \Spec(W(K)))\) over \(W(K)\).
Denote by \(\hat{T}\) the formal torus over \(W(K)\) attached to
\(T_0\); i.e. \(\hat{T}=X_*(T_0)\otimes_{\Zp}\widehat{\bbG_m}\),
where \(\widehat{\bbG_m}\) is the formal completion of 
\(\bbG_m\to \Spec(W(K))\) along its unit section.
\smallbreak

\begin{thm}\label{STcoord} Notation and assumption as above.
\begin{itemize}
\rmitem[(i)] Every deformation \(X\to \Spec(R)\) of \(X_0\) over an Artinian
local \(W(K)\)-algebra \(R\) is an ordinary BT-group over \(R\).
Therefore \(X\) sits in the middle of a short exact sequence
\[
0\to \tilde{T}\times_{\Spec(W(K))}\Spec(R)
\to X \to \tilde{E}\times_{\Spec(W(K))}\Spec(R)\to 0\,.
\]

\rmitem[(ii)] The deformation functor \(\,{\rm{Def}}(X_0/W(K)\,\)
has a natural structure, via the Baer sum construction, as
a functor from \(\,{\rm{Art}}_{W(K)}\,\) to the
category \(\,{\rm{AbG}}\,\) of abelian groups.
In particular the unit element in \(\,{\rm{Def}}(X_0/W(K)(R)\) 
corresponds to the BT-group 
\[\left(\tilde{T}\times_{\Spec(W(K))}\tilde{E}\right)
\times_{\Spec(W(K))}\Spec(R)\]
over \(R\).

\rmitem[(iii)] There is a natural isomorphism of functors
\begin{align*}
{\rm{Def}}(X_0/W(K))
\xleftarrow{\sim} \underline{\rm Hom}_{\Zp}({\rm T}_p(E_0),\hat{T})
&={\rm T}_p(E)^{\vee}\otimes_{\Zp}X_*(T_0)\otimes_{\Zp}\widehat{\bbG_m}\\
&=\underline{\rm Hom}_{\Zp}\left({\rm T}_p(E_0)\otimes_{\Zp}
X^*(T_0), \widehat{\bbG_m}\right)
\,.
\end{align*}
In other words, the deformation space \(\,{\rm{Def}}(X_0/W(K))\)
of \(X_0\) has a natural structure as a formal torus over \(W(K)\)
whose cocharacter group is isomorphic to the
\({\rm Gal}(K^{\rm alg}/K)\)-module
\({\rm T}_p(E)^{\vee}\otimes_{\Zp}X_*(T_0)\).
\end{itemize}
\end{thm}

\proof The statement (i) is follows from Prop.\ \ref{prop:ordinary},
so is (ii).  It remains to prove (iii).
\smallbreak

By \'etale descent, we may and do assume that \(K\) is algebraically
closed. By (i), over any Artinian local \(W(K)\)-algebra \(R\),
\({\rm Def}(X_0/W(K))(R)\) is the set of isomorphism classes
of extensions of \(\tilde{E}\times_{W(K)}\Spec(R)\) by
\(\tilde{T}\times_{W(K)}\Spec(R)\). 
Write \(T_0\) (resp.\ \(E_0\)) as a product of a finite number of copies of 
\(\bbG_m[p^{\infty}]\) (resp.\ \(\Qp/\Zp\)), we only need to verify
the statement (iii) in the case when \(T_0=\bbG_m[p^{\infty}]\)
and \(E_0=\Qp/\Zp\).
\smallbreak

Let \(R\) be an Artinian local \(W(k)\)-algebra.
We have seen that \(\,{\rm{Def}}(\Qp/\Zp,\bbG_m[p^{\infty}])(R)\,\)
is naturally isomorphic to the inverse limit 
\(\varprojlim_n\,{\rm Ext}_{\Spec(R)}^1(p^{-n}\bbZ/\bbZ,\mu_{p^n})\).
By Kummer theory, we have
\[
{\rm Ext}_{\Spec(R)}^1(p^{-n}\bbZ/\bbZ,\mu_{p^n})
=R^{\times}/(R^{\times})^{p^n}
= (1+\grm_R)/(1+\grm_R)^{p^n}\,;
\]
the second equality follows from the hypothesis that \(K\) is perfect.
We know that \(p\in\grm_R\) and \(\grm_R\) is nilpotent.
Hence there exists an \(n_0\) such that 
\((1+\grm_R)^{p^n}=1\) for all \(n\geq n_0\).
Taking the inverse limit as \(n\to\infty\), we see that
the natural map 
\[1+\grm_R\to \varprojlim_n\,
{\rm Ext}_{\Spec(R)}^1(p^{-n}\bbZ/\bbZ,\mu_{p^n})
\]
is an isomorphism.
\hfill\qed
\bigbreak

\begin{cor}\label{STcoordav}
Let \(K\supset \Fp\) be a perfect field, and let \(A_0\)
be an ordinary abelian variety.
Let \({\rm T}_p(A_0):={\rm T}_p(A_0[p^{\infty}]_{\rm et})\),
\({\rm T}_p(A_0^t):={\rm T}_p(A_0^t[p^{\infty}]_{\rm et})\).
Then 
\[\,{\rm Def}(A_0/W(K))
\cong \underline{\rm Hom}_{\Zp}({\rm T}_p(A_0)\otimes_{\Zp}
{\rm T}_p(A_0^t),\widehat{\bbG}_m)\,.
\]
\end{cor}
\bigbreak

\subsection{\bf Exercise}
Let \(R\) be a commutative ring with \(1\).
Compute
$$
\Ext^1_{\Spec(R), (\bbZ/n\bbZ)}(n^{-1}\bbZ/\bbZ,\mu_n),
$$
the group of isomorphism classes of extensions of the constant
group scheme \(n^{-1}\bbZ/\bbZ\) by \(\mu_n\) over \(\Spec(R)\)
in the category of finite flat group schemes over 
\(\Spec(R)\) which are killed by \(n\).

\vn
{\bf Notation.}\enspace
Let \(R\) be an Artinian local \(W(k)\)-algebra, where \(k\supset \Fp\)
is an algebraically closed field.
Let \(X\to \Spec(R)\) be an ordinary Barsotti-Tate group
such that the closed fiber \(X_0:=X\times_{\Spec(R)}\Spec(k)\) 
is an ordinary BT-group over \(k\).
Denote by \(\,q(X/R;\cdot,\cdot)\) the \(\Zp\)-bilinear map
\[
\,q(X/R;\cdot,\cdot): {\rm T}_p({X_0}_{\rm et})\times 
{\rm T}_p({X_0^t}_{\rm et})
\to 1+\grm_R
\]
correspond to the deformation \(X\to \Spec(R)\) of the 
BT-group \(X_0\) as in Cor.\ \ref{STcoordav}.
Here we have used the natural isomorphism 
\(X^*({X_0}_{\rm mult})\cong{\rm T}_p({X_0^t}_{\rm et})\),
so that the Serre-Tate coordinates for the BT-group \(X\to \Spec(R)\)
is a \(\Zp\)-bilinear map \(\,q(X/R;\cdot,\cdot)\,\)
on \({\rm T}_p({X_0}_{\rm et})\times {\rm T}_p({X_0^t}_{\rm et})\).
The the abelian group
\(1+\grm_R\subset R^{\times}\) is regarded as a \(\Zp\)-module,
so ``\(\Zp\)-bilinear'' makes sense.
Let \({\rm can}: X_0\xrightarrow{\sim} (X_0^t)^t\) 
be the canonical isomorphism 
from \(X_0\) to its double Serre dual, and 
let \({\rm can}_*:{\rm T}_p({X_0}_{\rm et})\xrightarrow{\sim}
{\rm T}_p((X_0^t)^t_{\rm et})\) be the isomorphism
induced by \({\rm can}\).
\medbreak

The relation between the Serre-Tate coordinate \(q(X/R;\cdot,\cdot)\) 
of a deformation of \(X_0\) and the Serre-Tate coordinates
\(q(X^t/R;\cdot,\cdot)\) of the Serre dual \(X^t\) of \(X\) is
given by \ref{ST_duality}. The proof is left as an exercise.

\begin{lemma}\label{ST_duality}
Let \(X\to \Spec(R)\) be an ordinary BT-group over an Artinian local 
\(W(k)\)-algebra \(R\).  
Then we have
\[
q(X;u,v_t)=q(X^t;v_t, {\rm can}_*(u))
\qquad
\forall u\in {\rm T}_p({X_0}_{\rm et}),\
\forall v\in {\rm T}_p({X_0^t}_{\rm et})\,.
\]
The same statement hold when the ordinary BT-group \(X\to \Spec(R)\)
is replaced by an ordinary abelian scheme \(A\to\Spec(R)\).
\end{lemma}
\medbreak

From the functoriality of the construction in \ref{STcoord},
it is not difficult to verify the following.

\begin{prop} Let \(X_0, Y_0\) be ordinary BT-groups
over a perfect field \(K\supset \Fpbar\). Let \(R\) be an
Artinian local ring over \(W(K)\).  Let \(X\to \Spec(R)\),
\(Y\to \Spec(R)\) be abelian schemes whose closed fibers
are \(X_0\) and \(Y_0\) respectively.  
Let \(q(X/R;\cdot,\cdot)\), \(q(Y/R;\cdot,\cdot)\) be the
Serre-Tate coordinates for \(X\) and \(Y\) respectively.
Let \(\beta:X_0\to Y_0\) be a homomorphism of abelian varieties
over \(k\).  Then \(\beta\) extends to a homomorphism 
from \(X\) to \(Y\) over \(\Spec(R)\) if and only if
\[
q(X/R;u,\beta^t(v_t))=q(Y/R;\beta(u), v_t)
\quad \forall u\in {\rm T}_p(X_0),\ \forall v_t\in{\rm T}_p(Y_0^t)\,.
\]
\end{prop}
\medbreak

\begin{cor}\label{nonpp}
Let \(A_0\) be an ordinary abelian variety over a perfect
field \(K\supset \Fp\). Let \(\lambda_0:A_0\to A_0^t\) be a polarization
on \(A_0\). Then 
\[
{\rm{Def}}((A_0,\lambda_0)/W(K))
\cong \underline{\rm Hom}_{\Zp}(S,\widehat{\bbG_m})\,,
\]
where 
\[
S:={\rm T}_p(A_0[p^{\infty}]_{\rm et})\otimes_{\Zp}
{\rm T}_p(A_0^t[p^{\infty}]_{\rm et})\,\left/\,
\big(u\otimes {\rm T}_p(\lambda_0)(v) -
v\otimes {\rm T}_p(\lambda_0)(u)\big)_{u,v\in
{\rm T}_p(A[p^{\infty}]_{\rm et})}\right.\,.
\]
\end{cor}
\medbreak

\subsection{\bf Exercise}\label{compdef}
Notation as in \ref{nonpp}.  Let \(p^{e_1}, \ldots, p^{e_g}\)
be the elementary divisors of the \(\Zp\)-linear map
\({\rm T}_p(\lambda_0): {\rm T}_p(A_0[p^{\infty}]_{\rm et})
\to {\rm T}_p(A_0^t[p^{\infty}]_{\rm et})\),
\(g=\dim(A_0)\), \(e_1\leq e_2\leq\cdots\leq e_g\).
The torsion submodule \(S_{\rm torsion}\) of \(S\)
is isomorphic to 
\(\,
\bigoplus_{1\leq i<j\leq g}\, (\Zp/p^{e_i}\Zp)\,.
\)

\begin{thm}\Bl\ \ {\rm (local rigidity)}  \label{localrig}
Let \(k\supset \Fp\) be an algebraically closed field.
Let \[T\cong (\widehat \bbG_m)^n=\Spf k[[u_1,\ldots,u_n]]\] 
be a formal torus,
with group law given by 
\[
u_i\ \mapsto \ u_i\otimes 1 + 1\otimes u_i + u_i\otimes u_i
\qquad
i=1,\ldots n\,.
\]
Let \(X=\Hom_k(\widehat \bbG_m,T)\cong \Zp^n\) be the 
cocharacter group of \(T\);
notice that \(\GL(X)\) operates naturally on \(T\).
Let \(G\) be a reductive linear algebraic subgroup of 
\(\GL(X\otimes_{\Zp} \Qp)\cong \GL_n\) over \(\Qp\).
Let \(Z\) be an irreducible closed formal subscheme of \(T\)
which is stable under the action of an open subgroup
\(U\) of \(G(\Qp)\cap \GL(X)\).
Then \(Z\) is a formal subtorus of \(T\).
\end{thm}

\noindent See Thm.\ 6.6 of \cite{Chai-Fam} for a proof of \ref{localrig};
see also \cite{Chai-Rigidity}.

\begin{cor}
Let \(x_0=[(A_0,\lambda_0,\eta_0)]\in\calA_{g,1,n}(\Fpbar)\)
be an \(\Fpbar\)-point of \(\calA_{g,1,n}\), where
\(\Fpbar\) is the algebraic closure of \(\Fp\).
Assume that the abelian variety \(A_0\) is ordinary.
Let \(Z(x_0)\) be the Zariski closure of the prime-to-\(p\)
Hecke orbit \(\calH^{(p)}_{\rm Sp_{2g}}(x_0)\) on 
\(\calA_{g,1,n}\).
The formal completion \(Z(x_0)^{/x_0}\) of 
\(Z(x_0)\) at \(x_0\) is a formal subtorus of the
Serre-Tate formal torus \(\calA_{g,1,n}^{/x_0}\).
\end{cor}

\proof This is immediate from \ref{localrig} and the
local stabilizer principal; see \ref{lsp} for the statement of
the local stabilizer principal.


\section{The Tate-conjecture: $\ell$-adic and $p$-adic}\label{3}
Most results of this section will not be used directly in our proofs.
However, this is such a beautiful part of mathematics that we 
like to tell more than we really need.\\
Basic references: \cite{Tate-Bourb} and \cite{Honda}; \cite{Tate-Endo}, 
\cite{J}, \cite{FO-isog}.

\subsection{} let $A$ be an abelian variety over a field $K$. 
The ring $\End(A)$ is an algebra over $\ZZ$, which has no torsion, and which 
is free of finite rank as $\ZZ$-module.  
We write $\End^0(A) = \End(A) \otimes_{\ZZ} \QQ$. 
Let $\mu: A \to A^t$ be a polarization. 
An endomorphism $x: A \to A$ defines $x^t: A^t \to A^t$. 
We define an anti-involution
$$
\dag: \End^0(A) \rightarrow \End^0(A), 
\quad\mbox{by}\quad x^t{\cdot}\mu = \mu{\cdot}x^{\dag},
$$
called the {\it Rosati-involution}. 
In case $\mu$ is a principal polarization this maps $\End(A)$  into itself.

The Rosati involution is {\it positive definite} on $D := \End^0(A)$, 
meaning that $x \mapsto {\rm Tr}(x{\cdot}x^{\dag})$ 
is a positive definite quadratic form on $\End^0(A)$. 
Such algebras have been classified by Albert, see \ref{Albert}.

\subsection{}{\bf Definition.} {\it A field $L$ is said to be a} 
CM-{\it field if $L$ is a finite extension of $\QQ$} 
(hence $L$ is a number field), 
{\it there is a subfield $L_0 \subset L$ such that $L_0/\QQ$ 
is totally real} 
(i.e. every $\psi_0: L_0 \to \CC$ gives $\psi_0(L_0) \subset \RR$) 
{\it and $L/L_0$ is quadratic totally imaginary} (i.e. 
$[L:L_0] = 2$ 
and for every $\psi: L \to \CC$ we have $\psi(L) \not\subset \RR$).\\
{\bf Remark.} The quadratic extension $L/L_0$ gives an involution 
$\iota \in \Aut(L/L_0)$. For every embedding $\psi: L \to \CC$ 
this involution corresponds with the restriction of complex conjugation 
on $\CC$ to $\psi(L)$.

\vn
See \ref{smCM} for details on the definition 
``sufficiently many Complex Multiplications''.

\vn
Even more is known about the endomorphism algebra of an abelian 
variety over a finite field. Tate showed that 
\subsection{}{\bf Theorem} (Tate)
{\it an abelian variety over a finite field admits 
sufficiently many Complex Multiplications.}
This is equivalent with: 
{\it Let $A$ be a simple abelian variety over a finite field. 
Then there is a} CM-{\it field of degree $2{\cdot}\dm(A)$ 
contained in $\End^0(A)$.}\\ 
See \cite{Tate-Endo}, \cite{Tate-Bourb}, see \ref{struct2}. 
In particular this implies the following.

\vn
Let $A$ be an abelian variety over $\FF = \overline{\FF_p}$.  
Suppose that $A$ is simple, and hence that $\End^0(A)$ is 
a division algebra; this algebra has finite rank over $\QQ$. Then
\begin{itemize}
\item either $A$ is a supersingular elliptic curve, 
and $D := \End^0(A) = \QQ_{p,\infty}$, which is the (unique) 
quaternion algebra central over $\QQ$, which is unramified for 
every finite prime $\ell \not= p$, i.e. $D \otimes \QQ_{\ell}$ is 
the $2 \times 2$ matrix algebra over  $\QQ_{\ell}$, and $D/\QQ$ 
is ramified at $p$ and at $\infty$; here $D$ is of Albert Type III(1); 
\item or $A$ is not a supersingular elliptic curve; in this case $D$ 
is of Albert Type IV($e_0,d$) with $e_0{\cdot}d = g := \dm(A)$.
\end{itemize}

\vn
In particular (to be used later). 
\subsection{}{\bf Corollary.}  {\it Let $A$ be an abelian variety over 
$\FF =   \overline{\FF_p}$. There exists $E = F_1 \times \cdots \times F_r$, 
a product of totally real fields, and an injective homomorphism 
$E \hookrightarrow \End^0(A)$ such that} $\dm_{\QQ}(E) = 
\dm(A)$.

\vn
Some examples.\\
(1) $E$ is a supersingular elliptic curve over $K = \FF_q$. 
Then either $D := \End^0(E)$ is isomorphic with $\QQ_{p,\infty}$, 
or $D$ is an imaginary quadratic field over $\QQ$ in which $p$ is not split.\\
(2) $E$ is a non-supersingular elliptic curve over $K = \FF_q$. 
Then  $D := \End^0(E)$ is an imaginary quadratic field over $\QQ$ 
in which $p$ is split.\\
(3) If $A$ is simple over $K = \FF_q$  such that $D := \End^0(A)$ is 
commutative, then $D = L = \End^0(A)$ is a CM-field of 
degree $2{\cdot}\dm(A)$ over $\QQ$.\\
(4) In characteristic zero the endomorphism algebra of a simple abelian 
variety which admits smCM is {\it commutative}. 
However in positive characteristic an Albert Type IV($e_0,d$) with 
$e_0>1$ can appear. For example, see \cite{Tate-Bourb}, 
page 67: for any prime number $p > 0$, and for any $g>2$ 
there exists a simple abelian variety over $\FF$ 
such that $D = \End^0(A)$ is a division algebra of rank $g^2$ 
over its center $L$, which is a quadratic imaginary field  over $\QQ$.

\subsection{} Let $D$ be an Albert algebra; i.e. $D$ is a division algebra, 
it is of finite rank over $\QQ$, and it has a positive definite 
$\dag: D \to D$ anti-involution. Suppose a characteristic is given. 
There exists a field $k$ of that characteristic, and an abelian variety 
$A$ over $k$ such that $\End^0(A) \cong D$, and such that $\dag$ 
is the Rosati-involution given by a polarization on $A$. 
This was proved by Albert, and by Shimura over $\CC$, 
see  \cite{Shim}, Theorem 5. In general this was proved by Gerritzen, 
see \cite{Ger}, for more references see \cite{FO-Endo}.

One can ask which possibilities we have for $\dm(A)$, once $D$ is given. 
This question is completely settled in characteristic zero. 
From properties of $D$ one can derive some properties of $\dm(A)$. 
However the answer to this question in positive characteristic is not yet 
complete.

\subsection{\bf Weil numbers and CM-fields.} {\bf Definition.}\label{Wn} 
{\it Let $p$ be a prime number, $n \in \ZZ_{>0}$; write $q = p^n$. 
A $q$-Weil number is an algebraic integer $\pi$ such that 
for every embedding $\psi: \QQ(\pi) \to \CC$ we have}
$$\mid\psi(\pi)\mid \quad = \quad \sqrt{q}.$$
We say that $\pi$ and $\pi'$ are {\it conjugated} if there exists 
an isomorphism $\QQ(\pi) \cong \QQ(\pi')$ mapping $\pi$ to $\pi'$.\\
{\bf Notation:}  $\pi \quad\sim\quad \pi'$. 
We write $W(q)$ for the set conjugacy classes of $q$-Weil numbers.

\subsection{}{\bf Proposition.}\label{P3} {\it Let $\pi$ be a $q$-Weil number. 
Then} \\
{\bf (I)} {\it either for at least one $\psi: \QQ(\pi) \to \CC$ 
we have  $\pm\sqrt{q} = \psi(\pi) \in \RR$; in this case we have:}\\
{\bf (Ie)} {\it $a$ is even, $\sqrt{q} \in \QQ$, 
and $\pi = + p^{n/2}$, or $\pi = - p^{n/2}$; or}\\
{\bf (Io)} {\it $a$ is odd, $\sqrt{q} \in \QQ(\sqrt{p})$, 
and $\psi(\pi) = \pm p^{n/2}$. 
\\
In particular in case} (I) 
{\it we have $\psi(\pi) \in \RR$ for every $\psi$}.\\
{\bf (II)} {\it Or for every  $\psi: \QQ(\pi) \to \CC$ 
we have $\psi(\pi) \not\in \RR$} (equivalently: for at least one $\psi$ 
we have $\psi(\pi) \not\in \RR$). 
{\it In case} (II) {\it the field $\QQ(\pi)$ is a} CM-{\it field}.\\
{\bf Proof.} Exercise.

\subsection{}{\bf Remark.} We see a characterization of $q$-Weil numbers. 
In case I we have $\pi = \pm\sqrt{q}$. If $\pi \not\in \RR$:
$$\beta := \pi + \frac{q}{\pi} \quad\mbox{\rm is totally real},$$
and $\pi$ is a zero of
$$T^2 - \beta{\cdot}T + q, \quad\mbox{\rm with}\quad  \beta < 2\sqrt{q}.$$
In this way it is easy to construct  $q$-Weil numbers. 

\subsection{}\label{pi} Let $K$ be a {\it finite field} and let $A$ be an  
abelian variety over $K$. Suppose $K = \FF_q$ with $q = p^n$.  
We have $F: A \to A^{(p)}$. Iterating this Frobenius map $n$ times, 
observing there is a canonical identification  $A^{(p^n)} = A$, 
we obtain $(\pi: A \to A) \in \End(A)$. If $A$ is simple, 
the subring $\QQ(\pi) \subset \End^0(A)$ is a subfield, 
and we can view $\pi$ as an algebraic integer.

\subsection{\bf Theorem}\label{WeilC} \ \Ex \ (Weil). {\it Let $K = \FF_q$ 
be a finite field, let $A$ be a simple abelian variety over $K$. 
Then $\pi$ is a $q$-Weil number.}\\
This is the famous ``Weil conjecture'' for an abelian variety over 
a finite field.\\
See \cite{Weil-C}, page 70; \cite{Weil-VA}, page 138; \cite{Mumford-AV}, 
Theorem 4 on page 206.

\subsection{}{\bf Exercise.} We indicate in which way this theorem, part of 
{\it the Weil conjecture for an abelian variety} can be proved.

{\it  Suppose that $A$ is a simple abelian variety over a field $K$, 
and let $L = {\rm Centre}(\End^0(A))$. A Rosati involution on $D$ induces 
complex conjugation on $L$ (for every embedding $L \hookrightarrow \CC$).}

{\it If moreover $K$ is a finite field, $\pi = \pi_A$ is a Weil-$q$-number.}

\vn
{\bf Proposition I.} {\it For a simple abelian variety $A$ 
over $K = \FF_q$ we have} 
$$\pi_A{\cdot}(\pi_A)^{\dag} \quad=\quad q.$$
Here $\dag: D \to D := \End^0(A)$ is the Rosati-involution.

\vspace{2mm}\n
One proof can be found in \cite{Mumford-AV}, formula (i) on page 206; 
also see \cite{C.S}, Coroll. 19.2 on page 144.

Another proof of (I)  can be given by duality: 
$$\left(F_{A/S}: A \to A^{(p)}\right)^t \quad=\quad 
V_{A^t/S}: (A^{(p)})^t \to A^t.$$
From this we see that 
$$\pi_{A^t} \cdot (\pi_A)^t = (F_{A^t})^n \cdot (V_{A^t})^n = p^n = q,$$
where we make the shorthand notation $F^n$ for the 
$n$ times iterated Frobenius morphism, and the same for $V^n$. 
See [GM], 5.21, 7.34 and Section 15.
\B

\vn
{\bf Proposition II.} {\it For any polarized abelian variety $A$ over
a field the Rosati-involution $\dag: D \to D := \End^0(A)$ 
is positive definite bilinear form on $D$,  
i.e. for any non-zero $x \in D$ we have ${\rm Tr}(x{\cdot}x^{\dag}) > 0$.}\\
See \cite{Mumford-AV}, Th. 1 on page 192, see \cite{C.S}, 
Th. 17.3 on page 138.
\B

\vn
Use Proposition I and Proposition II and prove the two statements above 
(and hence prove  Theorem \ref{WeilC}).

\subsection{\bf Remark.} Given $\pi = \pi_A$ of a  simple 
abelian variety over $\FF_q$ one can determine the structure 
of the division algebra $\End^0(A)$, see \cite{Tate-Bourb}, Th. 1. 
See \ref{struct2}.

\subsection{\bf Theorem}\label{H.T} \Ex \ (Honda and Tate). 
{\it By $A \mapsto \pi_A$ we obtain a bijective map
\begin{center}
\{abelian variety over $\FF_q$\}/$\sim_{\FF_q}$ 
$\quad\stackrel{\sim}{\longrightarrow}\quad  W_q/\sim$
\end{center}
between the set of $\FF_q$-isogeny classes of abelian varieties 
simple over $\FF_q$ and the set of conjugacy classes of $q$-Weil numbers.}

\subsection{}\label{cases}   Let $\pi$ be a Weil $q$-number. 
Let $\QQ \subset L \subset D$ be the central algebra determined by $\pi$. 
We remind the reader that 
$$
[L:\QQ] =: e, \quad [D:L] =: d^2, \quad 2g := e{\cdot}d. 
\quad\quad\quad\mbox{See  \ref{Sendo}.}
$$
 As we have seen in Proposition \ref{P3} there are three possibilities:

\vspace{2mm}\n
{\bf ($\RR$e)} {\it Either $\pi=\sqrt{q} \in \QQ$, and 
$q = p^n$ with $n$ an} {\bf even} {\it positive integer.}\quad 
\fbox{Type III(1), \ \ $g=1$}
\\
In this case  $\pi = + p^{n/2}$, or $\pi = - p^{n/2}$. 
Hence $L = L_0 = \QQ$. We see that $D/\QQ$ has rank 4, 
with ramification exactly at $\infty$ and at $p$. 
We obtain $g=1$, we have that $A=E$ is a supersingular elliptic curve, 
$\End^0(A)$ is of Type III(1), a definite quaternion algebra over $\QQ$. 
This algebra was denoted by Deuring as $\QQ_{p, \infty}$. 
Note that ``all endomorphisms of $E$ are defined over $K$'', i.e. for any 
$$\forall \quad K \subset K' \quad\mbox{we have}\quad 
\End(A) = \End(A \otimes K').$$

\vspace{2mm}\n
{\bf ($\RR$o)} {\it Or $q = p^n$ with $n$ an}\  {\bf odd} 
{\it positive integer, $\pi=\sqrt{q} \in \RR\notin \QQ$ 
.}\ \quad\ \  \ \fbox{Type III(2), \ $g=2$}
\\
In this case $L_0 = L = \QQ(\sqrt{p})$, a real quadratic field. 
We see that $D$ ramifies exactly at the two infinite places with 
invariants equal to $(n/2){\cdot}2/(2n) = 1/2$. 
Hence $D/L_0$ is a definite quaternion algebra over $L_0$; 
it is of Type III(2). We conclude $g=2$. 
If $K \subset K'$ is an extension of odd degree we 
have $\End(A) = \End(A \otimes K')$.  
If $K \subset K'$ is an extension of even degree $A \otimes K'$ 
is non-simple, it is $K'$-isogenous with a product of 
two supersingular elliptic curves, and $\End^0(A \otimes K')$ 
is a $2 \times 2$ matrix algebra over $\QQ_{p, \infty}$, and
$$\forall K' \quad\mbox{with }\quad 2 \mid [K':K] 
 \quad\mbox{we have}\quad \End(A) \not= \End(A \otimes K').$$

\vspace{2mm}\n
{\bf ($\CC$)} {\it For at least one embedding 
$\psi: \QQ(\pi) \to \CC$ we have $\psi(\pi) \not\in \RR$.}
\fbox{IV($e_0,d$), \ \ $g := e_0{\cdot}d$}
\\
In this case all conjugates of $\psi(\pi)$ are non-real. 
We can determine  $[D:L]$ knowing all $v(\pi)$ by \ref{struct2} (3); 
here $d$ is the greatest common divisor of all denominators of 
$[L_v:\QQ_p]{\cdot}v(\pi)/v(q)$, for all $v\mid p$. 
This determines $2g := e{\cdot}d$. 
The endomorphism algebra is of Type IV($e_0,d$). 
For $K = \FF_q \subset K' = \FF_{q^m}$ we have 
$$\End(A) = \End(A \otimes K') \quad\Longleftrightarrow\quad  
\QQ(\pi) = \QQ(\pi^m).$$

\subsection{\bf Exercise/Remark.} Let $m, n \in  \ZZ$ with $m > n > 0$; 
write $g = m+n$ and $q = p^g$. Consider the polynomial 
$T^2 + p^nT + p^g$, 
and let $\pi$ be a zero of this polynomial.\\
(a) {\it Show that $\pi$ is a $p^g$-Weil number; 
compute the $p$-adic values of all conjugates of 
$\pi$.}\\
(b) By the previous theorem we see that $\pi$ defines the isogeny class 
of an abelian variety $A$ over $\FF_q$. It can be shown that $A$ has 
dimension $g$, and that $\cN(A) = (m,n) + (n,m)$, 
see \cite{Tate-Bourb}, page 98.  
{\it This gives a proof of a conjecture by Manin,} see \ref{MC}. 

\subsection{\bf $\ell$-adic monodromy.} (Any characteristic.) 
Let $K$ be a base field, any characteristic.  
Write $G_K = \Gal(K^{\rm sep}/K)$. Let $\ell$ be a prime number, 
not equal to char$(K)$. Note that this implies 
that $T_{\ell}(A)= {\rm lim}_{\leftarrow j} \ A[\ell^j]$ 
can be considered as a group isomorphic with $(\ZZ_{\ell})^{2g}$ 
with a continuous $G_K$-action.

\subsection{\bf Theorem}\label{lmon} \Ex \ (Tate, Faltings, 
and many others). {\it Suppose $K$ is of finite type over its 
prime field.} 
(Any characteristic.) {\it The canonical map
$$
\End(A) \otimes_{\ZZ} \ZZ_{\ell} \quad\stackrel{\sim}{\longrightarrow}
\quad \End(T_{\ell}(A)) \cong \End_{G_K}((\ZZ_{\ell})^{2g})
$$
is an isomorphism.
 }\\
This was conjectured by Tate.  In 1966 Tate proved this in case $K$ 
is a finite field, see \cite{Tate-Endo}. The case of function field 
in characteristic $p$ was proved by Zarhin and by Mori, 
see \cite{Zarhin-Isog}, \cite{Zarhin-endo}, \cite{Mori}; 
also see \cite{MorB}, pp. 9/10 and VI.5 (pp. 154-161). 

The case $K$ is a number field was open for a long time; 
it was finally proved by Faltings in 1983, see \cite{Faltings-1983}. 
For the case of a function field in characteristic zero, see \cite{F.W},  
Th. 1 on page 204.

\subsection{\bf Remark.}\ \Ex \  The previous result holds over 
a number field, 
but the Tate map need not be an isomorphism for an abelian variety 
over a local field.

\vn
{\bf Example,}  Lubin and Tate, see \cite{L.T}, 3.5; 
see \cite{FO-Bowdoin}, 14.9. 
{\it There exists a finite extension $L \supset 
\QQ_p$ and an abelian variety over $L$ such that}
$$ 
\End(A) \otimes_{\ZZ} \ZZ_{\ell} \quad\subsetneqq\quad \End(T_{\ell}(A)).
$$
We give details of a proof of this fact (slightly more general than 
in the paper by Lubin and Tate). Choose a prime number $p$, 
and choose a supersingular elliptic curve $E_0$ over $K = \FF_q$ 
such that the endomorphism ring $R := \End(E_0)$ has rank $4$ over 
$\ZZ$. In that case and $R$  is a maximal order in the 
endomorphism algebra $D := \End^0(E_0)$, which is a 
quaternion division algebra central over $\QQ$. 
Let $I$ be the index set of all subfields $L_i$ of $D$, and let 
$$
\Lambda := \cup_{i \in I}  (L_i \otimes \QQ_p) 
\quad\subset\quad D \otimes \QQ_p.
$$
{\bf Claim.}  
$$
\Lambda  \quad\subsetneqq\quad D_p := D \otimes \QQ_p.
$$
Indeed, the set $I$ is countable, and $[L_i : \QQ] \leq 2$ for every $i$. 
Hence $\Lambda$ is a countable union of $2$-dimensional $\QQ_p$-vector 
spaces  inside $D_p \cong (\QQ_p)^4$. The claim follows.

\vn
Hence we can choose $\psi_0 \in R_p := R \otimes \ZZ_p$ 
such that $\psi_0 \not\in \Lambda$: first choose $\psi_0'$ in $D_p$ 
outside $\Lambda$, then multiply with a power of $p$ 
in order to make $\psi_0 = p^n{\cdot}\psi_0'$ integral.  

Consider $X_0 := E_0[p^{\infty}]$. The pair $(X_0,\psi_0)$ can be lifted 
to characteristic zero, see \cite{FO-Bowdoin}, Lemma 14.7, 
hence to $(X,\psi)$ defined over an order in a finite extension $L$ of 
$\QQ_p$. We see that $\End^0(X) = \QQ_p(\psi)$, which is 
a quadratic extension of $\QQ_p$.   By the theorem of Serre and Tate, 
see \ref{thm:ST}, we derive an elliptic curve $E$, 
which is a lifting of $E_0$, such that $E[p^{\infty}] = X$. 
Clearly $\End(X) \otimes \ZZ_p \subset \End(X)$. \\
{\bf Claim.} $\End(E) = \ZZ.$\\
In fact, if $\End(E)$ would be bigger, we would have 
$\End(E) \otimes \ZZ_p = \End(X)$. 
Hence $\psi \in \End^0(E) \subset \Lambda$, which is a contradiction. 
This finishes the proof of the example:
$$
\End(E) = \ZZ \quad\mbox{and}\quad 
\dm_{\QQ_p} \End^0(X) = 2 \quad\mbox{and}\quad 
\End(E) \otimes \ZZ_p \subsetneqq \End(X).
$$

\vn
However, surprise, in the ``anabelian situation'' of a 
hyperbolic curve over a $p$-adic field, the analogous situation, 
gives an isomorphism for fundamental groups, see \cite{Moch}.
We see:  the Tate conjecture as in \ref{lmon} does not hold over  
$p$-adic fields but the Grothendieck ``anabelian conjecture'' is true 
for hyperbolic curves over $p$-adic fields. Grothendieck took care to 
formulate his conjecture with a number field as base field, 
see \cite{Esq}, page 19; we see that this care is necessary for 
the original Tate conjecture for abelian varieties, 
but for algebraic curves this condition can be relaxed.

\subsection{}
 We like to have a $p$-adic analogue of \ref{lmon}. For this purpose 
it is convenient to have $p$-divisible groups instead of 
Tate-$\ell$-groups, 
and in fact the following theorem now has been proved to be true.

\subsection{\bf Theorem}\label{ext} \ \Th \ \& \Bl \ (Tate, De Jong).
 {\it Let $R$ be an integrally closed, Noetherian 
integral domain with field 
of fractions $K$.} (Any characteristic.) {\it Let $X, Y$ be 
$p$-divisible groups 
over $\Spec(R)$. Let $\beta_K: X_K \to Y_K$ be a homomorphism. 
There exists 
(uniquely) $\beta: X \to Y$ over $\Spec(R)$ extending  $\beta_K$.
}\\
This was proved by Tate, under the extra assumption that 
the characteristic of $K$ is zero. For the case char$(K) = p$, 
see \cite{J}, 1.2 and  \cite{AJdJ-ICM}, Th. 2 on page 261.

\subsection{\bf Theorem}\label{pmon}\ \Bl \ (Tate and De Jong).  
{\it Let $K$ be a field finitely generated over $\FF_p$. 
Let $A$ and $B$ 
be abelian varieties over $K$. The natural map
$$
\Hom(A,B) \otimes \ZZ_p \quad\stackrel{\sim}{\longrightarrow}\quad
\Hom(A[p^{\infty}], B[p^{\infty}])
$$
is an isomorphism.}\\
This was proved by Tate in case $K$ is a finite field; 
a proof was written up in \cite{W.M}. The case of a function field over 
$\FF_p$ was proved by Johan de Jong, see \cite{J}, Th. 2.6. 
This case follows from the result by Tate 
and from the proceeding result \ref{ext} on extending homomorphisms.

\subsection{}{\bf Ekedahl-Oort strata.}\label{EO}\Bl \  In \cite{FO-EO} 
a new technique is developed, which will be used below. We sketch some of 
the details of that method. We will only indicate details relevant for the 
polarized case (and we leave aside the much easier unpolarized case).

\vn
A finite group scheme $N$ (say over a perfect field) for which 
$N[V] = \Image(F_N)$ and $N[F] = \Image(V_N)$ is called a 
BT$_1$ group scheme (a Barsotti-Tate group scheme truncated at level 1). 
By a theorem of Kraft, independently observed by Oort, for a given rank 
over an algebraically closed field $k$ the number of isomorphism classes of 
BT$_1$ group schemes is finite, see \cite{Kraft}. 
For any abelian variety $A$, the group scheme $A[p]$ is a BT$_1$ 
group scheme. A principal polarization $\lambda$ on $A$ induces 
a form on $A[p]$, and the pair $(A,\lambda)[p]$ is a polarized BT$_1$ 
group scheme, see \cite{FO-EO}, Section 9  (there are subtleties 
in case $p=2$: the form has to be taken, over a perfect field, 
on the Dieudonn\'e module of $A[p]$). 

\vn
{\bf \ref{EO}.1} \quad  {\it The number of isomorphism classes of 
polarized BT$_1$ group schemes $(N,<,>)$ over $k$  
of given rank is finite; see the classification in} \cite{FO-EO}, 9.4.

\vn
Let $\va$ be the isomorphism type of a  polarized BT$_1$ group scheme. 
Consider $S_{\va} \subset \cA_{g,1}$, the set of all $[(A,\lambda)]$ 
such that $(A,\lambda)[p]$ is geometrically belongs to the isomorphism class  
$\va$. 

\vn
{\bf \ref{EO}.2} \quad  {\it It can be shown that this is a 
locally closed set}, called an EO-stratum. 
We obtain $\cA_{g,1} = \sqcup_{\va} \ \ S_{\va}$, 
a disjoint union of local closed sets. {\it This is a stratification}, 
in the sense that {\it the boundary of a stratum is a union of 
lower dimensional strata.}

\vn
On of the main theorems of this theory is that \\
{\bf \ref{EO}.3} \quad {\it for every $\va$ the set $S_{\va}$ is}  
quasi-affine (i.e. open in an affine scheme), see \cite{FO-EO}, 1.2. 

\vn
The finite set $\Phi_g$ of such isomorphism types has two partial orderings, 
see \cite{FO-EO}, 14.3. One of these, denoted by $\va \subset \va'$, 
is defined by the property that $S_{\va}$ is contained in 
the Zariski closure of  $S_{\va'}$.

\subsection{}{\bf An application.}\label{Appl} {\it Let $x \in \cA_{g,1}$. Let 
$$
\left(\cH_{\ell}(x)\right)^{\rm Zar} \quad\subset\quad \cA_{g,1}
$$
be the Zariski closure. This closed set in $\cA_{g,1}$  
contains a supersingular point.}\\
Use \ref{EO} and \ref{ssfinite}.
\B

\section{Dieudonn\'e modules and Cartier modules}\label{4}
\n
In this section we explain the theory of Cartier modules and Dieudonn\'e
modules.  These theories provide equivalence of categories of
geometric objects such as commutative smooth formal groups or
Barsotti-Tate groups on the one side, and modules over
certain non-commutative rings on the other side.
As a result, questions on commutative smooth formal groups or 
Barsotti-Tate groups, which are apparently non-linear in nature,
are translated into questions in linear algebras over rings.
Such results are essential for any serious computation.
\smallbreak

There are many versions and flavors of Dieudonn\'e theory.
We explain the Cartier theory for commutative smooth formal 
groups over general commutative rings, and the covariant
Dieudonn\'e modules for Barsotti-Tate groups
over perfect fields of characteristic \(p>0\).
Since the Cartier theory works over general commutative rings,
one can ``write down'' explicit deformations over complete
rings such as \(k[[x_1,\ldots,x_n]]\) or \(W(k)[[x_1,\ldots,x_n]]\),
something rarely feasible in algebraic geometry.
\medbreak

\noindent{\scshape Remarks on notation:}\enspace
(i) In the first part of this section, on Cartier theory,
\(k\) denotes a commutative ring with \(1\),
or a commutative \(\bbZ_{(p)}\)-algebra with \(1\).
This differs from the convention in the rest of this article
that \(k\) stands for an algebraically closed field.
\smallbreak

\noindent (ii) In this section, we used \(V\) and \(F\) as elements
in the Cartier ring \(\Cartp{k}\) or the smaller Dieudonn\'e
ring \(R_K\subset \Cartp{K}\) for a perfect field \(K\).
In the rest of this article, the notation \(\calV\) and \(\calF\)
are used; \(\calV\) corresponds to the relative Frobenius morphism
and \(\calF\) corresponds to the Verschiebung morphisms
for commutative smooth formal groups or Barsotti-Tate groups
over \(K\).
\bigbreak

\noindent{\bf A synopsis of Cartier theory}
\medbreak

Let \(k\) be a commutative ring with \(1\) over \(\bbZ_{(p)}\).
The main theorem of Cartier theory says that there is an
equivalence between the category of commutative smooth formal
groups over \(k\) and the category of left modules over
a non-commutative ring \(\Cartp{k}\) satisfying certain
conditions.  See \ref{thm:main_local} for a precise statement.
\smallbreak

The Cartier ring \(\Cartp{k}\) plays a crucial role.  
The is a topological ring which contains elements \(V\), \(F\) 
and \(\{\langle a\rangle\,\mid\,a\in k\).  These
elements form a set of topological generators, so that every element
of \(\Cartp{k}\) has a unique expression as a convergent
sum in the following form
\[
\sum_{m,n\geq 0}\, V^m \langle a_{mn}\rangle F^n\,,
\quad a_{mn}\in k\,, \forall\,m\ \exists\, C_m>0\ {\rm s.t.\ }
a_{mn}=0\ {\rm if\ } n\geq C_m\,.
\]
These topological generators satisfy the following commutation relations.
\begin{itemize}
\item \(F\,\langle a\rangle = \langle a^p\rangle\,F\)
for all \(a\in k\).

\item  \(\langle a \rangle V= V\,\langle a^p\rangle\)
for all \(a\in k\).

\item  \(\langle a\rangle \,\langle b\rangle=
\langle ab\rangle\) for all \(a,b\in k\).

\item \(FV=p\).

\item  \(V^m\langle a \rangle F^m\,V^n\langle b \rangle F^n
=p^r\,V^{m+n-r}\langle a^{p^{n-r}} b^{p^{m-r}}\rangle F^{m+n-r}\)
for all \(a,b\in k\) and all \(m,n\in \bbN\), where \(r=\min\{m,n\}\).
\end{itemize}
Moreover the ring of \(p\)-adic Witt vectors \(W_p(k)\) is embedded
in \(\Cartp{k}\) by the formula
\[
W_p(k)\ \ni\ \underline{c}=(c_0, c_1,c_2,\ldots)
\ \longmapsto\
\sum_{n\geq 0}V^n \langle c_n\rangle F^n
\ \in\ \Cartp{k}\,.
\]
The topology of \(\Cartp{k}\) is given by the decreasing filtration
\[{\rm Fil}^n\left(\Cartp{k}\right):=V^n\cdot \Cartp{k}\,,\] 
making \(\Cartp{k}\) a complete and separated topological ring.
Under the equivalence of categories mentioned above,
a left \(\Cartp{k}\) module corresponds to a finite dimensional
smooth commutative formal group \(G\) over \(k\)  if and only if
\begin{itemize}
\item \(V:M\to M\) is injective,

\item \(M\xrightarrow{\sim} \varprojlim_n M/V^n M\), and

\item \(M/VM\) is a projective \(k\)-module of finite type.

\end{itemize}
If so, then \(\Lie(G/k)\cong M/VM\), and \(M\) is a finitely
generated \(\Cartp{k}\)-module.
See \ref{defn_local_cartring} for the definition of \(\Cartp{k}\), 
\ref{exer_local_commutation} for the commutation relations in \(\Cartp{k}\),
and \ref{prop_elt_local} for some other properties of \(k\).
We strongly advice the readers with no prior experience
with Cartier theory to accept the ``big black box'' 
as described in the previous paragraph and use the materials in
\ref{def_grp} -- \ref{thm:main_local} as a dictionary 
only when necessary.
Instead, it would be more helpful to 
get familiar first with the structure of the ring \(\Cartp{k}\) 
in the case when \(kA\) is a perfect field and 
play with some examples of finitely generated modules over \(\Cartp{k}\),
in conjunction with the theory of covariant Dieudonn\'e modules
over perfect fields in characteristic \(p\).
\bigbreak

\noindent{\scshape References for Cartier theory.}\enspace
We highly recommend \cite{Zink-Cart}, where the approach
in \S 2 of \cite{raynaud:p-torsion} is fully developed.
Other references for Cartier theory are \cite{lazard:formal_group}
and \cite{hazewinkel}.
\bigbreak

\begin{defn}\label{def_grp}
Let \(k\) be a commutative ring with \(1\).
\begin{itemize}
\rmitem[(1)] Let \({\mathfrak{Nilp}}_k\) be the category of all
nilpotent \(k\)-algebras, consisting of all commutative 
\(k\)-algebras \(N\) without unit such that \(N^n=(0)\) for some
positive integer \(n\).

\rmitem[(2)] A \emph{commutative smooth formal group} over \(k\)
is a covariant functor \(G: {\mathfrak{Nilp}_k} \to {\mathfrak{Ab}}\) 
from \({\mathfrak{Nilp}_k}\) to the category of
all abelian groups such that the following properties are satisfied.
\begin{itemize}
\item[\(\bullet\)] \(G\) commutes with finite inverse limits;

\item[\(\bullet\)] \(G\) is formally smooth, i.e.\ every surjection
\(N_1\to N_2\) in \(\mathfrak{Nilp}_k\) induces a surjection
\(G(N_1)\to G(N_2)\);

\item[\(\bullet\)] \(G\) commutes with arbitrary direct limits.
\end{itemize}

\rmitem[(3)] The \emph{Lie algebra} of a commutative smooth
formal group \(G\) is defined to be
\(G(N_0)\), where \(N_0\) is the object in \(\mathfrak{Nilp}_k\)
whose underlying \(k\)-module is \(k\), and \(N_0^2=(0)\).
\end{itemize}
\end{defn}
\medbreak

\noindent{\bf Remark.}\enspace Let \(G\) be a commutative smooth formal group
over \(k\), then \(G\) extends uniquely to a functor \(\tilde{G}\)
on the category \(\mathfrak{ProNilp}_k\)
of all filtered projective system of nilpotent \(k\)-algebras
which commutes with filtered projective limits.
This functor \(\tilde{G}\) is often denoted \(G\) by abuse of notation.
\medbreak

\noindent{\bf Example.}\enspace
Let \(A\) be a commutative smooth group scheme over \(k\). 
For every nilpotent \(k\)-algebra \(N\), denote by \(k\oplus N\)
the commutative \(k\)-algebra with multiplication given by
\[(u_1,n_1)\cdot (u_2,n_2)= (u_1 u_2, u_1 n_2+u_2 n_1+ n_1 n_2)
\qquad \forall u_1,u_2\in k\ \ \forall n_1,n_2\in N\,.
\]
The functor which sends an object \(N\) in \(\mathfrak{Nilp}_k\)
to the abelian group 
\[{\rm Ker}\left(A(k\oplus N)\to A(k)\right)\] 
is a commutative smooth
formal group over \(k\), denoted by \(\widehat{A}\).
For instance we have \[\widehat{\bbG_a}(N)=N\qquad {\rm and}\qquad
\widehat{\bbG_m}(N)=1+N\subset (k\oplus N)^{\times}\]
for all \(N\in {\rm Ob}(\mathfrak{Nilp}_k)\).
\medbreak

\begin{defn} \label{def_Lambda}
We define a \emph{restricted version} of the smooth formal group
attached to the universal Witt vector group over \(k\),
denoted by \(\Lambda_k\), or \(\Lambda\) when the base ring \(k\)
is understood.
\[
\Lambda_k(N)= 1+ t\,k[t]\!\otimes_k\! N\subset ((k\oplus N)[t])^{\times}\,
\qquad \forall\  N\in {\rm Ob}({\mathfrak{Nilp}}_k)\,.
\] 
In other words, the elements of \(\Lambda(N)\) consists of all
polynomials of the form \(1+u_1\,t+u_2\,t^2+\cdots+u_r\,t^r\)
for some \(r\geq 0\), where \(u_i\in N\) for \(i=1,\ldots,r\).
The group law of \(\Lambda(N)\) comes from multiplication in
the polynomial ring \((k\oplus N)[t]\) in one variable \(t\).
\end{defn}
\smallbreak

\noindent{\bf Remark.}\enspace
(i) The formal group \(\Lambda\) will play the role of a free generator
in the category of (smooth) formal groups.\\
(ii) When we want to emphasize that the polynomial 
\(1+\sum_{i\geq 1}\,u_i\,t^i\)
is regarded as an element of \(\Lambda(N)\), we
denote it by \(\lambda(1+\sum_{i\geq 1}\,u_i\,t^i)\).
\medbreak

\subsection{\bf Exercise} Let \(k[[X]]^+=X\,k[[X]]\) 
be the set of all formal power series
over \(k\) with constant term \(0\); it is an object in 
\(\mathfrak{ProNilp}_k\).
Show that
\begin{eqnarray*}
\Lambda(k[[X]]^+)\ =\hspace{5in}\\
\hspace{2pt}\left\{ \left.
\prod_{m,n\geq 1}\,(1-a_{mn}\,X^m\,t^n)\,\right\vert\,a_{m,n}\in k,\
\forall\, m\,\exists\, C_m>0\ {\rm s.t.\ }a_{mn}=0{\ \rm if\ }n\geq C_m
\right\}
\end{eqnarray*}

\medbreak

\begin{thm}\Bl\label{yoneda}
Let 
\(H:{\mathfrak{Nilp}}_k \to {\mathfrak{Ab}}\) be a commutative
smooth formal group over \(k\).
Let \(\Lambda=\Lambda_k\) be the functor defined in 
\ref{def_Lambda}.
The map
\[
Y_H: \Hom(\Lambda_k, H) \to H(k[[X]]^+)
\]
which sends each homomorphism \(\alpha:\Lambda \to H\) 
of group-valued functors to the element 
\[\alpha_{_{k[[X]]^+}}(1-Xt) \in H(k[[X]]^+)\]
is a bijection.
\end{thm}
\medbreak

\noindent{\bf Remark.}\enspace
The formal group \(\Lambda\) is in some sense 
a free generator of the additive category of
commutative smooth formal groups, a phenomenon reflected in 
Thm. \ref{yoneda}.
\medbreak

\begin{defn}\label{defn_cartierring}
{\rm (i)} Define \(\Cart{k}\) to be \(\left(\End(\Lambda_k)\right)^{\rm op}\),
the opposite ring of the endomorphism ring of the smooth formal group
\(\Lambda_k\).
According to Thm.\ \ref{yoneda}, for every weakly symmetric 
functor \(H:{\mathfrak Nilp}_k\to {\mathfrak Ab}\),
the abelian group
\(H(k[[X]]^+)=\Hom(\Lambda_k, H)\) is a {\em left} module over
\(\Cart{k}\).

\noindent {\rm (ii)} We define some special elements of the Cartier ring
\(\Cart{k}\), naturally identified with \(\Lambda(k[[X]])\) via
the bijection
\(Y=Y_{\Lambda}:\End(\Lambda)\xrightarrow{\sim}\Lambda(k[[X]]^+)\)
in Thm.\ \ref{yoneda}.
\begin{itemize}
\item \(V_n:=Y^{-1}(1-X^n\,t)\,, n\geq 1\),

\item \(F_n:=Y^{-1}(1-X\,t^n)\,, n\geq 1\),

\item \([c]:= Y^{-1}(1-c\,X\,t)\,, c\in k\).
\end{itemize}
\end{defn}
\medbreak

\noindent{\bf Corollary.}\enspace
For every commutative ring with \(1\) we have
\[
\Cart{k}=\left\{ \sum_{m,n\geq 1} V_m\,[c_{mn}]\,F_n
\,\bigg\vert\, c_{mn}\in k,\ \forall\,m\, \exists\,C_m>0
{\ \rm s.t.\ }c_{mn}=0 {\ \rm if\ } n\geq C_m
\right\}
\]
\medbreak

\begin{prop}\Bl \label{prop_easy_id}
The following identities hold in \(\Cart{k}\).
\begin{itemize}
\rmitem[(1)] \(V_1=F_1=1\), \(F_n\,V_n=n\).
\rmitem[(2)] \([a]\,[b] =[ab]\) for all \(a,b\in k\)
\rmitem[(3)] \([c]\,V_n= V_n\,[c^n]\), \(F_n\,[c]= [c^n]\,F_n\)
for all \(c\in k\), all \(n\geq 1\).
\rmitem[(4)] \(V_m\,V_n=V_n\,V_m=V_{mn}\), 
\(F_m\,F_n=F_n\,F_m=F_{mn}\) for all \(m,n\geq 1\).
\rmitem[(5)] \(F_n\,V_m=V_m\,F_n\) if \((m,n)=1\).
\rmitem[(6)] \((V_n[a]F_n)\cdot (V_m[b]F_m)
= r\,V_{\frac{mn}{r}}\left[a^{\frac{m}{r}} b^{\frac{n}{r}}\right]
F_{\frac{mn}{r}}\), \(r=(m,n)\), for all \(a,b\in k\), \(m,n\geq 1\).
\end{itemize}
\end{prop}
\medbreak

\begin{defn}\label{defn_fil_ring}
The ring \(\Cart{k}\) has a natural filtration 
\({\rm Fil}^{\bullet}\Cart{k}\) by right ideals, where
\({\rm Fil}^j\Cart{k}\) is defined by
$$
{\rm Fil}^j\Cart{k}=\\ 
\left\{
\sum_{m\geq j} \sum_{n\geq 1}\,V_{m}[a_{mn}]F_{n}\,
\left\vert\,a_{mn}\in k,\ 
\forall\,m\geq j,\, \exists\, C_m>0\ {\rm s.t.\ }a_{mn}=0\ 
{\rm if}\ n\geq C_m
\right.
\right\}
$$
for every integer \(j\geq 1\).
The Cartier ring \(\Cart{k}\) is complete with respect to the
topology given by the above filtration. Moreover each right ideal
\({\rm Fil}^j\Cart{k}\) is open and closed in \(\Cart{k}\).
\end{defn}
\medbreak

\noindent{\bf Remark.}\enspace
The definition of the Cartier ring gives a functor 
\[
k\ \longmapsto \Cart{k}
\]
from the category of commutative rings with \(1\) to the category of
complete filtered rings with \(1\).
\medbreak

\begin{defn}\label{defn_red_mod}
Let \(k\) be a commutative ring with \(1\).
\begin{itemize}
\rmitem[(1)] A \(V\!\)-{\em reduced} left \(\Cart{k}\)-module is a 
left \(\Cart{k}\)-module \(M\) together with a separated decreasing 
filtration of \(M\) 
\[
M={\rm Fil}^1 M\supset {\rm Fil}^2 M\supset \cdots
{\rm Fil}^n M \supset {\rm Fil}^{n+1}\supset \cdots
\]
such that each \({\rm Fil}^n M\) is an abelian subgroup of \(M\) and
\begin{itemize}
\rmitem[(i)] \((M, {\rm Fil}^{\bullet}M)\) is complete with respect
to the topology given by the filtration \({\rm Fil}^{\bullet}M\).
In other words, the natural map 
\(\displaystyle {\rm Fil}^n M\to \varprojlim_{m\geq n}\,
({\rm Fil}^n M/{\rm Fil}^m M)\)
is a bijection for all \(n\geq 1\).

\rmitem[(ii)] \(V_m\cdot {\rm Fil}^n M\subset {\rm Fil}^{mn} M\)
for all \(m,n\geq 1\).

\rmitem[(iii)] The map \(V_n\) induces a bijection
\(V_n: M/{\rm Fil^2}M \xrightarrow{\sim} {\rm Fil}^n M/{\rm Fil}^{n+1}M\)
for every \(n\geq 1\).

\rmitem[(iv)] \([c]\cdot {\rm Fil}^n M\subset {\rm Fil}^n M\)
for all \(c\in k\) and all \(n\geq 1\).

\rmitem[(v)] For every \(m, n\geq 1\), there exists an \(r\geq 1\)
such that \(F_m\cdot {\rm Fil}^r M\subset {\rm Fil}^n M\).
\end{itemize}

\rmitem[(2)] A \(V\!\)-reduced left \(\Cart{k}\)-module 
\((M,{\rm Fil}^{\bullet}M)\)
is \(V\!\)-{\em flat} if \(M/{\rm Fil}^2 M\) is a flat \(k\)-module.
The \(k\)-module \(M/{\rm Fil}^2 M\) is defined to be the 
\emph{tangent space} of \((M,{\rm Fil}^{\bullet}M)\).
\end{itemize}
\end{defn}
\medbreak

\begin{defn}
Let \(H:\mathfrak{Nilp}_k\to \mathfrak{Ab}\) be a 
commutative smooth formal group over \(k\).
The abelian group \({\rm M}(H):=H(k[[X]]^+)\) 
has a natural structure as a
left \(\Cart{k}\)-module according to Thm.\ \ref{yoneda}
The \(\Cart{k}\)-module \({\rm M}(H)\) has a natural filtration, 
with 
\[{\rm Fil}^n{\rm M}(H):= \Ker(H(k[[X]]^+)\to H(k[[X]]^+/X^n k[[X]]))\,.\]
We call the pair \(({\rm M}(H), {\rm Fil}^{\bullet}{\rm M}(H))\) the
\emph{Cartier module attached to \(H\)}.
\end{defn}
\medbreak

\begin{defn}\label{defn_red_tensor}
Let \(M\) be a \(V\!\)-reduced left \(\Cart{k}\)-module and let 
\(Q\) be a right \(\Cart{k}\)-module.
\begin{itemize}
\rmitem[(i)] For every integer \(m\geq 1\), let
\(Q_m:={\rm Ann}_{Q}({\rm Fil}^m\Cart{k})\) be the subgroup of
\(Q\) consisting of all elements \(x\in Q\) such that
\(x\cdot {\rm Fil}^m\Cart{k}=(0)\). Clearly we have
\(Q_1\subseteq Q_2 \subseteq Q_3\subseteq \cdots\).

\rmitem[(ii)] For each \(m,r\geq 1\), define
\(Q_m\odot M^r\) to be the image of \(Q_m\otimes {\rm Fil}^r M\)
in \(Q\otimes_{\Cart{k}} M\). 
\smallbreak

Notice that if \(r\geq m\) and \(s\geq m\), then
\(Q_m\odot M^{r}= Q_m\odot M^{s}\). 
Hence \(Q_m\odot M^m\subseteq Q_n\odot M^n\) if
\(m\leq n\). 

\rmitem[(iii)] Define the \emph{reduced tensor product}
\(Q\overline{\otimes}_{\Cart{k}} M\) by
\[Q\overline{\otimes}_{\Cart{k}} M
 = Q\otimes_{\Cart{k}}M \left/\big(\bigcup_m (Q_m\odot M^m)\big)\right.
\,.\]

\end{itemize}
\end{defn}
\medbreak

\noindent{\bf Remark.} The reduced tensor product is used to construct
the arrow in the ``reverse direction'' in the equivalence of category
in \ref{thm_main} below.
\medbreak

\begin{thm}\Bl\label{thm_main}
Let \(k\) be a commutative ring with \(1\). There is
a canonical equivalence of categories, between the category
of smooth commutative formal groups over \(k\) 
as defined in {\rm \ref{def_grp}}
and the category of \(V\)\!-flat \(V\)\!-reduced left
\(\Cart{k}\)-modules, defined as follows. 
\[
\xymatrix@R=2pt{
\{{\rm smooth\ formal\ groups\ over\ } k\}\hspace{0.2cm}
\ar@1{->}[r]^{\sim} & 
\hspace{0.2cm}\{\mbox{V-}{\rm flat}\mbox{ V-}{\rm reduced\ left\ } 
\Cart{k}\mbox{-}{\rm mod}\}\\
\hspace{1.5cm}G\hspace{1.5cm} \ar@{|->}[r] 
& \hspace{1cm} {\rm M}(G)=\Hom(\Lambda,G)\\
\Lambda\overline{\otimes}_{\Cart{k}}M\hspace{0.7cm}\hspace{0.5cm}
&\hspace{2.3cm}M\hspace{1.2cm}\ar@{|->}[l]
}
\] 
Recall that \({\rm M}(G)=\Hom(\Lambda, G)\) is canonically isomorphic to
\(G(X\,k[[X]])\), the group of all formal curves in the smooth formal
group \(G\). The reduced tensor product 
\(\Lambda\overline{\otimes}_{\Cart{k}}M\) is the functor whose value
at any nilpotent \(k\)-algebra \(N\) is 
\(\Lambda(N)\overline{\otimes}_{\Cart{k}}M\).
\end{thm}
\medbreak

The Cartier ring \(\Cart{k}\) contains the ring of universal Witt
vectors \(\widetilde{W}(k)\) as a subring which contains the unit
element of \(\Cart{k}\).  

\begin{defn}\label{defn_univwitt}
(1) The \emph{universal Witt vector group} \(\widetilde W\) 
is defined as the
functor from the category of all commutative algebras with \(1\)
to the category of abelian groups such that
\[
\widetilde W(R)=1+T\,R[[T]]\subset R[[T]]^{\times}
\]
for every commutative ring \(R\) with \(1\). 
\smallbreak

When we regard a formal power series \(1+\sum_{m\geq 1}\,u_m\,T^m\) in
\(R[[T]]\) as an element of \(\widetilde W(R)\), we use the notation
\(\omega(1+\sum_{m\geq 1}\,u_m\,T^m)\).
It is easy to see that every element of 
\(\widetilde W(R)\) has a unique expression as
\[\omega\left(\prod_{m\geq 1}\,(1-a_m\,T^m)\right)\,.\]
Hence
\(\widetilde W\) is isomorphic to \(\Spec\,\bbZ[x_1,x_2,x_3,\ldots]\) as a 
scheme; the \(R\)-valued point such that \(x_i\mapsto a_i\) is denoted
by \(\omega(\underline{a})\), where \(\underline a\) is short for
\((a_1,a_2,a_3,\ldots)\).  In other words, 
\(\omega(\underline a)=\omega(\prod_{m\geq 1}\,(1-a_m\,T^m))\)
\medbreak

\noindent (2) The group scheme \(\widetilde W\) has a natural structure 
as a ring scheme, such that multiplication on
\(\widetilde W\) is determined by the formula
\[
\omega(1-a\,T^m)\cdot \omega(1-b\,T^n)=
\omega\left(\bigl(\,1-a^{\frac{n}{r}}\,b^{\frac{m}{r}}\,T^{\frac{mn}{r}}
\bigr)^r\right)\,,
\quad
{\rm where\ } r={\rm gcd}(m,n)\,.
\]
\medbreak

\noindent (3) There are two families of endomorphisms of the group scheme
\(\widetilde W\): \(V_n\) and \(F_n\), \(n\in \bbN_{\ge 1}\).
Also for each commutative ring \(R\) with \(1\) and each element
\(c\in R\) we have an endomorphism \([c]\)
of \(\widetilde W\times_{\Spec\,\bbZ}\Spec\,R\). 
These operators make \(\widetilde{W}(k)\) a left \(\Cart{k}\)-module;
they are defined as follows
\[
\begin{array}{rcccl}
V_n&:&\omega(f(T))&\mapsto& \omega(f(T^n))\\
\\
\,F_n&:&\omega(f(T))&\mapsto& 
\sum_{\zeta\in \mu_n} \omega(f(\zeta\,T^{\frac{1}{n}})) 
\qquad (\,{\rm formally}\,)\\
\\
{[c]}& :& \omega(f(T)) & \mapsto & \omega(f(cT))
\end{array}
\]
The formula for \(F_n(\omega(f(T)))\) 
means that \(F_n(\omega(f(T)))\) is defined as the unique element such that 
\(V_n(F_n(\omega(f(T))))=\sum_{\zeta\in \mu_n}\,\omega(f(\zeta\,T))\,.\)
\end{defn}
\medbreak

\subsection{\bf Exercise.} Show that the Cartier module of 
\(\widehat{\bbG_m}\) over
\(k\) is naturally isomorphic to
\(\widetilde{W}(k)\) as a module over \(\Cart{k}\).

\begin{prop}\Bl \label{prop_emb}
Let \(k\) be a commutative ring with \(1\). 
\begin{itemize}
\rmitem[(i)] The subset \(S\) of \(\Cart{k}\) consisting of all
elements of the form
\[
\sum_{n\geq 1}\,V_n [a_n] F_n\,, \quad a_n\in k\ \forall\,n\geq 1
\]
form a subring of \(\Cart{k}\).

\rmitem[(ii)] The injective map
\[
\widetilde{W}(k)\hookrightarrow \Cart{k},\qquad
\omega(\underline{a}) \mapsto \sum_{n\geq 1}\,V_n\,[a_n]\,F_n
\]
is an injective homomorphism of rings which sends \(1\) to \(1\);
its image is the subring \(S\) defined in (i).
\end{itemize}
\end{prop}
\medbreak

\begin{defn} \label{defn_idem}
It is a fact  
that every prime number
\(\ell\neq p\) is invertible in \(\Cart{\bbZ_{(p)}}\).
Define elements \(\epsilon_p\) and \(\epsilon_{p,n}\) 
of the Cartier ring \(\Cart{\bbZ_{(p)}}\) for \(n\in \bbN\),
\((n,p)=1\) by
\begin{eqnarray*}
\epsilon_p = \epsilon_{p,1}&=& \sum_{\tworows{(n,p)=1}{n\geq 1}}\,
\frac{\mu(n)}{n}V_n F_n
=\prod_{\tworows{\ell\neq p}{\ell\ {\rm prime}}}\,
\left(1-\frac{1}{\ell} V_{\ell} F_{\ell}\right)
\\
\epsilon_{p,n}&=&\frac{1}{n}V_n \epsilon_p F_n
\end{eqnarray*}
where \(\mu\) is the M\"obius function on \(\bbN_{\geq 1}\),
characterized by the following properties:
\(\mu(mn)=\mu(m)\,\mu(n)\) if \((m,n)=1\), and for every prime 
number \(\ell\) we have \(\mu(\ell)=-1\), \(\mu(\ell^i)=0\) if
\(i\geq 2\).
For every commutative with \(1\) over \(\bbZ_{(p)}\), the image of
\(\epsilon_p\) in \(\Cart{k}\) under the canonical ring homomorphism
\(\Cart{\bbZ_{(p)}}\to \Cart{k}\) is also denoted by \(\epsilon_p\).
\end{defn}
\medbreak


\subsection{\bf Exercise.} \label{exer:artin-hasse}
%
%
%
%
Let \(k\) be a \(\bbZ_{(p)}\)-algebra, and let 
\(\left(a_m\right)_{m\geq 0}\) be a sequence in \(k\).
Prove the equality
\begin{eqnarray*}
\epsilon_p\left(\omega\left(\prod_{m\geq 1}\,(1-a_m\,T^m) \right)\right)
=\epsilon_p\left(\omega\left(\prod_{n\geq 0}\,(1-a_{p^n}\,T^{p^n})\right)
\right)
\\
=\omega\left(\prod_{n\geq 0}\,E(a_{p^n}\,T^{p^n})\right)\,,
\end{eqnarray*}
in \(\tilde{W}(k)\),
where 
\[
E(X)=\prod_{(n,p)=1}\,(1-X^n)^{\frac{\mu(n)}{n}}
=\exp\left(-\sum_{n\geq 0}\,\frac{X^{p^n}}{p^n} \right)
\in 1+X\bbZ_{(p)}[[X]]
\]
is the inverse of the classical Artin-Hasse exponential.

\medbreak

\begin{prop} \Bl \label{prop_idempotent}
Let \(k\) be a commutative \(\bbZ_{(p)}\)-algebra with \(1\).
The following equalities hold in \(\Cart{k}\).
\begin{itemize}
\rmitem[(i)] \({\epsilon_p}^2 =\epsilon_p\).

\rmitem[(ii)] \(\displaystyle{\sum_{\tworows{(n,p)=1}{n\geq 1}}\,
\epsilon_{p,n} =1}\).

\rmitem[(iii)] \(\epsilon_p V_n=0\), \(F_n \epsilon_p=0\)
for all \(n\) with \({\rm gcd}(n,p)=1\).

\rmitem[(iv)] \({\epsilon_{p,n}}^2 =\epsilon_{p,n}\) for all \(n\geq 1\) with
\({\rm gcd}(n,p)=1\).

\rmitem[(v)] \(\epsilon_{p,n}\,\epsilon_{p,m}=0\)
for all \(m\neq n\) with \({\rm gcd}(mn, p)=1\).

\rmitem[(vi)] \([c]\,\epsilon_p =\epsilon_p\,[c]\) 
and \([c]\,\epsilon_{p,n} =\epsilon_{p,n}\,[c]\)
for all \(c\in k\) and all \(n\) with \({\rm gcd}(n,p)=1\).

\rmitem[(vii)] \(F_p\epsilon_{p,n}=\epsilon_{p,n} F_p\),
\(V_p\epsilon_{p,n}=\epsilon_{p,n} V_p\)
for all \(n\) with \({\rm gcd}(n,p)=1\).
\end{itemize}
\end{prop}
\medbreak

\begin{defn}\label{defn_local_cartring}
Let \(k\) be a commutative ring with \(1\) over \(\bbZ_{(p)}\). 
\begin{itemize}
\rmitem[(i)]
Denote by \(\Cartp{k}\) the subring \(\epsilon_p \Cart{k} \epsilon_p\)
of \(\Cart{k}\).  Note that \(\epsilon_p\) is the unit element
of \(\Cartp{k}\).

\rmitem[(ii)] Define elements \(F,V\in \Cartp{k}\) by
\[
F=\epsilon_p F_p = F_p \epsilon_p = \epsilon_p F_p \epsilon_p\,,
\quad
V=\epsilon_p V_p = V_p \epsilon_p = \epsilon_p V_p \epsilon_p\,.
\]

\rmitem[(iii)] For every element \(c\in k\), denote by
\(\langle c\rangle\) the element \
\(\epsilon_p [c] \epsilon_p =\epsilon_p [c] = [c] \epsilon_p \in \Cartp{k}\).
\end{itemize}
\end{defn}
\medbreak

\subsection{\bf Exercise.}\label{exer_local_commutation}
Prove the following identities in \(\Cartp{k}\).
\begin{itemize}
\rmitem[(1)] \(F\,\langle a\rangle = \langle a^p\rangle\,F\)
for all \(a\in k\).

\rmitem[(2)] \(\langle a \rangle V= V\,\langle a^p\rangle\)
for all \(a\in k\).

\rmitem[(3)] \(\langle a\rangle \,\langle b\rangle=
\langle ab\rangle\) for all \(a,b\in k\).

\rmitem[(4)] \(FV=p\).

\rmitem[(5)] \(VF=p\) if and only if \(p=0\) in \(k\).

\rmitem[(6)] Every prime number \(\ell\neq p\) is invertible
in \(\Cartp{k}\).  The prime number \(p\) is invertible in
\(\Cartp{k}\) if and only if \(p\) is invertible in \(k\).

\rmitem[(7)] \(V^m\langle a \rangle F^m\,V^n\langle b \rangle F^n
=p^r\,V^{m+n-r}\langle a^{p^{n-r}} b^{p^{m-r}}\rangle F^{m+n-r}\)
for all \(a,b\in k\) and all \(m,n\in \bbN\), where \(r=\min\{m,n\}\).
\end{itemize}

\begin{defn}\label{defn_local_lambda}
Let \(k\) be a commutative \(\bbZ_{(p)}\)-algebra with \(1\).
Denote by \(\Lambda_p\) the image of \(\epsilon_p\) in 
\(\Lambda\). In other words, \(\Lambda_p\) is the functor
from \(\mathfrak{Nilp}_k\) to \(\mathfrak{Ab}\) such that
\[
\Lambda_p(N)=\Lambda(N)\cdot \epsilon_p
\]
for any nilpotent \(k\)-algebra \(N\).
\end{defn}
\smallbreak

\begin{defn}
(1) Denote by \(W_p\) the image of \(\epsilon_p\), i.e.\ 
\(W_p(R):=\epsilon_p(\widetilde W(R))\) for every 
\(\bbZ_{(p)}\)-algebra \(R\).
Equivalently, \(W_p(R)\) is the intersection of the kernels
\({\rm Ker}(F_{\ell})\) of the operators \(F_{\ell}\)
on \(\widetilde{W}(R)\), where \(\ell\) runs through all
prime numbers different from \(p\).
\smallbreak

\noindent (2) Denote the element 
\[\omega(\prod_{n=0}^{\infty}\,E(c_n\,T^{p^n}))\in W_p(R)\]
by \(\omega_p(\underline{c})\). 
\smallbreak

\noindent (3) The endomorphism \(V_p, F_p\) of the group scheme
\(\widetilde{W}\) induces
endomorphisms of the group scheme \(W_p\),
denoted by \(V\) and \(F\) respectively.

\end{defn}
\medbreak

\noindent{\bf Remark.}\enspace 
The functor \(W_p\) has a natural structure
as a ring-valued functor induced from that of \(\widetilde W\);
it is represented by the
scheme \(\Spec\,\bbZ_{(p)}[y_0,y_1,y_2,\ldots, y_n,\ldots]\)
such that the element  \(\omega_p(\underline{c})\) has coordinates
\(\underline{c}=(c_0, c_1, c_2,\ldots)\).
\smallbreak

\subsection{\bf Exercise.}
Let \(k\) be a commutative \(\bbZ_{(p)}\)-algebra with \(1\).
Let \(E(T)\in \bbZ_{(p)}[[T]]\) be the inverse of the Artin-Hasse
exponential as in Exer.\ \ref{exer:artin-hasse}.
\begin{itemize}
\rmitem[(i)] Prove that for any nilpotent \(k\)-algebra \(N\),
every element of \(\Lambda_p(N)\) has a unique expression as a finite 
product
\[
\prod_{i=0}^m\,E(u_i\,t^{p^i})\,
\]
for some \(m\in \bbN\), and \(u_i\in N\)
for \(i=0,1,\ldots,m\).

\rmitem[(ii)] Prove that \(\Lambda_p\) is a smooth commutative formal
group over \(k\).

\rmitem[(iii)] Prove that every element of \(W_p(k)\) can 
be uniquely expressed as an infinite product
\[\omega(\prod_{n=0}^{\infty}\,E(c_n\,T^{p^n}))\in W_p(R)=:
\omega_p(\underline{c})\,.\]

\rmitem[(iv)] Show that the map 
from  \(W_p(k)\) to the product ring \(\,\prod_0^{\infty}k\,\)
defined by
\[
\omega_p(\underline{c})\ \ 
\longmapsto
\big(\,w_n(\underline{c})\,\bigr)_{n\geq 0}\qquad
{\rm where}\quad
w_n(\underline{c}):=\sum_{i=0}^{n}\,p^{n-i}\,c_{n-i}^{p^i}\,,
\]
is a ring homomorphism.
\end{itemize} 

\medbreak

\begin{prop}\label{prop_elt_local}
\begin{itemize}
\rmitem[(i)] The local Cartier ring \(\Cartp{k}\)
is complete with respect to the decreasing sequence of right ideals
\(V^i\Cartp{k}\).

\rmitem[(ii)] Every element of \(\Cartp{k}\) can be expressed 
as a convergent sum in the form
\[
\sum_{m,n\geq 0}\, V^m \langle a_{mn}\rangle F^n\,,
\quad a_{mn}\in k\,, \forall\,m\ \exists\, C_m>0\ {\rm s.t.\ }
a_{mn}=0\ {\rm if\ } n\geq C_m\,
\]
in a unique way.

\rmitem[(iii)] The set of all elements of \(\Cartp{k}\) 
which can be represented
as a convergent sum of the form
\[
\sum_{m\geq 0}\, V^m \langle a_{m}\rangle F^m\,,
\quad a_{m}\in k
\]
is a subring of \(\Cartp{k}\).
The map
\[
w_p(\underline{a})\ \mapsto 
\sum_{m\geq 0}\, V^m \langle a_{m}\rangle F^m\,\quad
\underline{a}=(a_0,a_1,a_2,\ldots),\
a_i\in k\ \forall\,i\geq 0
\]
establishes an isomorphism from the ring of \(p\)-adic Witt vectors
\(W_p(k)\) to the above subring of \(\Cartp{k}\).
\end{itemize}
\end{prop}
\medbreak

\subsection{\bf Exercise.}
Prove that \(\Cartp{k}\) is naturally
isomorphic to \(\End(\Lambda_{p})^{\rm op}\), the opposite
ring of the endomorphism ring of \(\End(\Lambda_{p})\).

\medbreak

\begin{defn}\label{defn_local_vred}
Let \(k\) be a commutative \(\bbZ_{(p)}\)-algebra.
\begin{itemize}
\rmitem[(i)] A \(V\)\!-{\emph{reduced}} left \(\Cartp{k}\)-module \(M\)
is a left \(\Cartp{k}\)-module such that 
the map \(V:M\to M\) is injective and the canonical map
\(\displaystyle M\to \varprojlim_n\,(M/V^n M)\) is an isomorphism.

\rmitem[(ii)] A \(V\)\!-{reduced} left \(\Cartp{k}\)-module \(M\) is
\(V\)\!-{\emph{flat}} if \(M/VM\) is a flat \(k\)-module.
\end{itemize}
\end{defn}
\medbreak

\begin{thm}\label{thm_modules}
Let \(k\) be a commutative \(\bbZ_{(p)}\)-algebra with \(1\).
\begin{itemize}
\rmitem[(i)]
There is an equivalence of categories
between the category of \(V\)\!-reduced  left \(\Cart{k}\)-modules 
and the category of \(V\)\!-reduced left \(\Cartp{k}\)-modules, defined
as follows.
\[
\xymatrix@R=2pt{
\{\mbox{ V-}{\rm reduced\ left\ }
\Cart{k}\mbox{-}{\rm mod}\}\hspace{0.2cm}
\ar@1{->}[r]^{\sim} & 
\hspace{0.2cm}\{\mbox{ V-}{\rm reduced\ left\ }
\Cartp{k}\mbox{-}{\rm mod}\}
\\
\hspace{1.5cm}M\hspace{1.5cm} \ar@{|->}[r] & \hspace{1cm} 
\epsilon_p M\\
\Cart{k}\epsilon_p\widehat{\otimes}_{\Cartp{k}}M_p\hspace{0.7cm}\hspace{0.5cm}
&\hspace{2.3cm}M_p\hspace{1.2cm}\ar@{|->}[l]
}
\]

\rmitem[(ii)]
Let \(M\) be a  \(V\)\!-reduced left \(\Cart{k}\)-module \(M\), and 
let \(M_p\) be the \(V\)\!-reduced let \(\Cartp{k}\)-module \(M_p\) 
attached to \(M\) as  in {\rm (i)} above.
There is a canonical isomorphism \(M/{\rm Fil}^2 M\cong M_p/V M_p\).
In particular \(M\) is \(V\)\!-flat if and only if 
\(M_p\) is \(V\)\!-flat.
Similarly \(M\) is a finitely generated \(\Cart{k}\)-module 
if and only if \(M_p\) is
a finitely generated \(\Cartp{k}\)-module.
\end{itemize}
\end{thm}
\bigbreak

\begin{thm}\label{thm:main_local}
Let \(k\) be a commutative \(\bbZ_{(p)}\)-algebra with \(1\). 
There is
a canonical equivalence of categories, between the category
of smooth commutative formal groups over \(k\) as defined in \ref{def_grp}
and the category of \(V\)\!-flat \(V\)\!-reduced left
\(\Cartp{k}\)-modules, defined as follows. 
\[
\xymatrix@R=2pt{
\{{\rm smooth\ formal\ groups\ over\ } k\}\hspace{0.2cm}
\ar@1{->}[r]^{\sim} & 
\hspace{0.2cm}\{\mbox{V-}{\rm flat}\mbox{ V-}{\rm reduced\ left\ } 
\Cartp{k}\mbox{-}{\rm mod}\}\\
\hspace{1.5cm}G\hspace{1.5cm} \ar@{|->}[r] & \hspace{1cm} 
{\rm M}_{p}(G)=\epsilon_p\Hom(\Lambda,G)\\
\Lambda_p{\otimes}_{\Cartp{k}}M\hspace{0.7cm}\hspace{0.5cm}
&\hspace{2.3cm}M\hspace{1.2cm}\ar@{|->}[l]
}
\]
\end{thm}

\newpage
\n
{\large\bf\scshape Dieudonn\'e modules.}
\medbreak

In the rest of this section, \(K\) stands for a perfect field of
characteristic \(p>0\). We have \(FV=VF=p\) in \(\Cartp{K}\).
It is well-known that the ring of \(p\)-adic Witt vectors \(W(K)\)
is a complete discrete valuation ring with residue field \(K\),
whose maximal ideal is generated by \(p\).
Denote by \(\sigma:W(K)\to W(K)\) the Teichm\"uller lift
of the automorphism \(\,x\mapsto x^p\,\) of \(K\).
With the Witt coordinates we have \(\,\sigma: (c_0,c_1,c_2,\ldots)
\mapsto (c_0^p,c_1^p,c_2^p,\ldots)\).
Denote by \(L=B(K)\) the field of fractions of \(W(K)\).
\medbreak

\begin{defn} Denote by \(R_K\) the (non-commutative) ring generated by 
\(W(K)\), \(F\) and \(V\), subject to the following relations
\[
F\cdot V=V\cdot F=p,\quad F\cdot x\, =\, {}^{\sigma}x\cdot F,\quad 
x \cdot V=V\cdot {}^{\sigma}x
\qquad \forall x\in W(K)\,.
\]
\end{defn}

\noindent{\bf Remark.}\enspace
There is a natural embedding \(\,R_K\hookrightarrow \Cartp{K}\);
we use it to identify \(\,R_K\) as a dense subring of the Cartier ring
\(\Cartp{K}\). For every continuous left \(\Cartp{K}\)-module \(M\),
the \(\Cartp{K}\)-module structure on \(M\) is determined by the induced 
left \(R_K\)-module structure on \(M\).
\medbreak

\noindent{\bf Lemma/Exercise.}\enspace
(i) The ring \(R_K\) is naturally identified with the ring 
\[
W(K)[V,F]
:=\left(\bigoplus_{i<0} p^{-i} V^i W(K)\right)
\bigoplus \left(\bigoplus_{i\geq 0} V^i W(K)\right)
\,,
\]
i.e.\ elements of \(W(K)[V,F]\) are sums of the form
\(\sum_{i\in\bbZ} a_i V^i\), where \(a_i\in L\) for all \(i\in\bbZ\),
\({\rm ord}_p(a_i)\geq 
\max(0,-i)\ \ \forall i\in\bbZ\),
and \(a_i=0\) for all but finitely many \(i\)'s.
The commutation relation between \(W(K)\) and \(V^i\) is 
\[\,x\cdot V^i=\ V^i\cdot {^{\sigma^i}\! x}\quad {\rm for\ all\ } 
x\in W(K)\ {\rm and\ all\ } 
i\in\bbZ\,.\]
\smallbreak

\noindent (ii) The ring \(\Cartp{K}\) is naturally
identified with the set \(\,W(K)[[V,F\rangle\rangle\,\), 
consisting of all non-commutative formal power
series of the form
\(\,
\sum_{i\in\bbZ} a_i V^i\,
\) 
such that \(a_i\in L\ \forall i\in\bbZ\),  
\({\rm ord}_p(a_i)\geq  \max(0,-i)\  \forall i\in\bbZ\), and
\({\rm ord}_p(a_i)+i\to\infty\ {\rm as}\ |i|\to\infty\).
\smallbreak

\noindent (iii) Check that ring structure on \(W(K)[V,F]\)
extends to \(\,W(K)[[V,F\rangle\rangle\,\) by continuity.
In other words, the inclusion 
\(\,W(K)[V,F]\hookrightarrow W(K)[[V,F\rangle\rangle\,\)
is a ring homomorphism, and \(W(K)[V,F]\) is 
dense in  \(\,W(K)[[V,F\rangle\rangle\,\) with respect to the 
\(V\)-adic topology on \(\,W(K)[[V,F\rangle\rangle\cong \Cartp{K}\,\).
The latter topology on \(\,W(K)[[V,F\rangle\rangle\,\) is equivalent
to the topology given by the discrete valuation \(v\) on 
\(\,W(K)[[V,F\rangle\rangle\,\) defined by
\[
v\left(\sum_{i\in \bbZ}\,a_i V^i\right)=
{\rm Min}\left\{{\rm ord}_p(a_i)+i\,\vert\,i\in\bbZ\right\}\,.
\]
\medbreak

\begin{defn}
(1) A \emph{Dieudonn\'e module} is a left \(R_K\)-module \(M\)
such that \(M\) is a free \(W(K)\)-module of finite rank.
\smallbreak

\noindent (2) Let \(M\) be a Dieudonn\'e module over \(K\).
Define the \(\alpha\)-rank of \(M\) to be the
natural number \(\,a(M)=\dim_K(M/(VM+FM)\).
\end{defn}
\medbreak

\begin{defn} 
\begin{itemize} 
\rmitem[(i)] For any natural number \(n\geq 1\) and 
any scheme \(S\), denote by \(\underline{(\bbZ/n\bbZ)}_S\)
the constant group scheme over \(S\) attached to the finite
group \(\bbZ/n\bbZ\). The scheme underlying \(\underline{(\bbZ/n\bbZ)}_S\)
is the disjoint union of \(n\) copies of \(S\), indexed
by the finite group \(\bbZ/n\bbZ\), see \ref{grsch}.

\rmitem[(ii)]  For any natural number \(n\geq 1\) and 
any scheme \(S\), denote by \(\mu_{n,S}\)
the kernel of \([n]:{\bbG_m}_{/S}\to {\bbG_m}_{/S}\).
The group scheme \(\mu_{n,S}\) is finite and locally free over
\(S\) of rank \(n\); it is the Cartier dual of \(\underline{(\bbZ/n\bbZ)}_S\).

\rmitem[(iii)] For any field \(K\supset \Fp\), define
a finite group scheme
\(\alpha_p\) over \(K\) to be
the kernel of the endomorphism
\[
{\rm Fr}_p: {\bbG_a}_{/K}=\Spec(K[X])\to {\bbG_a}_{/K}=\Spec(K[X])
\]
of \(\bbG_a\) over \(K\)
defined by the \(K\)-homomorphism from the \(K\)-algebra
\(K[X]\) to itself which sends \(X\) to \(X^p\).
We have \(\alpha_p=\Spec(K[X]/(X^p))\) as a scheme.
The comultiplication is induced by \(\,X\mapsto X\otimes X\).
\end{itemize}
\end{defn}

\begin{prop}\label{btdecomp}\Bl\ 
Let \(X\) be a Barsotti-Tate group over a perfect field \(K\supset \Fp\).
Then there exists a canonical splitting
\[
X\cong X_{\rm mult}\times_{\Spec(K)} X_{\rm ll}\times_{\Spec(K)}X_{\rm et}\,
\]
where \(X_{\rm et}\) is the maximal \'etale quotient of \(X\),
\(X_{\rm mult}\) is the maximal multiplicative Barsotti-Tate
subgroup of \(X\),
and \(X_{\rm ll}\) is a Barsotti-Tate group with no non-trivial
\'etale quotient nor non-trivial multiplicative Barsotti-Tate subgroup.
\end{prop}
\smallbreak

\noindent{\bf Remark.}\enspace
The analogous statement for finite group schemes over
\(K\) can be found in \cite[Chap.\ 1]{Manin}, from which
\ref{btdecomp} follows.  See also \cite{Demazure}, \cite{D.G}.
\bigbreak

\begin{defn}
Let \(m,n\) be non-negative integers such that \({\rm gcd}(m,n)=1\).
Let \(k\supset \Fp\) be an algebraically closed field.
Let \(G_{m,n}\) be the Barsotti-Tate group whose
Dieudonn\'e module is 
\[
{\bbD}(G_{m,n})=R_K/R_K\cdot(V^n-F^m)\,.
\]
\end{defn}
\bigbreak

\begin{thm}\Bl\label{main_dieu}
\begin{itemize}
\rmitem[(1)] There is an equivalence of categories between
the category of Barsotti-Tate groups over \(K\)
and the category of Dieudonn\'e modules over \(R_K\).
Denote by \({\bbD}(X)\) the covariant Dieudonn\'e module
attached to a Barsotti-Tate group over \(K\).
This equivalence is compatible with direct product and exactness,
i.e.\ short exact sequences correspond under the above equivalence
of categories.

\rmitem[(2)] Let \(X\) be a Barsotti-Tate group over \(K\)
such that \(X\) is a \(p\)-divisible formal group in the
sense that the maximal \'etale quotient of \(X\) is trivial.
Denote by \(\hat{X}\) the formal group attached to \(X\),
i.e.\ \(\hat{X}\) is the formal completion of
\(X\) along the zero section of \(X\).
Then there is a canonical isomorphism 
\({\bbD}(X)\xrightarrow{\sim} {\rm M}_p(\hat{X})\)
between the Dieudonn\'e module of \(X\)
and the Cartier module of \(\hat{X}\).

\rmitem[(3)] Let \(X\) be a Barsotti-Tate group over \(K\),
and let \({\bbD}(X)\) be the covariant Dieudonn\'e module
of \(X\).  Then \({\rm ht}(X)={\rm rank}_{W(K)}({\bbD}(X))\),
and we have a functorial isomorphism
\({\rm Lie}(X)\cong {\bbD}(X)/V\cdot{\bbD}(X)\).

\rmitem[(4)] Let \(X^t\) be the Serre-dual of the Barsotti-Tate
group of \(X\).  Then the Dieudonn\'e module \({\bbD}(X^t)\)
can be described in terms of \({\bbD}(X)\) as follows.
The underlying \(W(K)\)-module is the linear dual
\(\bbD(X)^{\vee}:={\rm Hom}_{W(K)}({\bbD}(X),W(K))\) of \({\bbD}(X)\).
The actions of \(V\) and \(F\) on \(\bbD(X)^{\vee}\)
are defined as follows.
\[
(V\cdot h)(m)\,=\,^{\sigma^{-1}}(h(Fm))\,,
\quad
(F\cdot h)(m)\,=\,^{\sigma}(h(Vm))
\]
for all \(h\in \bbD(X)^{\vee}={\rm Hom}_{W(K)}({\bbD}(X),W(K))\) and
all \(m\in {\bbD}(X)\).

\rmitem[(5)] A Barsotti-Tate group \(X\) over \(K\) is \'etale
 if and only if 
\(V:{\bbD}(X)\to {\bbD}(X)\) is bijective,
or equivalently, 
\(F:{\bbD}(X)\to {\bbD}(X)\)
is divisible by \(p\).
A Barsotti-Tate group \(X\) over \(K\) is multiplicative
 if and only if 
\(V:{\bbD}(X)\to {\bbD}(X)\) is divisible by \(p\),
or equivalently, 
\(F:{\bbD}(X)\to {\bbD}(X)\)
is bijective.
A Barsotti-Tate group \(X\) over \(K\) has no non-trivial 
\'etale quotient nor non-trivial multiplicative Barsotti-Tate 
subgroup if and only if
both \(F\) and \(V\) are topologically nilpotent on 
\(\bbD(X)\).

\end{itemize}
\end{thm}
\smallbreak

\noindent{\bf Remark.} \enspace
See \cite{Oda.FO} for Thm.\ \ref{main_dieu}.  
The Dieunonn\'e module \(\bbD(X)\) attached to a \(p\)-divisible group 
over \(K\) can also be defined in terms of the covariant Dieudonn\'e
crystal attached to \(X\) described in \ref{GM}.
In short, \(\bbD(X)\) ``is'' \(\bbD(X/W(K))_{W(K)}\), the
limit of the ``values'' of the Dieudonn\'e crystal 
at the divided power structures \((W(K)/p^m W(K), pW(K)/p^mW(K),\gamma)\)
as \(m\to\infty\), where \((W(K)/p^m W(K), pW(K)/p^mW(K),\gamma)\) is the
reduction modulo \(p\) of the natural DP-structure on
\((W(K),pW(K))\).  However there is a technical complication if \(p=2\),
because the natural DP-structure on \((W(K)/p^mW(K), pW(K)/p^mW(K))\) 
is not nilpotent.  In the rest of this article we will not use
the crystalline interpretation of \(\bbD(X)\).
\smallbreak

\begin{prop} \Bl\  Let \(X\) be a Barsotti-Tate group over \(K\).
We have a natural isomorphism
\[
\Hom_K(\alpha_p, X[p])\cong 
\Hom_{W(K)}\left(
{\bbD}(X)^{\vee}/(V{\bbD}(X)^{\vee}+F{\bbD}(X)^{\vee}), B(K)/W(K)\right)\,,
\]
where \(B(K)={\rm frac}(W(K))\) is the fraction field of \(W(K)\).
In particular we have 
\[\dim_K(\Hom_K(\alpha_p, X[p]))=a({\bbD}(X))\,.\]
The natural number \(a({\bbD}(X))\) of a BT-group \(X\) over \(K\)
is zero if and only if \(X\) is an extension of an 
\'etale BT-group by a multiplicative BT-group.
\end{prop}
\medbreak

\noindent{\bf Definition.}\enspace 
The \(a\)-number of a Barsotti-Tate group \(X\) over
\(K\) is the natural number 
\[a(X):=\dim_K(\Hom_K(\alpha_p, X[p]))=a({\bbD}(X))\,.\]
\bigbreak

\subsection{\bf Exercise.} (i) Prove that \({\rm ht}(G_{m,n})=m+n\).
\smallbreak

\noindent (ii) Prove that \(\dim(G_{m,n})=m\).
\smallbreak

\noindent (iii) Show that \(G_{0,1}\) is isomorphic
to the \'etale Barsotti-Tate group \(\Qp/\Zp\),
and \(G_{1,0}\) is isomorphic to the multiplicative
Barsotti-Tate group 
\(\mu_{\infty}=\bbG_m[p^{\infty}]\).

\noindent (iv) Show that \(\End(G_{m,n})\otimes_{\Zp}\Qp\)
is a central division algebra over \(\Qp\) of dimension
\((m+n)^2\), and compute the Brauer invariant of this
central division algebra.

\begin{thm}\Bl \label{DManin}
Let \(k\supset \Fp\) be an algebraically closed field.
Let \(X\) be a simple Barsotti-Tate group over \(k\), i.e.
\(X\) has no non-trivial quotient Barsotti-Tate groups.
Then \(X\) is isogenous to \(G_{m,n}\) for a
uniquely determined pair of natural numbers \(m,n\) with
\({\rm gcd}(m,n)=1\), i.e.\
there exists a surjective homomorphism
\(X\to G_{m,n}\) with finite kernel.
\end{thm}
\bigbreak

\begin{defn}
{\rm (i)} The slope of \(G_{m,n}\) is \(m/(m+n)\)
with multiplicity \(m+n\). The Newton polygon of
\(G_{m,n}\) is the line segment on the plane from 
\((0,0)\) to \(m+n,m\). The slope sequence of \(G_{m,n}\)
is the finite sequence \((m/(m+n),\ldots,m/(m+n))\) 
with \(m+n\) entries.
\smallbreak

\noindent {\rm (ii)} Let \(X\) be a Barsotti-Tate group
over a field \(K\supset \Fp\), and let \(k\) be an 
algebraically closed field containing \(K\).
Suppose that \(X\) is isogenous to
\[
G_{m_1,n_1}\times_{\Spec(k)}\cdots\times_{\Spec(k)}G_{m_r,n_r}\,
\]
\({\rm gcd}(m_i,n_i)=1\) for \(i=1,\ldots,r\),
and \(m_i/(m_i+n_i)\leq m_{i+1}/(m_{i+1}+n_{i+1})\)
for \(i=1,\ldots,r-1\).
Then the Newton polygon of \(X\) is defined by the 
data \(\sum_{i=1}^r (m_i, n_i)\).
Its slope sequence is the concatenation of the slope sequence
for \(G_{m_1,n_1},\ldots,G_{m_r,n_r}\).
\end{defn}
\bigbreak

\noindent{\bf Example.}\enspace
A Barsotti-Tate group \(X\) over \(K\) is \'etale (resp.\ multiplicative)
if and only if all of its slopes are equal to \(0\) (resp.\ \(1\)).
\bigbreak

\subsection{\bf Exercise.}
Suppose that \(X\) is a Barsotti-Tate group over \(K\)
such that \(X\) is isogenous to \(G_{1,n}\) (resp.\ \(G_{m,1}\)).
Show that \(X\) is isomorphic to \(G_{1,n}\) (resp.\ \(G_{m,1}\)).

\subsection{\bf Exercise.}
Show that there are infinitely many non-isomorphic Barsotti-Tate
groups with slope sequence \((1/2,1/2,1/2,1/2)\)
(resp.\ \((1/3,1/3,1/3,2/3,2/3,2/3)\).

\subsection{\bf Exercise.} 
Determine all Newton polygons attached to
a Barsotti-Tate group of height \(6\), and the symmetric Newton polygons
among them.
\smallbreak

Recall that the set of all Newton polygons is a partially ordered;
\(\zeta_1\prec \zeta_2\) if and only if \(\zeta_1,\zeta_2\) have
the same end points, and \(\zeta_2\) lies below \(\zeta_1\).
Show that this poset is \emph{ranked}, i.e.\ any two
maximal chains between two elements of this poset have the same length.
\medbreak

\begin{prop}\Bl\ \label{dieumod_fingrp}
\begin{itemize}
\rmitem[(i)] There is an equivalence of categories between 
the category of finite group schemes over the perfect base field \(K\)
and the category of left \(R_K\)-modules which are \(W(K)\)-modules  of 
finite length. 
Denote by \(\bbD(G)\) the left \(R_K\)-module attached to a
finite group scheme \(G\) over \(K\).

\rmitem[(ii)] Suppose that 
\(\,0\to G\to X\xrightarrow{\beta} Y\to 0\,\)
is a short exact sequence, where \(G\) is a finite group scheme over 
\(K\), and \(\beta:X\to Y\) is an isogeny between BT-groups over \(K\).
Then we have a natural isomorphism
\[
\bbD(G)\xrightarrow{\sim} {\rm Ker}
\left(\bbD(X)\otimes_{W(K)}B(K)/W(K)\xrightarrow{\beta}
\bbD(Y)\otimes_{W(K)}B(K)/W(K)
\right)
\]
of left \(R_K\)-modules.
\end{itemize}
\end{prop}

\noindent{\bf Remark.}\enspace
(i) We say that \(\bbD(G)\) is the \mbox{covariant Dieudonn\'e module}
of \(G\), abusing the terminology, because \(\bbD(G)\) is not
a free \(W(K)\)-module.

(ii) Prop.\ \ref{dieumod_fingrp} is a covariant version of the
classical contravariant Dieudonn\'e theory in \cite{Demazure}
and \cite{Manin}.  See also \cite{FO-reduced}.
\medbreak

\subsection{\bf ``Dieudonn\'e modules'' over non-perfect fields?}
\label{nonperfect}  This is a difficult topic. 
However, in one special case statements and results are easy.
\medbreak

{\bf $p$-Lie algebras.} Basic reference \cite{D.G}. 
We will need this theory only in the commutative case. 
For more general statements see \cite{D.G}, II.7.\\
Let $K \supset \FF_p$ be a field. 
A commutative finite group scheme of {\it height one} over $K$ 
is a finite commutative group scheme $N$ over $K$ such that 
$(F: N \to N^{(p)}) = 0$ is the zero map. Denote the category 
of such objects by ${\rm GF}_K$.\\
A commutative finite dimensional $p$-Lie algebra $M$ over $K$ 
is a pair $(M,g)$, where $M$ is a finite dimensional vector space
 over $K$, and $g: M \to M$ is a homomorphism of additive groups 
with the property  
$$
g(b{\cdot}x) \quad=\quad b^p {\cdot} g(x).
$$
Denote the category of such objects by ${\rm Liep}_K$.
\medbreak

\begin{thm} \  \Bl 
There is an equivalence of categories 
$$\cD_K: {\rm GF}_K \stackrel{\sim}{\longrightarrow} {\rm Liep}_K.$$ 
This equivalence commutes with base change. If $K$ is a perfect field 
this functor coincides with the Dieudonn\'e module functor: 
$\cD_K = {\rm M}_p$, with ${\rm M}_p(V) = g$.
\end{thm}

See \cite{D.G}, II.7.4

\subsection{}{\bf Exercise.} 
{\it Classify all commutative group schemes of rank $p$ 
over $k$, an algebraically closed field of characteristic $p$. 

Classify all commutative group schemes of rank $p$ 
over a perfect field $K \supset \FF_p$.}

\subsection{} 
For group schemes in characteristic $p$ we have seen the 
Frobenius homomorphism, and for  commutative group schemes the Verschiebung. 
Also for Dieudonn\'e modules such homomorphisms are studied. 
However some care has to be taken. In the {\it covariant} 
Dieudonn\'e module theory the Frobenius on group schemes is related to the 
Verschiebung on the related module, and the same phenomenon for the 
Verschiebung on group schemes giving the Frobenius on modules; 
for details see \cite{FO-EO}, 15.3. In case confusion is possible 
we can write $F$ for the Frobenius on group schemes and 
$\cV = \DD(F)$ for the Verschiebung on modules, 
and $V$ on group schemes and $\cF = \DD(V)$ on modules. See \ref{FV}.


\section{Cayley-Hamilton: a conjecture by Manin and the 
weak Grothendieck conjecture}\label{5}
\n
Main reference: \cite{FO-CH}. 

\subsection{}
For a matrix $\cF$  over a {\it commutative} integral domain $R$ 
we have the Cayley-Hamilton theorem: {\it let 
$$\det(\cF - T{\cdot}I) =: g \in R[T]$$ 
be the characteristic polynomial of this matrix; 
then $g(\cF) = 0$, i.e.  ``a matrix is a zero of its 
own characteristic polynomial''.}

\subsection{}{\bf Exercise.} Show the classical Cayley-Hamilton theorem for 
a matrix over a commutative ring $R$: 
\begin{center}
{\it let $X$ be a  $n \times n$ matrix with entries in $R$; \\
let $g = {\rm Det}(X - T{\cdot}{\bf 1}_n) \in R[T]$; 
the matrix $g(X)$ is the zero matrix.} 
\end{center}
Here are some suggestions for a proof:\\
(a)  For any commutative ring $R$, and a $n \times n$ matrix $X$ 
with entries in $R$ there exists a ring homomorphism 
$\ZZ[t_{1,1}, \cdots , t_{i,j}, \cdots, t_{n,n}] \to R$ such that 
the matrix $(t) = ( t_{i,j} \mid 1 \leq i, j \leq n)$ is mapped to $X$. \\
(b) Construct $\ZZ[ t_{i,j} \mid   1 \leq i, j \leq n] \hookrightarrow \CC$. 
Show that the matrix $(t)$ considered over $\CC$ has mutually different 
eigenvalues. \\
(c) Conclude the theorem for $(t)$, over $\CC$ and for $(t)$ over  
$\ZZ[ t_{i,j} \mid   1 \leq i, j \leq n]$. Conclude the theorem for 
$X$ over $R$.\\
Here are suggestions for a different proof:\\ 
(1) Show it suffices to prove this for an algebraically closed field of 
characteristic zero.\\
(2) Show the  classical Cayley-Hamilton theorem holds for a matrix which 
is in diagonal form with all diagonal elements mutually different.\\
(3) Show that the set of all conjugates of matrices as in (2) is 
Zariski dense in Mat$(n \times n)$. Finish the proof.

\subsection{}
We will develop a very useful analog of this over the Dieudonn\'e ring. 
Note that over a non-commutative ring there is no reason 
any straight analog of Cayley-Hamilton should be true. However, 
given a specific element in a special situation, 
we construct an operator  $g(\cF)$ which  
annihilates that specific  element in the Dieudonn\'e module. 
Warning: In general  $g(\cF)$ does not  annihilate all elements of 
the Dieudonn\'e module.  

\subsection{\bf Notation.}\label{alpha}  Let $G$ be a group scheme 
over a field $K \supset \FF_p$. 
Consider $\alpha_p = \Ker(F: \GG_a \to \GG_a)$. 
Choose a perfect field $L$ containing $K$. 
Note that $\Hom(\alpha_p,G_L)$ is a right-module over 
$\End(\alpha_p \otimes_{\FF_p} L) = L$. We define
$$a(G) = \dim_L\left(\Hom(\alpha_p,G_L)\right).$$
{\bf Remarks.} For any field $L$ we write $\alpha_p$ instead of  
$\alpha_p \otimes_{\FF_p} L$ if no confusion is possible.

The group scheme $\alpha_{p,K}$ over a field $K$ corresponds under 
\ref{nonperfect} (in any case) or by Dieudonn\'e theory 
(in case $K$ is perfect) to the module $K^+$ with operators $\cF = 0$  
and $\cV = 0$.

If $K$ is not perfect it might happen that 
$\dim_K\left(\Hom(\alpha_p,G)\right) < 
\dim_L\left(\Hom(\alpha_p,G_L)\right)$, see the Exercise \ref{Exc} below. 

However if $L$ is perfect, and $L \subset L'$ is any field extension 
then $\dim_L\left(\Hom(\alpha_p,G_L)\right) = 
\dim_{L'}\left(\Hom(\alpha_p,G_{L'})\right)$; 
hence the definition of $a(G)$ is independent of the chosen 
perfect extension $L$.

\subsection{\bf Exercise.} Let $N$ be a finite group scheme over 
a perfect field $K$. Assume that $F$ and $V$ on $N$ are nilpotent on $N$, 
and suppose that $a(N) = 1$. {\it Show that the 
Dieudonn\'e module $\DD(N)$ 
is generated by one element over the Dieudonn\'e ring.}

Let $A$ be an abelian variety  over a perfect field $K$. 
Assume that the $p$-rank of $A$ is zero, and that $a(A) = 1$. 
{\it Show that the Dieudonn\'e module $\DD(A[p^{\infty}])$ 
is generated by one element over the Dieudonn\'e ring.}

\vn
{\bf Remark.} We will see that if $a(X_0) = 1$, then 
\begin{center}
$\cW_{\cN(X_0)}({\rm Def}(X_0))$ {\it is non-singular.}
\end{center} 
Let $(A,\lambda)$ be a principally polarized abelian variety, 
$\xi= \cN(A)$. 
The Newton polygon stratum $\cW_{\xi}(\cA_{g,1,n})$ will be shown to be 
regular at the point $(A,\lambda)$ (here we work with a fine moduli 
scheme: assume $n \geq 3$).

We see that we can a priori consider a set of points in which we know 
the Newton polygon stratum is non-singular. That is the main result 
of this section. Then, in Section \ref{7} we show that such points 
are dense in both cases considered, $p$-divisible groups and 
{\it principally} polarized abelian varieties.

\subsection{\bf Remark/Exercise.}\label{Exc} We give an example 
where $\dim_K\left(\Hom(\alpha_p,G)\right) < a(G)$; 
we see that the condition ``$L$ is perfect'' is necessary 
in \ref{alpha} .\\
(1) Let $K$ be a non-perfect field, with $b \in K$ and 
$\sqrt[p]{b} \not\in K$. Let $(M,g)$ be the commutative finite 
dimensional $p$-Lie algebra defined by:
$$
M = K{\cdot}x \oplus K{\cdot}y \oplus K{\cdot}z, 
\quad g(x) = bz, \ g(y) = z, \  g(z) = 0.
$$
Let $N$ be the finite group scheme of height one defined by this 
$p$-Lie algebra, i.e. such that $\cD_K(N) = (M,g)$, 
see \ref{nonperfect}. {\it Show}:
$$
\dim_K\left(\Hom(\alpha_p,N)\right) = 1, \quad a(N) = 2.
$$
\\
(2) Let $N_2 = W_2[F]$ be the kernel of $F: W_2 \to W_2$ over $\FF_p$; 
here $W_2$ is the 2-dimensional group scheme of Witt vectors of 
length $2$. In fact one can define $N_2$ by 
$\DD(N_2) = \FF_p{\cdot}r \oplus   \FF_p{\cdot}s$, \  
$\cV(r) = 0 = \cV(s) = \cF(s)$ and   $\cF(r) = s$. 
Let $L = K(\sqrt[p]{b})$. {\it Show that} 
$$
N \not\cong_K  (\alpha_p \oplus W_2[F]) \otimes K \quad\mbox{and}
\quad N \otimes L \cong_L (\alpha_p \oplus W_2[F]) \otimes L.
$$
{\bf Remark.} In \cite{Li.FO}, I.5  Definition (1.5.1) should 
be given over a perfect field $K$. We thank Chia-Fu Yu for 
drawing our attention to this flaw.

\subsection{} We fix integers $h \geq d \geq 0$, and we write $c := h-d$. 
We consider Newton polygons ending at $(h,d)$. 
For such a Newton polygon $\beta$ we write:
$$\diam(\beta) := \{(x,y) \in \ZZ \times \ZZ\mid  y < d, \ 
\ y < x, \ \ (x,y) \prec \beta\};$$
here we denote by $(x,y) \prec \beta$ the property ``$(x,y)$ 
is on or above $\beta$;"
we write

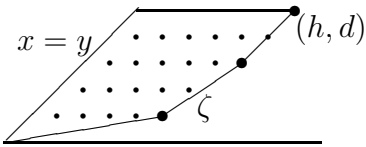
\begin{center}
\fbox{$\dm(\zeta) := \#(\diam(\zeta)).$}
\end{center}

\begin{picture}(1300,80)(0,0)
\put(-18,70){\bf Example:}
\put(0,0){\line(1,0){120}}
\put(0,0){\line(6,1){60}}
\put(60,10){\line(3,2){30}}
\put(90,30){\line(1,1){20}}

\put(60,10){\circle*{4}}
\put(90,30){\circle*{4}}
\put(110,50){\circle*{4}}

\put(1,1){\line(1,1){50}}
\put(50,50){\line(1,0){60}}

\put(20,10){\circle*{2}}
\put(30,10){\circle*{2}}
\put(40,10){\circle*{2}}
\put(50,10){\circle*{2}}

\put(30,20){\circle*{2}}
\put(40,20){\circle*{2}}
\put(50,20){\circle*{2}}
\put(60,20){\circle*{2}}
\put(70,20){\circle*{2}}

\put(40,30){\circle*{2}}
\put(50,30){\circle*{2}}
\put(60,30){\circle*{2}}
\put(70,30){\circle*{2}}
\put(80,30){\circle*{2}}

\put(50,40){\circle*{2}}
\put(60,40){\circle*{2}}
\put(70,40){\circle*{2}}
\put(80,40){\circle*{2}}
\put(90,40){\circle*{2}}
\put(100,40){\circle*{2}}

\put(5,35){$x=y$}
\put(110,40){$(h,d)$}
\put(73,10){$\zeta$}

\put(180,50){$\zeta =  2  \times(1,0) + (2,1) + (1,5)=$}
\put(180,30){$= 6 \times \frac{1}{6} + 3 \times \frac{2}{3} + 
2  \times \frac{1}{1}; \ \  h=11.$}
\put(180,10){Here $\dm(\zeta) =   \#(\diam(\zeta)) = 22$.}

\end{picture}


\vn
Note that for $\rho = d{\cdot}(1,0) + c{\cdot}(0,1)$ 
we have $\dm(\rho) = dc$.

\subsection{}\label{G}{\bf Theorem}   (Newton polygon-strata 
for $p$-divisible groups). 
{\it Suppose $a(X_0) \leq 1$. Write $\cD = {\rm Def}(X_0)$.  
For every $\beta \succ \gamma = \cN(X_0)$ 
we have: $\dm(\cW_{\beta}(D)) = \dm({\beta})$. 
The strata $\cW_{\beta}(D)$ are nested as given by the partial ordering on 
Newton polygon, i.e. 
$$\cW_{\beta}(D) \subset \cW_{\delta}(D) 
\quad\Longleftrightarrow \quad\diam(\beta) 
\subset \diam(\delta) \quad\Longleftrightarrow\quad  \beta \prec \delta.$$
Generically on $\cW_{\beta}(D)$ the fibers have Newton polygon 
equal to $\beta$.}\\
For the notion ``generic'' for a $p$-divisible group over a formal 
scheme, see \ref{algebr}.

\subsection{}  In fact, this can be visualized and made more precise 
as follows.  Choose variables $T_{r,s}$, with 
$1 \leq r \leq d = \dim(X_0)$, \ $1 \leq s \leq h = {\rm height}(X_0)$ 
and write these in a diagram
\[ \begin{array}{lllllllllll}
 & & & &   &  &  &0& \cdots& 0 &-1 \\
&&&&&&                      T_{d,h} &  \cdot  &  \cdots  &T_{1,h}& \\
&&&&  . &   \vdots &                     \vdots  &  \vdots &  . && \\
&&     T_{d,d+2}&   \cdots&  T_{i,d+2}&      \cdots&                    
T_{2,d+2}       &   T_{1,d+2} &&& \\
&T_{d,d+1} &   \cdots   &  T_{i,d+1} &  \cdots  & \cdots  
&                   T_{1,d+1}               &&&& \\
\end{array} \]

\vn
We show that 
$$\cD :={\rm Def}(X_0) = \Spf(k[[Z_{(x,y)}\mid (x,y) \in \diam]]), 
\quad\quad T_{r,s} = Z_{(s-r,s-1-d)}.$$
Moreover, for any  $\beta \succ \cN(G_0)$ we write 
$$R_{\beta}= \frac{k[[Z_{(x,y)}\mid (x,y) \in \diam ]]}{(  
 Z_{(x,y)} \ \ \forall (x,y) \not\in \diam(\beta))} 
\quad\cong\quad k[[Z_{(x,y)} \mid (x,y) \in \diam(\beta)]].$$
{\it Claim}: $$\left(\Spec(R_{\beta}) \subset \Spec(R)\right) 
\quad=\quad  \left( \cW_{\beta}(D) \subset \cD\right).$$

\vn
Clearly this claim proves the theorem. We will give a proof of the claim, 
and hence of this theorem  by using the theory of displays 
(see the paper by Messing in this volume), 
and by using the following tools. 

\subsection{\bf Definition.}\label{NF} We consider matrices which 
can appear as $F$-matrices associated with a display. 
Let $d, c \in \ZZ_{\geq 0}$, 
and $h = d + c.$ Let $W$ be a ring. 
We say that a display-matrix $(a_{i,j})$ of size $h \times h$ 
is in {\it normal form} form over $W$ if the $F$-matrix is 
of the following form:

\[\left(  \begin{array}{cccclclcccl}
0  & 0  & \cdots  & 0 & a_{1d}&  & pa_{1,d+1}  
&  \cdots &  \cdots& \cdots &   pa_{1,h}  \\
1  &  0 & \cdots  & 0  &  a_{2d}& & \cdots  &  &  pa_{i,j} &  & \cdots \\
0  & 1  & \cdots  &  0 &  a_{3d}& &   &  & 1 \leq i \leq d &  &  \\
\vdots  & \vdots  & \ddots  &  \ddots &  \vdots & &  &  & 
d \leq j \leq h  
&  &  \\
0 &  0 & \cdots   & 1  &  a_{dd} & & pa_{d,d+1}  & \cdots & \cdots & 
\cdots & pa_{d,h} \\
&&&&&&&&&&\\
 0 & \cdots  & \cdots  & 0  & 1 & & 0  & \cdots & \cdots & \cdots &  0\\  
 0 & \cdots  &   &  \cdots &  0 & &  p & 0 &\cdots  & \cdots & 0 \\
 0 &  \cdots &   & \cdots  &  0 & &  0 & p & 0 & \cdots & 0 \\ 
 0 & \cdots  &   &  \cdots &  0& &  0 &  0&\ddots  & 0 &  0\\
 0 &  \cdots &   &  \cdots &  0& &  0 &  \cdots &  \cdots & p & 0 \\
\end{array}  \right),   \quad\quad({\cal F})  \]
$a_{i,j} \in W, \ \ a_{1,h} \in W^{\ast}$; i.e. it consists of blocks 
of sizes $(d$ or $c) \times (d$ or $c)$; in the left hand upper corner, 
which is of size $d\times d$, there are entries on the last column, 
named $a_{i,d}$, and the entries immediately below the diagonal are 
equal to 1; the left and lower block has only element unequal to zero, 
and it is 1; the right hand upper corner is unspecified, 
entries are called  $pa_{i,j}$;  the right hand lower corner, 
which is of size $c \times c$, has only entries immediately 
below the diagonal, and they are all equal to $p$.

Note that if a Dieudonn\'e module is defined by a matrix in displayed 
normal form then either the $p$-rank $f$ is maximal, $f=d$, 
and this happens if and only if $a_{1,d}$ is not divisible by $p$, or 
$f<d$, and  in that case  $a=1$. The $p$-rank is zero if and only if 
$a_{i,d} \equiv 0 \ \ (\bmod \ \ p), \ \ \forall 1 \leq i \leq d$.

\subsection{\bf Lemma.}\label{N} \  \Bl \  {\it Let $M$ be the 
Dieudonn\'e module of a p-divisible group $G$ over $k$ with $f(G)=0$. 
Suppose $a(G)=1$.  Then there exists a $W$-basis for 
$M$ on which $\cF$ has a matrix 
which is in normal form.} In this case the entries 
$a_{1,d}, \cdots , a_{d,d}$ are divisible by $p$, 
they can be chosen to be equal to zero.\\

\subsection{\bf Lemma}\label{CH} (of Cayley-Hamilton type). 
{\it Let $L$ be a field of characteristic $p$, 
let $W = W_{\infty}(L)$ be its 
ring of infinite Witt vectors. Let $X$ be a $p$-divisible group,  
with $\dm(G)=d$, and height$(G)=h$, with Dieudonn\'e module $\DD(X) = M$. 
Suppose there is a $W$-basis of $M$, such that the display-matrix  
$(a_{i,j})$  on this base gives an $\cF$-matrix in normal form as in} 
\ref{NF}.  {\it We write $e = X_1 = e_1$ for the first base vector.
Then for the expression}
$$P := \sum\limits_{i=1}^d \sum\limits_{j=d}^h \\  
p^{j-d} a_{i,j}^{\sigma^{h-j}} \cF^{h+i-j-1}  
\quad\mbox{\it we have}\quad
\cF^h{\cdot}e = P{\cdot}e.$$

\vspace{5mm}
\noindent
Note that we take powers of $\cF$ in the $\sigma$-linear sense, 
i.e. if the display matrix is $(a)$, i.e. $\cF$ is given by 
the matrix $(pa)$ as above, 
$$\cF^n \quad\mbox{is given by the matrix}
\quad (pa){\cdot}(pa^{\sigma}){\cdot} \cdots {\cdot}(pa^{\sigma^{n-1}}).$$
The exponent $h+i-j-1$ runs from $0= h+1-h-1$ to $h-1 =h+d-d-1$. 

\vn
Note that we do not claim that $P$ and $\cF^h$ have the same effect 
on all elements of $M$.

\vn
{\bf Proof.} Note that $\cF^{i-1}e_1 = e_i$ for $i \leq d$.\\
{\bf Claim.} {\it For $d \leq s < h$ we have}: 
$$\cF^sX = \left(\sum\limits_{i=1}^d \sum\limits_{j=d}^s \\ 
\cF^{s-j} p^{j-d} a_{i,j} \cF^{i-1}\right) X + p^{s-d}e_{s+1}.$$
This is correct for $s=d$. The induction step from $s$ to $s+1 <h$ 
follows from 
$\cF e_{s+1} = \left(\sum_{i=1}^d p \ a_{i,s+1}F^{i-1} \right) 
X + pe_{s+2}$. 
This proves the claim. Computing $\cF(\cF^{h-1}X)$ gives the 
desired formula.   \hfill \qed

\subsection{\bf Proposition.}\label{cN} 
{\it Let $k$ be an algebraically closed field of 
characteristic $p$, let $W = W_{\infty}(K)$ 
be its ring of infinite Witt vectors. 
Suppose $G$ is a $p$-divisible group over $k$ such that for its 
Dieudonn\'e module the map $\cF$ is given by a matrix in normal form. 
Let $P$ be the polynomial given in the previous proposition. 
The Newton polygon ${\cal N}(G)$ of this $p$-divisible group equals 
the Newton polygon given by the polynomial $P$.}\\
{\bf Proof.} Consider the $W[F]$-sub-module $M' \subset M$ generated by 
$X = e_1$. Note that $M'$ contains $X=e_1, e_2, \cdots , e_d$. 
Also  it contains $\cF e_d$, which equals $e_{d+1}$ plus a  
linear combination of the previous ones; hence  $e_{d+1} \in M'$. 
In the same way we see: $pe_{d+2} \in M'$, 
and $p^2e_{d+3} \in M'$ and so on. 
This shows that $M' \subset M = \oplus_{i \leq h} \ \  W{\cdot}e_i$ is 
of finite index. We see that  $M' = W[F]/W[F]{\cdot}(F^h -P)$. 
From this we see by the classification of $p$-divisible groups 
up to isogeny, 
that the result follows by \cite{Manin}, II.1. also see \cite{Demazure}, 
pp. 82-84. By \cite{Demazure}, page 82, Lemma 2 we conclude that the 
Newton polygon of $M'$ in case of the monic polynomial 
$\cF^h -  \sum_0^m \ \ b_i\cF^{m-i}$ is given by the lower convex hull of 
the pairs $\{(i,v(b_i)) \mid i\}$. 
Hence the proposition is proved. \hfill \qed

\subsection{\bf Corollary.}\label{diagram} {\it We take the 
notation as above. 
Suppose that every element 
$a_{i,j}, \ \ 1 \leq i \leq c, \ \ c\leq j \leq h$, 
is either equal to zero, or is a unit in $W(k)$.  Let $S$ be the set of 
pairs $(i,j)$ with $0 \leq i \leq c$ and $c \leq j \leq h$  for which 
the corresponding element is non-zero:
$$(i,j) \in S \quad\Longleftrightarrow\quad    a_{i,j} \not= 0.$$
Consider the image $T$ under  
$$S \rightarrow  T \subset \ZZ \times \ZZ \quad\mbox{given by}\quad (i,j) 
\mapsto (j+1-i,j-c).$$ 
Then ${\cal N}(X)$ is the lower convex hull of the set 
$T \subset \ZZ \times \ZZ$ and the point $(0,0)$; 
note that  $a_{1,h} \in W^{\ast}$, hence $(h,h-c=d) \in T$.
This can be visualized in the following diagram 
(we have pictured the case $d \leq h-d$):
\[ \begin{array}{llllllllllll}
 & & & &   &  &  && a_{c,h} & \cdots & &a_{1,h} \\
&&&&&&&  . &   \cdots&   & . & \\
&&&&&&                      a_{c,2c+2} &  \cdot  &  \cdot & \cdot&& \\
&&&&&          a_{c,2c+1}    &. &  .  &  .  &a_{1,2c+1}&& \\
&&&&  . &   \vdots &                     \vdots  &  \vdots &  . && &\\
&&     a_{c,c+1}&   \cdots&  a_{i,c+1}&      \cdots&                    
a_{2,c+1}       &   a_{1,c+1} &&&& \\
&a_{c,c} &   \cdots   &  a_{i,c} &  \cdots  & \cdots  &                   
a_{1,c}               &&&&& \\
\end{array} \]
Here the element $a_{c,c}$ is in the plane with coordinates
$(x=1,y=0)$ and $a_{1,h}$ has coordinates $(x=h,y=h-c=d)$. 
One erases the spots where 
$a_{i,j} = 0$, and one leaves the places where   $a_{i,j}$ is as unit. 
The lower convex hull of these points and $(0,0)$ (and $(h,h-c)$) equals 
${\cal N}(X)$}.\\

\vn
Theorem \ref{G} proves the following statement.\\
{\bf The weak Grothendieck conjecture.} {\it Given Newton polygons 
$\beta \prec \delta$ there exists a family of 
$p$-divisible groups over 
an integral base having 
$\delta$ as Newton polygon for the generic fiber, and $\beta$ as 
Newton polygon for for a closed fiber.} \\
However we will prove a much stronger result later.

\subsection{} For principally quasi-polarized $p$-divisible groups 
and for principally polarized abelian varieties we have an 
analogous method.

\subsection{} We fix an integer $g$. For every {\it symmetric} 
Newton polygon $\xi$ of height $2g$ we define:
$$\triangle(\xi) = \{(x,y) \in \ZZ \times \ZZ\mid  y < g, \ 
\ y < x \leq g, \ \ (x,y) \prec \xi\},$$
and
we write
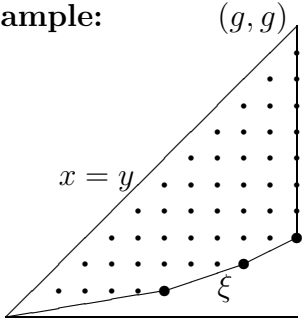
\begin{center}
\fbox{
$\sdm(\xi) := \#( \triangle(\xi)).$}
\end{center}

\begin{picture}(1000,120)(0,0)
\put(-18,110){\bf Example:}

\put(60,10){\circle*{4}}
\put(90,20){\circle*{4}}
\put(110,30){\circle*{4}}

\put(0,0){\line(1,0){110}}

\put(0,0){\line(1,1){110}}
\put(20,50){$x=y$}
\put(80,110){$(g,g)$}

\put(0,0){\line(6,1){60}}
\put(60,10){\line(3,1){30}}
\put(90,20){\line(2,1){20}}
\put(110,30){\line(0,1){80}}

\put(20,10){\circle*{2}}
\put(30,10){\circle*{2}}
\put(40,10){\circle*{2}}
\put(50,10){\circle*{2}}

\put(30,20){\circle*{2}}
\put(40,20){\circle*{2}}
\put(50,20){\circle*{2}}
\put(60,20){\circle*{2}}
\put(70,20){\circle*{2}}
\put(80,20){\circle*{2}}

\put(40,30){\circle*{2}}
\put(50,30){\circle*{2}}
\put(60,30){\circle*{2}}
\put(70,30){\circle*{2}}
\put(80,30){\circle*{2}}
\put(90,30){\circle*{2}}
\put(100,30){\circle*{2}}
\put(110,30){\circle*{2}}

\put(50,40){\circle*{2}}
\put(60,40){\circle*{2}}
\put(70,40){\circle*{2}}
\put(80,40){\circle*{2}}
\put(90,40){\circle*{2}}
\put(100,40){\circle*{2}}
\put(110,40){\circle*{2}}

\put(60,50){\circle*{2}}
\put(70,50){\circle*{2}}
\put(80,50){\circle*{2}}
\put(90,50){\circle*{2}}
\put(100,50){\circle*{2}}
\put(110,50){\circle*{2}}

\put(70,60){\circle*{2}}
\put(80,60){\circle*{2}}
\put(90,60){\circle*{2}}
\put(100,60){\circle*{2}}
\put(110,60){\circle*{2}}

\put(80,70){\circle*{2}}
\put(90,70){\circle*{2}}
\put(100,70){\circle*{2}}
\put(110,70){\circle*{2}}

\put(90,80){\circle*{2}}
\put(100,80){\circle*{2}}
\put(110,80){\circle*{2}}

\put(100,90){\circle*{2}}
\put(110,90){\circle*{2}}

\put(110,100){\circle*{2}}

\put(80,8){$\xi$}
\put(130,80){$\dm(\cW_{\xi}(\cA_{g,1} \otimes \FF_p)) = 
\#( \triangle(\xi))$}

\put(150,60){$\xi = (5,1)+(2,1)+2{\cdot}(1,1)+(1,2)+(1,5),  $}
\put(170,40){ g=11; slopes: 
$\{ 6 \times \frac{5}{6}, 3 \times \frac{2}{3}, 
4 \times \frac{1}{2},  3 \times \frac{1}{3}, 6 \times \frac{1}{6}\}$.}
\put(150,20){This case: 
$\dm(\cW_{\xi}(\cA_{g,1} \otimes \FF_p))= \sdm(\xi) 
= 48,$ see \ref{dimp}}
\end{picture}


\vn 
Suppose given a $p$-divisible group $X_0$ over $k$ of dimension $g$ with 
a principal quasi-polarization. We write $\cN(X_0)=\gamma$; this is 
a symmetric Newton polygon.  We write $\cD = {\rm Def}(X_0,\lambda)$ 
for the universal deformation space. For every symmetric Newton polygon 
$\xi$ with $\xi \succ \gamma$ we define $W_{\xi} \subset \cD$ as 
the maximal closed, reduced formal subscheme 
carrying all fibers with Newton polygon 
equal or above $\xi$; this space exists by Grothendieck-Katz, 
see \cite{Katz-Slope}, Th. 2.3.1 on page 143. Note that $W_{\rho} = \cD$, 
where $\rho = g{\cdot}((1,0)+(0,1))$.

\subsection{}\label{pG}{\bf Theorem} (NP-strata for principally 
quasi-polarized formal groups).
{\it Suppose $a(X_0) \leq 1$. For every symmetric 
$\xi \succ \gamma := \cN(G_0)$ 
we have: $\dm(W_{\xi}) = \sdm({\xi})$. 
The strata $W_{\xi}$ are nested as given by the partial ordering 
on symmetric Newton polygons, i.e. 
$$W_{\xi} \quad\subset\quad W_{\delta} \quad\Longleftrightarrow\quad 
\triangle(\xi) \subset  \triangle(\delta) \quad\Longleftrightarrow\quad 
\xi \prec \delta.$$
Generically on $W_{\xi}$ the fibers have Newton polygon equal to $\xi$.
We can choose a coordinate system on $\cD$ in which all $W_{\xi}$ are 
given by linear equations.}

\subsection{}\label{CpG}{\bf Corollary}. {\it Suppose given a 
principally polarized abelian variety $(A_0,\lambda_0)$ 
over $k$. Strata in 
${\rm Def}(X_0,\lambda_0)$  according to Newton polygons are exactly 
as in \ref{pG}. In particular, the  fiber above the  generic point of  
$W_{\xi}$ is a principally polarized abelian scheme over $\Spec(B_{\xi})$ 
having Newton polygon equal to} $\xi$ (for $B_{\xi}$, see the proof 
of  \ref{pG} below; for the notion ``generic point of  $W_{\xi}$ " 
see the proof of \ref{CpG} below; also see \ref{algebr}).\\
{\bf Proof.} We write $(A_0,\lambda_0)[p^{\infty}] =: (X_0,\lambda_0)$.  
By  Serre-Tate theory, see \cite{Katz-ST}, Section 1, 
the formal deformation 
spaces of $(A_0,\lambda_0)$  and of $(X_0,\lambda_0)$ are canonically 
isomorphic, say $(A,\lambda) \to \Spf(R)$ and $(\cX,\lambda) \to \Spf(R)$ 
and $(A,\lambda)[p^{\infty}] \cong (\cX,\lambda)$.
By Chow-Grothendieck, see \cite{EGA}, III$^1$.5.4 (this is also called 
a theorem of ``GAGA-type"),  the  formal polarized abelian scheme  
is algebraizable, and we obtain the universal deformation as a polarized 
abelian scheme over $\Spec(R)$, and we can consider the generic point of  
$W_{\xi} \subset \Spec(R)$. Hence the Newton polygon of fibers can be 
read off from the fibers in  $(\cX,\lambda) \to \Spec(R)$. 
This proves that \ref{CpG} follows from \ref{pG}. \hfill$\Box$

\vn
{\bf Proof of}\ \  \ref{pG}. The proof of this theorem is 
analogous to the proof of \ref{G}. We use  the diagram 
\[ \begin{array}{lllllll}
&&&&&&-1 \\
&&&    T_{g,g} &  \cdots&  T_{1,g}& \\
&&   . &       \vdots &      \cdot      &          & \\
1&   T_{g,1}&  \cdots&   T_{1,1} &&&\\

\end{array} \]
Here $X_{i,j}$, \ \ $1 \leq i, j \leq g$,  is written on the place with 
coordinates $(g-i+j,j-1)$.
We use the ring 
$$B := \frac{k[[T_{i,j}; 1 \leq i,j \leq g]]}{(T_{k\ell}-T_{\ell k})}, 
\quad  T_{i,j} = Z_{(g-i+j,j-1)}, \ \ (g-i+j,j-1) \in \triangle.$$
Note that $B=k[[T_{i,j}\mid 1 \leq i \leq j \leq g]] = 
k[[Z_{x,y}\mid (x,y) 
\in \triangle]]$. 
For a symmetric $\xi$ with $\xi \succ \cN(X_0)$ we consider
$$B_{\xi} = \frac{k[[T_{i,j}; 1 \leq i,j \leq g]]}{(T_{k\ell}-T_{\ell k}, 
\quad\mbox{and}\quad Z_{(x,y)} \ \ \forall (x,y) \not\in \triangle(\xi))} 
\quad\cong\quad k[[Z_{(x,y)} \mid (x,y) \in \triangle(\xi))]] .$$
With these notations, applying \ref{CH} and \ref{diagram} we finish the 
proof of \ref{pG} as we did in the proof of \ref{G} above.  \hfill$\Box$

\subsection{\bf A conjecture by Manin.}\label{MC} 
Let $A$ be an abelian variety. 
The Newton polygon $\cN(A)$ is symmetric \ref{sym}. A conjecture by Manin 
expects the converse to hold:\\
{\bf Conjecture,} see \cite{Manin}, page 76, Conjecture 2. 
{\it For any symmetric Newton polygon $\xi$ there exists  
an abelian variety $A$ such that $\cN(A) = \xi$.}

\vn
This was proved in the Honda-Tate theory, see \ref{H.T}. 
{\it We sketch a pure characteristic $p$ proof,} 
see \cite{FO-CH}, Section 5. 
It is not difficult to show 
that there exists a principally polarized supersingular abelian variety 
$(A_0,\lambda_0)$, see \cite{FO-CH}, Section 4; this also follows from   
\cite{Li.FO}, 4.9.  By \ref{pG} it follows that 
$\cW^0_{\xi}({\rm Def}(A_0,\lambda_0))$ is non-empty, 
{\it which proves the Manin conjecture.} \hfill \qed


\section{Hilbert modular varieties}\label{6}
\medbreak

We discuss Hilbert modular varieties over \(\Fpbar\) in this
section. (Recall that \(\Fpbar\) is the algebraic closure
of \(\Fp\).) 
A Hilbert modular variety attached to a totally real number
field \(F\) classifies  ``abelian varieties with real multiplication
by \(\ringO_F\)''. An abelian variety \(A\) is said to have 
``real multiplication by \(\ringO_F\) if \(\dim(A)=[F:\bbQ]\)
and there is an embedding \(\ringO_F\hookrightarrow \End(A)\);
the terminology ``fake elliptic curve'' was used by some
authors.  The moduli space of such objects behave very much
like the modular curve, except that its dimension is equal
to \([F:\bbQ]\).  Similar to the modular curve, a Hilbert
modular variety attached to a totally real number field \(F\) 
has a family of Hecke correspondences 
coming from the group \(\SL_2(F\otimes_{\bbQ}\bbA_f^{(p)})\) or 
\(\GL_2(\bbA_f^{(p)})\) depending on the definition one uses.
Hilbert modular varieties are closely related to 
modular forms for \(\GL_2\) over totally real fields
and the arithmetic of totally real fields.
\smallbreak

Besides their intrinsic interest, Hilbert modular varieties
plays an essential role in the Hecke orbit problem for
Siegel modular varieties.  This connection results from a
special property of \(\calA_{g,1,n}\) which is not shared
by all modular varieties of PEL type: For every 
\(\Fpbar\)-point \(x_0\) of \(\calA_{g,1,n}\), 
there exists a Hilbert modular variety \(\calM\)
and an isogeny correspondence \(R\) on \(\calA_{g,1,n}\)
such that \(x_0\) is contained in the image of \(\calM\)
under the isogeny correspondence \(R\).
See \ref{thm_hilb} for a precise formulation, and also the
beginning of \S\ref{8}.
\bigbreak

\noindent{\scshape References.}\enspace 
\cite{Rap}, \cite{D.P}, \cite{Geer-HMS} Chap X,
\cite{Deligne.Ribet}, \cite{Goren.FO}, \cite{Yu}.
\bigbreak


Let \(F_1,\ldots,F_r\) are totally real number fields,
and let \(E:=F_1\times\cdots\times F_r\).
Let \(\ringO_E=\ringO_{F_1}\times\cdots\times\ringO_{F_r}\) be
the product of the rings of integers of \(F_1,\ldots,F_r\).
Let \(\calL_i\) be an invertible \(\ringO_{F_i}\)-module,
and let \(\calL\) be the invertible \(\ringO_E\)-module
\(\calL=\calL_1\times\cdots\times\calL_r\).
\medbreak

\begin{defn}
Notation as above.
A \emph{notion of positivity}
on an invertible \(\ringO_E\)-module \(\,\calL\,\)
is a union \(\calL^{+}\)
of connected components of 
\(\calL\otimes_{\bbQ}\bbR\) such that 
\(\calL\otimes_{\bbQ}\bbR\) is the disjoint union of 
\(\calL^{+}\) and \(-\calL^{+}\).
\end{defn}

\medbreak

\begin{defn}
\begin{itemize}
\rmitem[(i)] An \(\ringO_E\)-linear abelian
scheme is a pair \((A\to S, \iota)\), where 
\(A\to S\) is an abelian scheme, and
\(\iota:\ringO_E\to \End_S(A)\) is an injective ring
homomorphism such that \(\iota(1)={\rm Id}_A\).
Notice that every \(\ringO_E\)-linear abelian
scheme \((A\to S,\iota)\) as above decomposes as a product
\((A_1\to S,\iota_1)\times\cdots\times(A_r\to S,\iota_r)\).
Here \((A_i,\iota_i)\) is an \(\ringO_{F_i}\)-linear abelian
scheme for \(i=1,\ldots,r\), and 
\(A=A_1\times_S\cdots\times_S A_r\).

\rmitem [(ii)] An \(\ringO_E\)-linear abelian scheme \((A\to S, \iota)\)
is said to be of HB-type if 
\(\dim(A/S)=\dim_{\bbQ}(E)\).

\rmitem[(iii)] An \(\ringO_E\)-linear polarization of an
\(\ringO_E\)-linear abelian scheme is a polarization
\(\lambda: A\to A^t\) 
such that \(\lambda\circ \iota(u)=\iota(u)^t\circ \lambda\)
for all \(u\in\ringO_E\).
\end{itemize}
\end{defn}
\medbreak

\subsection{\bf Exercise.}
Suppose that \((A\to S, \iota)\) is an
\(\ringO_E\)-linear abelian scheme, and
\((A\to S,\iota)=(A_1\to S,\iota_1)\times\cdots\times(A_r\to S,\iota_r)\)
as in (i).
Show that \((A_1\to S,\iota_1)\) is an \(\ringO_{F_i}\)-linear
abelian scheme of HB-type for \(i=1,\ldots,r\).

\medbreak

\subsection{\bf Exercise.}
Show that every \(\ringO_E\)-linear abelian
variety of HB-type over a field 
admits an \(\ringO_E\)-linear polarization.

\bigbreak

\begin{defn}
Let \(E_p=\prod_{j=1}^s F_{v_j}\) be a product of finite
extension fields \(F_{v_j}\) of \(\Qp\). Let
\(\ringO_{E_p}=\prod_{j=1}^s \ringO_{F_{v_j}}\) be the product
of the rings of elements in \(F_{v_j}\) which are integral over \(\Zp\).
\begin{itemize}

\rmitem[(i)] An \(\ringO_{E_p}\)-linear BT-group
is a pair \((X\to S,\iota)\), where 
\(X\to S\) is a BT-group, and 
\(\iota:\ringO_E\otimes_{\bbZ}\Zp \to {\rm End}_{S}(X)\)
is an injective ring homomorphism such that \(\iota(1)={\rm Id}_X\).
Every \((\ringO_{E_p})\)-linear BT-group
\((X\to S,\iota)\) decomposes canonically into a product
\((X\to S,\iota)=\prod_{j=1}^s\, (X_j,\iota_j)\),
where \((X_j,\iota_j)\) is an \(\ringO_{F_{v_j}}\)-linear
BT-group, defined to be the image of the idempotent in 
\(\ringO_{E_p}\) corresponding to the factor \(\ringO_{F_{v_j}}\) of
\(\ringO_{E_p}\).

\rmitem[(ii)] An \(\ringO_{E_p}\)-linear BT-group
\((X\to S,\iota)\) is said to have rank two if in the decomposition
\((X\to S,\iota)=\prod_{j=1}^s\, (X_j,\iota_j)\)
in (i) above we have
\({\rm ht}(X_j/S)=2\,[F_{v_j}:\Qp]\) for all \(j=1,\ldots,s\).

\rmitem[(iii)] An \(\ringO_{E_p}\)-linear polarization 
\((\ringO_E\otimes_{\bbZ}\Zp)\)-linear BT-group
\((X\to S,\iota)\) is a symmetric isogeny
\(\lambda: X\to X^t\) 
such that \(\lambda\circ \iota(u)=\iota(u)^t\circ \lambda\)
for all \(u\in\ringO_{E_p}\).

\rmitem[(iv)] A rank-two \(\ringO_{E_p}\)-linear BT-group
\((X\to S,\iota)\) is of HB-type if it admits an
\(\ringO_{E_p}\)-linear polarization.
\end{itemize}

\end{defn}
\medbreak

\subsection{\bf Exercise.} 
Show that for every \(\ringO_E\)-linear
abelian scheme of HB-type \((A\to S,\iota)\), the associated
\((\ringO_E\otimes_{\bbZ}\Zp)\)-linear BT-group
\((A[p^{\infty}],\iota[p^{\infty}])\) is of HB-type.

\begin{defn}\label{defn-hilb}
Let \(E=F_1\times\cdots\times F_r\), where \(F_1,\ldots,F_r\)
are totally real number fields.
Let \(\ringO_E=\ringO_{F_1}\times\cdots\times\ringO_{F_r}\) be
the product of the ring of integers of \(F_1,\ldots,F_r\).
Let \(k\supset \Fp\) be an algebraically closed field as before.
Let \(n\geq 3\) be an integer such that \((n,p)=1\).
Let \((\calL,\calL^+)\) be an invertible \(\ringO_E\)-module
with a notion of positivity.
The Hilbert modular variety \(\calM_{E,\calL,\calL^+,n}\) over \(k\) is a 
smooth scheme over \(k\) of dimension \([E:\bbQ]\) such that
for every \(k\)-scheme \(S\) the set of \(S\)-valued points
of \(\calM_{E,\calL,\calL^+_,n}\) is the set of isomorphism class of 
5-tuples
\((A\to S, \iota, \calL, \calL^{+}, \lambda, \eta)\),
where 
\begin{itemize}
\rmitem[(i)] \((A\to S, \iota)\) is an \(\ringO_E\)-linear abelian scheme
of HB-type;
%
%
%

\rmitem[(ii)] \(\lambda:\calL\to \Hom_{\ringO_E}^{\rm sym}(A, A^t)\)
is an \(\ringO_E\)-linear homomorphism such that
\(\lambda(u)\) is an \(\ringO_E\)-linear polarization of 
\(A\) for every \(u\in\calL\cap \calL^+\),
and the homomorphism 
\( A\otimes_{\ringO_E}\calL
\xrightarrow{\sim}A^t\) 
induced by \(\lambda\) is an isomorphism of abelian schemes.

\item[(iii)] \(\eta\) is an \(\ringO_E\)-linear level-\(n\) 
structure for \(A\to S\), i.e.\
an \(\ringO_E\)-linear isomorphism from the constant group
scheme \((\ringO_E/n\ringO_E)^2_S\) to \(A[n]\).
\end{itemize}
\end{defn}
\medbreak

\begin{rem} Let \((A\to S, \iota, \lambda,\eta)\) be
an \(\ringO_E\)-linear abelian scheme with polarization sheaf
by \((\calL,\calL^+)\) and a level-\(n\) 
satisfying the condition in (ii) above, 
Then the \(\ringO_E\)-linear polarization \(\lambda\) induces
an \(\ringO_E/n\ringO_E\)-linear isomorphism 
\[
(\ringO_E/n\ringO_E)=\bigwedge^2 (\ringO_E/n\ringO_E)^2\xrightarrow{\sim}
\calL^{-1}\calD_E^{-1}\otimes_{\bbZ} \mu_n
\]
over \(S\), where \(\calD_E\) denotes the invertible \(\ringO_E\)-module
\(\calD_{F_1}\times\cdots\times\calD_{F_r}\).
This isomorphism is a discrete invariant of 
the quadruple \((A\to S, \iota, \lambda,\eta)\).
The above invariant
defines a morphism \(f_n\) from the
Hilbert modular variety \(\calM_{E,\calL,\calL^+,n}\) to the 
finite \'etale scheme \(\Xi_{E,\calL,n}\) over \(k\), 
where the finite \'etale \(k\)-scheme \(\Xi_{E,\calL,n}\) is characterized by
\(\Xi_{E,\calL,n}(k):={\rm Isom}(\ringO_E/n\ringO_E, 
\calL^{-1}\calD_E^{-1}\otimes_{\bbZ} \mu_n)\) 
Notice that \(\Xi_{E,\calL,n}\) is an 
\((\ringO_E/n\ringO_E)^{\times}\)-torsor; it is constant over
\(k\) because \(k\) is algebraically closed.
The morphism \(f_n\) is faithfully flat.
\smallbreak

Although we defined the Hilbert modular variety \(\calM_{E,\calL,\calL^+,n}\)
over an algebraically closed field \(k\supset\Fp\),
we could have defined it over \(\Fp\).
Then we should use the \'etale \((\ringO_E/n\ringO_E)^{\times}\)-torsor
\(\Xi_{E,\calL,n}:=\underline{\rm Isom}(\ringO_E/n\ringO_E, 
\calL^{-1}\calD_E^{-1}\otimes_{\bbZ} \mu_n)\) over \(\Fp\),
and we have a faithfully flat
morphism \(f_n:\calM_{E,\calL,\calL^+,n}\to \Xi_{E,\calL,n}\) over \(\Fp\).
\end{rem}
\medbreak

\begin{rem}
(i) We have followed \cite{D.P} in the definition
of Hilbert Modular varieties, except that \(E\) is a product
of totally real number fields, rather than a totally real number
field as in \cite{D.P}.
\smallbreak

\noindent
(ii) The product decompositions
\(\ringO_E=\ringO_{F_1}\times\cdots\times \ringO_{F_r}\) and
\((\calL,\calL^+)=(\calL_1,\calL_1^+)\times\cdots\times (\calL_r,\calL_r^+)\)
induces a natural isomorphism
\[
\calM_{E,\calL,\calL^+,n}\xrightarrow{\sim} 
\calM_{F_1,\calL_1,\calL_1^+,n}\times\cdots\times
\calM_{F_r,\calL_r,\calL_r^+n}\,.
\]
\end{rem}
\medbreak

\begin{rem}
The \(\ringO_E\)-linear homomorphism \(\lambda\) in
Def.\ \ref{defn-hilb} should be thought of as specifying 
a family of \(\ringO_E\)-linear polarizations, not just one 
polarization: every element \(u\in \calL\cap \calL^+\)
gives a polarization \(\lambda(u)\) on \(A\to S\).
Notice that given a point \(x_0=[(A,\iota,\lambda,\eta)]\) 
in \(\calM_{E,\calL,\calL^+,n}(k)\), 
there may not exist an \(\ringO_E\)-linear principal
polarization on \(A\), because  that means that
the element of the strict ideal class group
represented by \((\calL,\calL^+)\) is trivial.
However every point \([(A,\iota,\lambda,\eta)]\)
of \(\calM_{E,\calL,\calL^+,n}\) admits an \(\ringO_E\)-linear
polarization of degree prime to \(p\), because
there exists an element \(u\in \calL^{+}\) such that
\({\rm Card}(\calL/\ringO_E\cdot u)\) is not divisible by \(p\).
In \cite{chiafuyu:4fold} and \cite{Yu} a version of Hilbert modular
varieties was defined by specifying a polarization degree \(d\)
which is prime to \(p\).
The resulting Hilbert modular variety is not necessarily irreducible
over \(\Fpbar\);
rather it is a disjoint union of modular varieties
of the form \(\calM_{E,\calL,\calL^+,n}\).
\end{rem}

\begin{thm}\Bl \label{thm:irredHMV}
Notation as above.
\begin{itemize}

\rmitem[(i)] The modular variety \(\calM_{E,\calL,\calL^+,n}\)
over the algebraically closed field \(k\supset \Fp\) is
normal and is a local complete intersection. 
Its dimension is equal to \(\dim_{\bbQ}(E)\). 

\rmitem[(ii)] Every fiber of \(f_n:\calM_{E,\calL,\calL^+,n}\to
\Xi_{E,\calL,n}\) is irreducible.

\rmitem[(iii)] The morphism \(f_n\) is smooth outside a closed
subscheme of \(\calM_{E,\calL,\calL^+,n}\) of codimension at least two.
\end{itemize}

\end{thm}

\noindent{\bf Remark.} 

\noindent (i) See \cite{D.P} and 
for a proof of Thm.\ \ref{thm:irredHMV}
which uses the arithmetic toroidal compactification 
constructed in \cite{Rap}.

\noindent (ii) The modular variety \(\calM_{E,\calL,\calL^+,n}\)
is not smooth over \(k\) if any one of the totally real fields
\(F_i\) is ramified above \(p\).
\bigbreak

\subsection{\scshape Hecke orbits on Hilbert modular varieties}
\smallbreak

Let \(E\), \(\calL\) and \(\calL^+\) be as before.
Denote by \(\,\widetilde{\calM_{E,\calL,\calL^+}}\)
the projective system \(\,\left(\calM_{E,\calL,\calL^+,n}\right)_n\) of
Hilbert modular varieties over \(\Fpbar\), where
\(n\) runs through all positive integer such that 
\(n\geq 3\) and \({\rm gcd}(n,p)=1\).
It is clear that the profinite group 
\(\SL_2(\ringO_E\otimes_{\bbZ}\widehat{\bbZ}^{(p)})\) operates
on the tower \(\,\widetilde{\calM_{E,\calL,\calL^+}}\), by
pre-composing with the \(\ringO_E\)-linear level structures.
Here \(\widehat{\bbZ}^{(p)}=\prod_{\ell\neq p}\Zl\).
The transition maps in the projective system are
\[
\pi_{mn,n}:\calM_{E,\calL,\calL^+,mn}\to
\calM_{E,\calL,\calL^+,n}
\qquad
(mn,p)=1, n\geq 3, m\geq 1\,.
\]
The map \(\pi_{mn,n}\) is defined by the following construction.
Let \([m]:(\ringO_E/n\ringO_E)^2\to (\ringO_E/mn\ringO_E)^2\) be the
injection induced by ``multiplication by \(m\)''.
Given a point \((A,\iota,\lambda,\eta)\) of 
\(\calM_{E,\calL,\calL^+,mn}\),
the composition \(\eta\circ [m]\) factors through the inclusion
\(i_{m,n}:A[m]\hookrightarrow A[mn]\) to give a level-\(n\) structure
\(\eta'\) such that \(\eta\circ [m]=i_{m,n}\circ \eta'\).
\smallbreak

Let \(\,\widetilde{\Xi_E}\,\) be the projective system
\(\,\left(\Xi_{E,n}\right)_{n}\), 
where \(n\) also runs through all positive integer such that 
\(n\geq 3\) and \((n,p)=1\).
The transition maps are defined similarly.
The maps \(f_n:\calM_{E,\calL,\calL^+,n}\to \Xi_{E,n}\)
define a map \(\tilde{f}:\widetilde{\calM_{E,\calL,\calL^+}}
\to \widetilde{\Xi_E}\,\) between projective systems.
\smallbreak

It is clear that the profinite group 
\(\SL_2(\ringO_E\otimes_{\bbZ}\widehat{\bbZ}^{(p)})\) operates
on the right of the tower \(\,\widetilde{\calM_{E,\calL,\calL^+}}\), by
pre-composing with the \(\ringO_E\)-linear level structures.
Moreover this action is compatible with the map
\(\tilde{f}:\widetilde{\calM_{E,\calL,\calL^+}}
\to \widetilde{\Xi_E}\,\) between projective systems.
\medbreak

The above right action of the compact group 
\(\SL_2(\ringO_E\otimes_{\bbZ}\widehat{\bbZ}^{(p)})\)
on the projective system \(\,\widetilde{\calM_{E,\calL,\calL^+}}\)
extends to a right action of 
\(\SL_2(E\otimes_{\bbQ}\bbA_f^{(p)})\) on 
\(\,\widetilde{\calM_{E,\calL,\calL^+}}\).
Again this action is compatible with 
\(\tilde{f}:\widetilde{\calM_{E,\calL,\calL^+}}
\to \widetilde{\Xi_E}\,\).
This action can be described as follows.
A geometric point of \(\,\widetilde{\calM_{E,\calL,\calL^+}}\,\)
is a quadruple \((A,\iota_A,\lambda_A,\tilde{\eta}_A)\),
where the infinite prime-to-\(p\) level structure 
\[
\tilde{\eta}_A:\coprod_{\ell\neq p}(\ringO_E[1/\ell]/\ringO_E)
\xrightarrow{\sim}\coprod_{\ell\neq p}A[\ell^{\infty}]
\]
is induced by a compatible system of level-\(n\)-structures,
\(n\) running through integers such that \((n,p)=1\) and \(n\geq 3\).
Suppose that we have an element
\(\gamma\in \SL_2(E\otimes_{\bbQ}\bbA_f^{(p)})\),
and \(m\,\gamma\)
belongs to \({\rm M}_2(\ringO_E\otimes_{\bbZ}\widehat{\bbZ}^{(p)})\),
where \(m\) is a non-zero integer which is prime to \(p\).
Then the image of the point \((A,\iota_A,\lambda_A,\tilde{\eta}_A)\)
under \(\gamma\) is a quadruple
\((B,\iota_B,\lambda_B,\tilde{\eta}_B)\) 
such that there exists an \(\ringO_E\)-linear 
prime-to-\(p\) isogeny \(m\beta:B\to A\)
such that the diagram
\[
\xymatrix{
\coprod_{\ell\neq p}(\ringO_E[1/\ell]/\ringO_E)^2\ar[r]^{\qquad\tilde{\eta}_A}
&\coprod_{\ell\neq p} A[\ell^{\infty}]\\
\coprod_{\ell\neq p}(\ringO_E[1/\ell]/\ringO_E)^2\ar[r]_{\qquad\tilde{\eta}_B}
\ar[u]^{m\gamma}
&\coprod_{\ell\neq p} B[\ell^{\infty}]\ar[u]_{m\beta}
}
\]
commutes. In the above notation, as the point 
\((A,\iota_A, \lambda_A,\tilde\eta_A)\) varies, we get a
prime-to-\(p\) quasi-isogeny \(\beta=m^{-1}\cdot (m\beta)\)
attached to \(\gamma\), between the universal abelian schemes.
\smallbreak

On a fixed level \(\calM_{E,\calL,\calL^+,n}\), the action of 
\(\SL_2(E\otimes_{\bbQ}\bbA_f^{(p)})\) on the projective system
\(\,\widetilde{\calM_{E,\calL,\calL^+}}\,\) induce
a family of finite \'etale correspondences, which
will be called \(\SL_2(E\otimes_{\bbQ}\bbA_f^{(p)})\)-Hecke 
correspondences on \(\calM_{E,\calL,\calL^+,n}\),
or prime-to-\(p\) \(\SL_2\)-Hecke correspondences for short.
Suppose \(x_0\) is a geometric point of \(\calM_{E,\calL,\calL^+,n}\),
and \(\tilde{x}\) is point of \(\,\widetilde{\calM_{E,\calL,\calL^+}}\,\)
lifting \(x_0\).  Then the prime-to-\(p\) \(\SL_2\)-Hecke orbit of
\(x_0\), denoted \(\calH^{(p)}_{\SL_2}(x_0)\),
is the image in \(\calM_{E,\calL,\calL^+,n}\)
of the orbit \(\SL_2(E\otimes_{\bbQ}\bbA_f^{(p)})\cdot x_0\).
The set \(\calH^{(p)}_{\SL_2}(x_0)\) is countable.

\begin{thm}\label{STHMV}
Let \(x_0=[(A_0,\iota_0,\lambda_0,\eta_0)]\in\calM_{E,\calL,\calL^+,n}(k)\)
be a closed point of \(\calM_{E,\calL,\calL^+,n}\) such that 
\(A_0\) is an ordinary abelian scheme.
Let \(\,\Sigma_{E,p}=\{\wp_1,\ldots,\wp_s\}\,\) be the set of 
all prime ideals of \(\ringO_E\) containing \(p\).
Then we have a natural isomorphism
\[
\calM_{E,\calL,\calL^+,n}^{/x_0}
\cong
\prod_{j=1}^s\,\underline{\rm Hom}_{\Zp}
\left({\rm T}_p(A_0[\wp_j^{\infty}]_{\rm et})\otimes_{(\ringO_E\otimes\Zp)}
\!{\rm T}_p(A_0^t[\wp_j^{\infty}]_{\rm et}),\,
\widehat{\bbG_m}
\right)\,.
\]
In particular, the formal completion of the Hilbert modular variety
\(\calM_{E,\calL,\calL^+,n}\) at the ordinary point \(x_0\)
has a natural structure as a \([E:\bbQ]\)-dimensional 
\((\ringO_E\otimes_{\bbZ}\Zp)\)-linear formal torus,
non-canonically isomorphic to 
\((\ringO_E\otimes_{\bbZ}{\Zp})\otimes_{\Zp}\widehat{\bbG_m}\).
\end{thm}

\proof By the Serre-Tate theorem, we have
\[
\calM_{E,\calL,\calL^+,n}^{/x_0}
\cong
\prod_{j=1}^s\,\underline{\rm Hom}_{\ringO_E\otimes\Zp}
\left(
{\rm T}_p(A_0[\wp_j^{\infty}]_{\rm et}),\,
\widehat{A_0[\wp_j^{\infty}]_{\rm mult}}
\right)\,,
\]
where \(\widehat{A_0[\wp_j^{\infty}]_{\rm mult}}\) is the formal
torus attached to \(A_0[\wp_j^{\infty}]_{\rm mult}\), or equivalently
the formal completion of \(A_0\).
The character group of the last formal torus is naturally isomorphic
to the \(p\)-adic Tate module 
\({\rm T}_p(A_0^t[\wp_j^{\infty}]_{\rm et})\)
attached to the maximal \'etale quotient of 
\(A_0^t[\wp_j^{\infty}]_{\rm et})\).
\hfill\qed
\medbreak

\begin{prop}\label{HMVlocalaction}  Notation as in \ref{STHMV}.
Assume that \(k=\Fpbar\), so 
\(x_0=[(A_0,\iota_0,\lambda_0,\eta_0)]
\in\calM_{E,\calL,\calL^+,n}(\Fpbar)\)
and \(A_0\) is an ordinary \(\ringO_E\)-linear abelian variety
of HB-type over \(\Fpbar\).
\begin{itemize}
\rmitem[(i)] There exists totally imaginary quadratic extensions
\(K_i\) of \(F_i\), \(i=1,\ldots,r\) such that
\(\End_{\ringO_E}^0(A_0)\cong K_1\times\cdots \times K_r=:K\).
Moreover, for every prime ideal \(\wp_j\) of \(\ringO_E\) containing \(p\),
we have
\begin{multline*}
\End_{\ringO_E}(A_0)\otimes_{\ringO_E}\ringO_{E_{\wp_j}}
\xrightarrow{\sim}
\End_{\ringO_{E_{\wp_j}}}(A_0[\wp_j^{\infty}]_{\rm mult})
\times 
\End_{\ringO_{E_{\wp_j}}}(A_0[\wp_j^{\infty}]_{\rm et})\\
\cong \ringO_{E_{\wp_j}}\times \ringO_{E_{\wp_j}}
\xleftarrow{\sim} \ringO_K\otimes_{\ringO_E}\ringO_{E_{\wp_j}}\,.
\end{multline*}
In particular, the quadratic extension \(K_i/F_i\) is split above
every place of \(F_i\) above \(p\), for all \(i=1,\ldots,r\).

\rmitem[(ii)] Let \(H_{x_0}=
\left\{u\in(\ringO_E\otimes\Zp)^{\times}
\,\mid\, u\cdot \bar u=1\right\}\,\).
Then both projections
\[{\rm pr}_1:H_{x_0}\to \prod_{\wp\in\Sigma_{E,p}}
\left(\End_{\ringO_{E_{\wp}}}(A_0[\wp_j^{\infty}]_{\rm mult})
\right)^{\times}\cong\prod_{\wp\in\Sigma_{E,p}}\ringO_{E_{\wp}}^{\times}\]
and 
\[{\rm pr}_2:H_{x_0}\to \prod_{\wp \in\Sigma_{E,p}}
\left(\End_{\ringO_{E_{\wp}}}(A_0[\wp^{\infty}]_{\rm et})
\right)^{\times}\cong\prod_{\wp\in\Sigma_{E,p}}\ringO_{E_{\wp}}^{\times}\]
are isomorphisms.
Here \(\Sigma_{E,p}\) denotes the set consisting 
of all prime ideals of \(\ringO_E\) which contain \(p\).

\rmitem[(iii)] The group \(H\) operates on the
\((\ringO_E\otimes_{\bbZ}\Zp)\)-linear formal torus
\(\calM_{E,\calL,\calL^+,n}^{/x_0}\)
through the character 
\[
H\ \ni t\longmapsto {\rm pr}_1(t)^2\in (\ringO_E\otimes_{\bbZ}\Zp)^{\times}
\,.
\]

\rmitem[(iv)] Notation as in {\rm (ii)} above.
Let \(Z\) be a reduced, irreducible closed formal
subscheme of the formal scheme \(\calM_{E,\calL,\calL^+,n}^{/x_0}\)
which is stable under the natural action of an open subgroup \(U_{x_0}\)
of \(H_{x_0}\) on \(\calM_{E,\calL,\calL^+,n}^{/x_0}\).
Then there exists a subset \(S\subset \Sigma_{E,p}\)
such that
\[
Z=\prod_{\wp\in S}
\underline{\rm Hom}_{\Zp}
\left({\rm T}_p(A_0[\wp^{\infty}]_{\rm et})\otimes_{(\ringO_E\otimes\Zp)}
\!{\rm T}_p(A_0^t[\wp^{\infty}]_{\rm et}),\,
\widehat{\bbG_m}
\right)
\]
\end{itemize}
\end{prop}

\proof The statement (i) is a consequence of Tate's theorem on
endomorphisms of abelian varieties over finite field, see \cite{Tate-Bourb}.
The statement (ii) follows from (i).
The statement (iii) is immediate from the displayed canonical 
isomorphism in Thm.\ \ref{STHMV}.
It remains to prove (iv).
\medbreak

By Thm.\ \ref{localrig} and Thm.\ \ref{STHMV}, 
we know that \(Z\) is a formal subtorus of the
formal torus 
\[
\calM_{E,\calL,\calL^+,n}^{/x_0}=\prod_{\wp\in \Sigma_{E,p}}
\underline{\rm Hom}_{\Zp}
\left({\rm T}_p(A_0[\wp^{\infty}]_{\rm et})\otimes_{(\ringO_E\otimes\Zp)}
\!{\rm T}_p(A_0^t[\wp^{\infty}]_{\rm et}),\,
\widehat{\bbG_m}
\right)\,.
\]
Let \(X_*(Z)\) be the group of formal cocharacters of the formal torus 
\(Z\). We know that \(X_*(Z)\) is a \(\Zp\)-submodule of
the cocharacter group
\[
\prod_{\wp\in\Sigma_{E,p}}
\left({\rm T}_p(A_0[\wp^{\infty}]_{\rm et})\right)^{\vee}
\otimes_{(\ringO_E\otimes\Zp)}
\!\left({\rm T}_p(A_0^t[\wp^{\infty}]_{\rm et})\right)^{\vee}
\]
of \(\calM_{E,\calL,\calL^+,n}^{/x_0}\), which is co-torsion free.
Moreover \(X_*(Z)\) is stable under the action of \(U_{x_0}\).
Denote by \(\,\ringO\,\) the closed subring of 
\(\,\prod_{\wp\in\Sigma_{E,p}}\!\ringO_{\wp}\,\)
generated by the image of the projection \({\rm pr}_1\) in (ii).
Since the image of \({\rm pr_1}\) is an open subgroup of 
\(\,\prod_{\wp\in\Sigma_{E,p}}\!\ringO_{\wp}^{\times}\),
the subring \(\,\ringO\,\) of 
\(\,\prod_{\wp\in\Sigma_{E,p}}\!\ringO_{\wp}\,\)
is an order of \(\prod_{\wp\in\Sigma_{E,p}}\!\ringO_{\wp}\).
So \(X_*(Z)\otimes \bbQ\) is stable under the action of 
\(\,\prod_{\wp\in\Sigma_{E,p}}\! E_{\wp}\). 
It follows that there exists a subset \(S\subset \Sigma_{E,p}\)
such that 
\(\,X_*(Z)\otimes_{\Zp}\Qp\,\) is equal to
\[
\prod_{\wp\in S}
\left({\rm T}_p(A_0[\wp^{\infty}]_{\rm et})\right)^{\vee}
\otimes_{\ringO_E\otimes\Zp}
\!\left({\rm T}_p(A_0^t[\wp^{\infty}]_{\rm et})\right)^{\vee}
\,.
\]
Since \(X_*(Z)\) is a co-torsion free \(\Zp\)-submodule
of 
\[\prod_{\wp\in\Sigma_{E,p}}
\left({\rm T}_p(A_0[\wp^{\infty}]_{\rm et})\right)^{\vee}
\otimes_{(\ringO_E\otimes\Zp)}
\!\left({\rm T}_p(A_0^t[\wp^{\infty}]_{\rm et})\right)^{\vee}\,,
\]
we see that \(X_*(Z)\) is equal to 
\(\left(\prod_{\wp\in S}
\left({\rm T}_p(A_0[\wp^{\infty}]_{\rm et})\right)^{\vee}
\otimes_{\ringO_E\otimes\Zp}
\!\left({\rm T}_p(A_0^t[\wp^{\infty}]_{\rm et})\right)^{\vee}
\right)
\).
\hfill\qed
\medbreak

\begin{cor}\label{cor:HMVlocalorbit}
Let \(x_0=[(A_0,\iota_0,\lambda_0,\eta_0)]
\in\calM_{E,\calL,\calL^+,n}(\Fpbar)\) be an ordinary \(\Fpbar\)-point of 
the Hilbert modular variety \(\calM_{E,\calL,\calL^+,n}\) as in 
\ref{HMVlocalaction}.  
Let \(Z\) be a reduced closed subscheme of \(\calM_{E,\calL,\calL^+,n}\)
such that \(x_0\in Z(\Fpbar)\). Assume that \(Z\) 
is stable under all 
\(\SL_2(\bbA_f^{(p)})\)-Hecke correspondences on
\(\calM_{E,\calL,\calL^+,n}(\Fpbar)\).
Then there exists a subset \(S_{x_0}\) of the set \(\Sigma_{E,p}\)
of prime ideals of \(\ringO_E\) containing \(p\)
such that 
\[
Z^{/x_0}=\prod_{\wp\in S}
\underline{\rm Hom}_{\Zp}
\left({\rm T}_p(A_0[\wp^{\infty}]_{\rm et})\otimes_{(\ringO_E\otimes\Zp)}
\!{\rm T}_p(A_0^t[\wp^{\infty}]_{\rm et}),\,
\widehat{\bbG_m}
\right)\,.
\]
Here \(Z^{/x_0}\) is the formal completion of 
\(Z\) at the closed point \(x_0\).
\end{cor}

\proof
Notation as in \ref{HMVlocalaction}. Recall that 
\(K=\End_{\ringO_E}^0(A_0)\).
Denote by \({\rm U}_K\) the unitary group attached to \(K\);
\({\rm U}_K\) is a linear algebraic group over \(\bbQ\)
such that \({\rm U}_K(\bbQ)=\{u\in K^{\times}\,\mid\,u\cdot \bar{u}=1\}\).
By \ref{HMVlocalaction} {\rm (i)}, \({\rm U}_K(\Qp)\) is isomorphic
to \((E\otimes\Qp)^{\times}\).  
Denote by \({\rm U}_K(\Zp)\) the compact open subgroup 
of \({\rm U}_K(\Qp)\) corresponding to the subgroup 
\((\ringO_E\otimes\Zp)^{\times}
\subset (\ringO_E\otimes\Qp)^{\times}\).  
This group \({\rm U}_K(\Zp)\) is 
isomorphic to the group \(H_{x_0}\) in \ref{HMVlocalaction} (ii),
via the projection to the first factor in the
displayed formula in \ref{HMVlocalaction} (i).
We have a natural action of \({\rm U}_K(\Zp)\) on
\[
\,\rm{Def}((A_0,\iota_0,\lambda_0)[p^{\infty}])/\Fpbar)
\cong
\calM_{E,\calL,\calL^+,n}^{/x_0}
\]
by the deformation of the deformation functor
\(\rm{Def}((A_0,\iota_0,\lambda_0)[p^{\infty}])\)
\smallbreak

Denote by \({\rm U}_K(\bbZ_{(p)})\) the subgroup 
\({\rm U}_K(\bbQ)\cap {\rm U}_K(\Zp)\) of \({\rm U}_K(\bbQ)\);
in other words \({\rm U}_K(\bbZ_{(p)})\)
consisting of all elements \(u\in {\rm U}_K(\bbQ)\) such that
\(u\) induces an automorphism of \(A_0[p^{\infty}]\).
Since \(Z\) is stable under all \(\SL_2(\bbA_f^{(p)})\)-Hecke
correspondences, the formal completion
\(Z^{/x_0}\) at \(x_0\) of 
the subvariety \(Z\subset \calM_{E,\calL,\calL^+,n}^{/x_0}\)
is stable under the natural action of the subgroup 
\({\rm U}_K(\bbZ_{(p)})\) of \({\rm U}_K(\bbQ)\).
By the weak approximation theorem for linear algebraic groups
(see \cite{PR}, 7.3, Theorem 7.7 on page 415),
\({\rm U}_K(\bbZ_{(p)})\) is \(p\)-adically dense in
\({\rm U}_K(\Zp)\).
So \(Z^{/x_0}\subset\calM_{E,\calL,\calL^+,n}^{/x_0}\) 
is stable under the action of \({\rm U}_K(\Zp)\)
by continuity.
We conclude the proof by invoking \ref{HMVlocalaction} (iii)
and (iv).
\hfill\qed
\bigbreak

\subsection{\bf Exercise.}
Let \((A,\iota)\) be an \(\ringO_E\)-linear abelian variety of
HB-type over a perfect field \(K\supset \Fp\).
Show that \({\rm M}_p((A,\iota)[p^{\infty}])\) is 
a free \((\ringO_E\otimes_{\bbZ}\Zp)\)-module of rank two.

\smallbreak

\subsection{\bf Exercise.}
Let \([x=(A,\iota,\lambda,\eta)]\in\calM_{E,\calL,\calL^+,n}(k)\)
be a geometric point of a Hilbert modular variety
\(\calM_{E,\calL,\calL^+,n}\), where \(k\supset\Fp\) is an
algebraically closed field.
Assume that \(\Lie(A/k)\) is a free \((\ringO_E\otimes_{\bbZ}k)\)-module
of rank one.  Show that \(\calM_{E,\calL,\calL^+,n}\) is smooth 
at \(x\) over \(k\).

\smallbreak

\subsection{\bf Exercise.}
Let \(k\supset \Fp\) be an algebraically closed field.
Assume that \(p\) is unramified in \(E\), i.e.\ \(E\otimes_{\bbZ}\Zp\)
is a product of unramified extension of \(\Qp\).
Show that \(\Lie(A/k)\) is a free \((\ringO_E\otimes_{\bbZ}k)\)-module
of rank one for every geometric point 
\((A,\iota,\lambda,\eta)]\in\calM_{E,\calL,\calL^+,n}(k)\).

\smallbreak

\subsection{\bf Exercise.}
Give an example of a geometric point 
\([x=(A,\iota,\lambda,\eta)]\in\calM_{E,\calL,\calL^+,n}(k)\)
such that \(\Lie(A/k)\) is not a free \((\ringO_E\otimes_{\bbZ}k)\)-module
of rank one.


\section{Deformations of $p$-divisible groups to $a \leq 1$}\label{7}
Main references: \cite{AJdJ.FO}, \cite{FO-NP}\\
In this section we will prove and use the following 
rather technical result.

\begin{thm}\label{Defoa1} \Th \ 
{\rm ({\bf Deformation to $a \leq 1$.})} {\it Let $X_0$ be a 
$p$-divisible group over a field $K$. 
There exists an integral scheme $S$, 
a point $0 \in S(K)$ and a $p$-divisible group $\cX \to S$ 
such that the fiber $\cX_0$ is isomorphic with $X_0$, 
and for the generic point 
$\eta \in S$ we have:
$$
\cN(X_0) = \cN(X_{\eta}) \quad\mbox{\it and}\quad a(X_{\eta}) \leq 1.
$$}
\end{thm}
See \cite{AJdJ.FO}, 5.12 and \cite{FO-NP}, 2.8. 

Note that if $X_0$ is ordinary (i.e. every  slope of $\cN(X_0)$ is 
either $1$ or $0$), there is not much to prove: 
$a(X_0) = 0 = a(X_{\eta})$; 
if however $X_0$ is not ordinary, the theorem says something non-trivial 
and in that case we end with  $a(X_{\eta}) = 1$.

At the end of this section we discuss the quasi-polarized case. 

\subsection{} {\it In this section we prove}  Theorem \ref{Defoa1} 
{\it in case $X_0$ is simple.} 
Surprisingly, this is the most difficult step. 
We will see, in Section 9, that once we have the theorem 
in this special case,  \ref{Defoa1} and  \ref{Defopa1} will follow 
without much trouble.

The proof (and the only one I know) of this special case given here is a 
combination of general theory, and a  computation. 
We start with one of the tools.

\begin{thm}\label{Purity} \Bl \  
{\rm ({\bf Purity  of the Newton polygon stratification.})}
Let $S$ be an integral scheme, and let $X \to S$ be a 
$p$-divisible group. 
Let $\gamma = \cN(X_{\eta})$ be the Newton polygon of the generic fiber. 
Let $S \supset D =  S_{\not=\gamma}
 := \{s \mid \cN(A_s) \precneqq \gamma\}$ {\rm (Note  that $D$ is closed 
in $S$ by Grothendieck-Katz.)} Then either $D$ is empty or 
${\rm codim}(D \subset S) = 1.$ 
\end{thm}
\n
I know two proofs of this theorem, and both proofs are non-trivial. 
See \cite{AJdJ.FO}, Th. 4.1. Also see \cite{Vasiu-bound}, th. 6.1; 
this second proof of Purity was analyzed in \cite{FO-Purity}.

When this result was first mentioned, it met unbelief. Why? 
If you follow the proof by Katz, see \cite{Katz-Slope}, 2.3.2, 
you see that $D = S_{\not=\gamma} \subset S$ is given by ``many'' 
defining equations. {}From that point of view ``codimension one'' 
seems unlikely. In fact it is not known (to my knowledge) whether 
there exists a scheme structure on $D = S_{\not=\gamma}$ such that 
$(D,\cO_{D}) \subset S$ is a Cartier divisor (locally principal) (i.e.
locally complete a intersection, or locally a 
set-theoretic complete intersection).

\subsection{\bf The simple case, notation.} 
We follow \cite{AJdJ.FO}, \S 5, \S 6. 
In order to prove \ref{Defoa1} in case $X_0$ is simple we fix notations, 
to be used for the rest of this section. Let $m \geq n > 0$ be 
relatively prime integers. We will write $r = (m-1)(n-1)/2$. 
We write $\delta$ for the isoclinic Newton polygon with slope 
$m/(m+n)$ with multiplicity  $m+n$.\\

\vn
{\it We define the $p$-divisible group} $H_{m,n}$, 
see \cite{AJdJ.FO}, 5.3. 
Let $K$ be a perfect field, $W = W_{\infty}(K)$. 
Consider symbols $\{e_i \mid i \in \ZZ_{\geq 0}\}$. Define 
$$
\cF{\cdot}e_i = e_{i+m}, \quad \cV{\cdot}e_{i+ n}, \quad p{\cdot}e_i 
= e_{i_m+n}, \quad M_{m,n} := \oplus_{0 \leq i < m+n} \  W.{\cdot}e_i.
$$
This is a Dieudonn\'e module, and the $p$-divisible group determined by 
this we write as $H_{m,n}$. Note that 
if $K \supset \FF_{p^{m+n}}$ 
then $D:= \End^0(H_{m,n})$ is a division algebra of rank $(m+n)^2$ over 
$\QQ_p$, and $R:= \End(H_{m,n})$ is the maximal order in
$\End^0(H_{m,n})$; note that $\pi \in \End(H_{m,n})$ 
defined by $\pi(e_i) = e_{i+1}$ indeed is 
an endomorphism,  it is a ''uniformizer`` in $R$, and $R$ is a 
``non-commutative discrete valuation ring''.\\
{\bf Remark.} We see that $H_{m,n}$ is defined over $\FF_p$; 
for any field $L$ 
we will write $H_{m,n}$ instead of $H_{m,n} \otimes L$ 
if no confusion can occur.\\
{\bf Remark.} The $p$-divisible group  $H_{m,n}$ is the 
``minimal $p$-divisible group'' with Newton polygon equal 
to $\delta$. For properties of minimal $p$-divisible groups 
see \cite{FO-Minp}. Such groups are of importance in 
understanding various stratifications of $\cA_g$.

\vn
We want to understand all   $p$-divisible groups isogenous with 
$H := H_{m,n}$ ($m$ and $n$ will remain fixed).

\begin{lemma} \ \Bl \ Work over a perfect field $K$. For every 
$X \sim H$ there is an isogeny $\va: H \to X$ of degree $p^r$.
\end{lemma}
\n
A proof of this lemma is not difficult. {\bf Exercise:}  give a proof.

\subsection{} {\bf Construction.} For every scheme $S$ over 
$\FF_p$ consider
$$
S \quad\mapsto\quad \{(\va,X) \mid \va: H \times S \to X, \ \ 
\degree(\va) = p^r\}.
$$ 
This functor is representable; denote the representing object by 
$(T = T_{m,n}, H_T \to \cG) \to \Spec(\FF_p)$. Note, using the lemma,  
that for any $X \sim H$ over a perfect field $K$  
there exists a point $x \in T(K)$ such that $X \cong \cG_x$.

\vn
{\bf Discussion.} The scheme $T$ constructed is just 
a special case of the theory of Rapoport-Zink spaces, 
see \cite{R.Z}, Th. 2.16. In the special 
case considered here, the space $T$ is a connected component of the 
reduction-mod-$p$ fiber of a Rapoport-Zink space. 

\begin{thm}\label{cat} \ \Th \ \ The scheme $T$ is geometrically 
irreducible of dimension $r$ over $\FF_p$. The set $T(a=1) \subset T$
 is open and dense in $T$.
\end{thm}
See \cite{AJdJ.FO}, Th. 5.11.
{\it Note that} \ref{Defoa1} {\it follows from this theorem 
in case} $X_0 \sim H_{m,n}$.
 We  focus on a proof of \ref{cat}. \\
{\bf Remark.} Suppose we have proved the case that $X_0 \sim H_{m,n}$. 
Then by duality we have $X_0^t \sim H_{m,n}^t = H_{n,m}$, and this case 
follows also. Hence it suffices to consider 
only the case $m \geq n > 0$. 

\vn
Notational remark. In this section we will not consider abelian varieties. 
The letters $A$, $B$, etc. in this section will not be used for 
abelian varieties. And then, semi-modules will only be considered 
in this section and in later sections these letters again 
will be used for abelian varieties.

\subsection{}
{\bf Definition.} We say that $A \subset \ZZ$ is a {\it semi-module} or
 more precisely, a $(m,n)$-semi-module, if 
\begin{itemize}
\item  $A$ is bounded from below, and if 
\item for every $x \in A$ we have $a+m, a+n \in A$. 
\end{itemize}
We write $A = \{a_1, a_2, \cdots\}$ with $a_j < a_{j+1} \ \ \forall j$.
We say that semi-modules $A, B$ are {\it equivalent} if there 
exists $t \in \ZZ$ such that $B = A + t := \{x + t \mid x \in A\}$.
\\
We say that $A$ is {\it normalized} if:

(1) $A \subset \ZZ_{\geq 0}$,

(2) $a_1 < \cdots < a_r \leq 2r$,

(3) $A = \{a_1, \cdots , a_r\} \cup [2r,\infty)$;\\
notation: $[y,\infty) := \ZZ_{\geq y}$.

\vn
Write $A^t = \ZZ \backslash (2r-1-A) = \{y \in\ZZ \mid 
2r-1-y \not\in A \}$.\\
{\bf Explanation.} For a semi-module $A$ the set $\ZZ \backslash A$ 
of course is a ``$(-m,-n)$-semimodule''. Hence $\{y \mid y \not\in A\}$
is a semi-module; then normalize.

\vn
{\bf Example.} Write $<0>$ for the  semi-module generated by $0$, 
i.e. consisting of all integers of the form $im + jn$ for  $i, j \geq 0$.

\vn 
{\bf Exercise.} 
(4) {\it Note that  $<0>$ indeed is normalized. 
Show that $2r-1 \not\in <0>$.}\\
(5) {\it Show: if $A$ is normalized then $A^t$ is normalized.}\\
(6) $A^{tt} = A$.\\
(7) {\it For every $B$ there is a unique normalized $A$ 
such that $A \sim B$.}\\
(8) {\it If  $A$ is  normalized, then: $A = <0> 
\quad\Longleftrightarrow\quad 0 \in A \quad\Longleftrightarrow\quad 2r-1 
\not\in A$. 
}

\subsection
{\bf Construction.} {\it Work over a perfect field. 
For every $X \sim H_{m,n}$ 
there exists a semi-module.}\\ 
An isogeny $X \to H$ gives an inclusion 
$$\DD(X) \quad\hookrightarrow\quad \DD(H) = M =  
\oplus_{0 \leq i < m+n} \  W.{\cdot}e_i.$$
Write $M^{(i)} = \pi^i{\cdot}M$. Define
$$B \quad:=\quad \{j \mid \DD(X) 
\cap M^{(j)} \not= \DD(X) \cap M^{(j+1)}\},$$
i.e. $B$ is the set of values where the filtration induced on 
$\DD(X)$ jumps. 
It is clear that $B$ is a semi-module. Let $A$ be the unique normalized  
semi-module equivalent to $B$.\\
{\bf Notation.} {\it The normalized semi-module constructed 
in this way will be called the}  
type {\it of $X$, denoted as ${\rm Type}(X)$ }.

\vn
We denote by $U_A \subset T$ the set where the semi-module $A$ 
is realized:
$$
U_A \quad=\quad \{t \in T \mid  {\rm Type}(\cG_t) = A\}.
$$

\begin{prop}\label{prop}
{\rm  (1)} $U_A \hookrightarrow T$ is locally closed, 
$T =  \sqcup _A \ U_A$.
\\
{\rm  (2)} $A = <0> \quad
\Longleftrightarrow\quad a(X) = 1$.
\\
{\rm  (3)} $U_{<0>}$ is geometrically irreducible and has dimension $r$.
\\
{\rm  (4)} If $A \not= <0>$ then every component of $U_A$ has 
dimension strictly less than $r$.
\end{prop}
For a proof see \cite{AJdJ.FO}, the proof on page 233, and 6.5 and 6.15. 
\B
\\
Note that a proof for this proposition is not very deep but somewhat 
involved (combinatorics and studying explicit equations). 

\subsection{}\label{def}\  \Bl \ Let $Y_0$ be any 
$p$-divisible group 
over a field $K$, of dimension $d$ and let $c$ 
be the dimension of $Y_0^t$. 
The  universal deformation space is isomorphic with 
$\Spf(K[[t_1, \cdots , t_{cd}]])$ and the generic fiber of that 
universal deformation is ordinary; in this case its Newton polygon 
$\rho$ has $c$ slopes equal to $1$ and $d$ slopes equal to $0$. 
See \ref{cor:defbt}.
See \cite{Ill}, 4.8, \cite{AJdJ.FO}, 5.15.

\subsection{\bf We prove \ref{cat},}  using \ref{Purity} and \ref{prop}. 
Note that the Zariski closure $(U_{<0>})^{\rm Zar} \subset T$ is 
geometrically irreducible, and has dimension $r$; we want to show 
equality $(U_{<0>})^{\rm Zar} = T$. Suppose there would be an 
irreducible component $T'$ of $T$ not contained in  
$(U_{<0>})^{\rm Zar}$. 
By \ref{prop} (3) and (4) we see that $\dim(T') < r$. Let $y \in T'$, 
with corresponding $p$-divisible group $Y_0$. 

Consider the formal completion $T^{/y}$ of $T$ at $y$.  
Write $D = {\rm Def}(Y_0)$ for the  universal deformation space of $Y_0$. 
The moduli map $T^{/y} \to  D = {\rm Def}(Y_0)$ is an immersion, 
see \cite{AJdJ.FO}, 5.19. Let  
 $T'' \subset D$ be the image of $(T')^{/y}$ in $D$; we conclude that 
no irreducible component of $T''$ is contained in 
any irreducible component of the image of $T^{/y} \to D$ in $D$, 
i.e. every component of $T''$ is an component of $\cW_{\delta}(D)$. 
Clearly $\dim(T')  = \dim(T'') < r$. 

\vn
{\bf Obvious, but crucial observation.} 
\\
\begin{picture}(100,70)(0,0)
\put(0,50){\it Consider the graph of all Newton polygons}
\put(40,30){$\zeta \quad\mbox{with}\quad  \delta \prec \zeta \prec \rho.$}
\put(-0,10){\it The longest path in this graph has 
length $\leq mn-r$.}
\put(320,0){\line(1,0){40}}
\put(360,0){\line(1,1){60}}
\put(320,0){\line(5,3){100}}
\put(320,0){\line(1,1){60}}
\put(380,60){\line(1,0){40}}

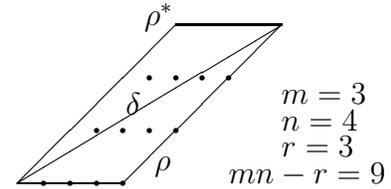
\put(340,0){\circle*{2}}
\put(360,0){\circle*{2}}

\put(360,20){\circle*{2}}
\put(380,20){\circle*{2}}

\put(380,40){\circle*{2}}
\put(400,40){\circle*{2}}

\put(330,0){\circle*{2}}
\put(350,0){\circle*{2}}

\put(350,20){\circle*{2}}
\put(370,20){\circle*{2}}

\put(370,40){\circle*{2}}
\put(390,40){\circle*{2}}

\put(372,5){$\rho$}
\put(361,26){$\delta$}
\put(368,60){$\rho^{\ast}$}

\put(420,30){$m=3$}
\put(420,20){$n=4$}
\put(420,10){$r=3$}
\put(400,0){$mn-r=9$}

\end{picture}

\vn
{\bf Proof.} Consider the Newton polygon $\rho$, 
in this case given by $n$ 
slopes equal to $0$ and $m$ slopes equal to $1$. Note that $\gcd(m,n)=1$, 
hence the Newton polygon $\delta$ does not contain integral points except 
beginning and end point. Consider the interior of the parallelogram 
given by $\rho$ and by $\rho^{\ast}$, the upper convex polygon given by: 
first $m$ slopes equal to $1$ and then $n$ slopes equal to $0$. 
The number of interior points of this parallelogram equals $(m-1)(n-1)$.
Half of these are above $\delta$, and half of these are below $\delta$. 
Write $\delta \precneqq (i,j)$ for the property ``$(i,j)$ is 
strictly below $\delta$'', and $(i,j) \prec \rho$ for ``$(i,j)$ is 
upon or above $\rho$''. We see:

$
\#\left(\{(i,j) \mid \delta \precneqq (i,j) \prec \rho\}\right) 
\quad=\quad  (m-1)(n-1)/2 + (m+n-1)\quad=\quad  mn-r.$ \quad\qed
\n
We used the fact: If $\zeta_1 \precneqq \zeta_2$, 
then there is an integral 
point on $\zeta_2$ strictly below $\zeta_1$. One can even show that 
all maximal chains of Newton polygons in the fact above have the 
same length, and in fact equal to 
$$\#\left(\{(i,j) \mid \delta \precneqq (i,j) \prec \rho\}\right).$$ 
  
\vn
As $\dim({\rm Def}(Y_0)) = mn$ this observation implies by Purity, 
see \ref{Purity}, that every irreducible component of $\cW_{\delta}(D)$ 
had dimension at least $r$. 
This is a contradiction with the assumption of 
the existence of $T'$, i.e.  $\dim(T')  = \dim(T'') < r$. 
Hence $(U_{<0>})^{\rm Zar} = T$. This proves \ref{cat}. 
Hence we have proved \ref{Defoa1} in case $X_0$ is isogenous with 
$H_{m,n}$. \hfill\qed
\medbreak

\begin{thm}\label{Defopa1} \Th \ 
{\rm ({\bf Deformation to $a \leq 1$ in 
the principally quasi-polarized case.})} 
{\it Let $X_0$ be a $p$-divisible group over a field $K$ 
with a {\it principal} 
quasi-polarization $\lambda_0: X_0 \to X_0^t$. There exists 
an integral scheme $S$, a point $0 \in S(K)$ and a 
principally quasi-polarized $p$-divisible group $(\cX,\lambda) \to S$ 
such that there is an isomorphism 
$(X_0,\lambda_0) \cong (\cX,\lambda)_0$, 
and for the generic point $\eta \in S$ we have:
$$
\cN(X_0) = \cN(X_{\eta}) \quad\mbox{\it and}\quad a(X_{\eta}) \leq 1.
$$}
\end{thm}
See \cite{AJdJ.FO}, 5.12 and \cite{FO-NP}, 3.10.
\B

\begin{cor}\label{DefAVtoa1}  \Th \  
{\rm ({\bf Deformation to 
$a \leq 1$ in the case of  principally polarized abelian varieties.})} 
{\it Let $(A_0,\lambda_0)$ be a 
principally polarized abelian variety over $K$.  
There exists an integral scheme $S$, a point $0 \in S(K)$ 
and a principally polarized abelian scheme $(A,\lambda) \to S$ 
such that there is an isomorphism 
$(A_0,\lambda_0) \cong (A,\lambda)_0$, 
and for the generic point $\eta \in S$ we have:
$$
\cN(A_0) = \cN(A_{\eta}) \quad\mbox{\it and}\quad a(X_{\eta}) \leq 1.
$$}
\end{cor}

\subsection{\bf The non-principally polarized case.} 
Note that the analog of the theorem 
and of the corollary is not correct in general in the 
{\it non-principally polarized} case. Here is an example, 
see \cite{KOIII}, 6.10, and also see \cite{Li.FO}, 12.4 and 12.5 
where more examples are given. {\it Consider $g=3$, let $\sigma$ 
be the supersingular Newton polygon;  
it can be proved that for any 
$x \in \cW_{\sigma}(\cA_{3,p})$ we have $a(A_x) \geq 2$.}

\vn
We will show that for $\xi_1 \prec \xi_2$ we have in the principally 
polarized case: 
$$\cW_{\xi_1}^0(\cA_{g,1}) =: W_{\xi_1}^0 \quad\subset\quad  
(W_{\xi_2}^0)^{\rm Zar}  = W_{\xi_2} := \cW_{\xi_2}(\cA_{g,1}).$$
In the non-principally polarized case this inclusion and the equality 
$(W_{\xi_2}^0)^{\rm Zar}  = W_{\xi_2}$
does not hold in general as is showed by the following example. 
Let $g=3$, and $\xi_1 = \sigma$ the supersingular Newton polygon, 
and $\xi_2 = (2,1) + (1,2)$. Clearly $\xi_1 \prec \xi_2$. 
By \cite{KOIII}, 6.10, 
there is a component of $\cW_{\sigma}(\cA_{g,p^2})$ of dimension $3$;
more generally  see  \cite{Li.FO}, Th. 10.5 (ii) 
for the case of $\cW_{\sigma}(\cA_{g,p^{[(g-1)^2/2]}})$ 
and components of dimension equal to $g(g-1)/2$. 
As the $p$-rank $0$ locus in $\cA_g$ has pure dimension equal 
to $g(g+1)/2  + (f-g)= g(g-1)/2$, see \cite{Norman.FO}, Th. 4.1, 
this shows the existence of a polarized supersingular abelian variety 
(of dimension 3, respectively of any dimension at least 3)  which 
cannot be deformed to a non-supersingular abelian variety with 
$p$-rank equal to zero.

Many more examples where $(W_{\xi}^0)^{\rm Zar}  \not= W_{\xi}$  
follow from \cite{FO-Dimleaves}, Section 3.


\section{Proof of the Grothendieck conjecture}\label{8}
\n
Main reference: \cite{FO-NP}.

\subsection{\bf Definition.}\ \Ex \  
Let $X$ be a $p$-divisible group over a 
base $S$. We say that 
$0 =  X^{(0)} \subset X^{(1)} \subset \cdots \subset X^{(s)} = X$ is  
{\it the slope filtration} of $X$ if $Y_i := X^{(i)}/ X^{(i-1)}$ for 
$1 \leq i \leq s$ is isoclinic of slope $\tau_i$ with 
$\tau_1 > \tau_2 > \cdots > \tau_s$. \\
{\bf Remarks.} Clearly, if a slope filtration exists, it is unique.

From the Dieudonn\'e-Manin classification it follows that 
the slope filtration on $X$ exists if $K$ is perfect.

By Grothendieck and Zink we know that for every $p$-divisible group over 
any field $K$ the slope filtration exists, see \cite{Zink-Slopef}, Coroll. 13.

In general for a $p$-divisible group $X \to S$ over a base 
a slope filtration on $X/S$ does not exists. 
Even if the Newton polygon is constant in a family, 
in general the slope filtration does not exist.

\subsection{\bf Definition.} We say that  
$0 =  X^{(0)} \subset X^{(1)} \subset \cdots \subset X^{(s)} = X$ 
is a {\it maximal filtration} of $X \to S$  
if every geometric fiber of $Y^{(i)} := X^{(i)}/ X^{(i-1)}$ for 
$1 \leq i \leq s$ 
is {\it simple} and isoclinic  of slope $\tau_i$ with 
$\tau_1 \geq  \tau_2 \geq  \cdots \geq \tau_s$. \\
{\bf Lemma.}\ \Bl \  {\it For every $X$ over $k$ 
a maximal filtration exists.}\\
See \cite{FO-NP}, 2.2.  \B

\subsection{\bf Lemma.} \ \Bl \  {\it Let $\{X_0^{(i)}\}$ be a 
$p$-divisible group $X_0$ with maximal filtration over $k$. 
There exists an integral scheme $S$ and a $p$-divisible group $X/S$ 
with a maximal filtration $\{X^{(i)}\} \to S$ and a closed point 
$0 \in S(k)$  such that $\cN(Y^{(i)})$ is constant for $1 \leq i \leq s$, 
such that    $\{X^{(i)}\}_0 = \{X_0^{(i)}\}$   
and such that for the generic point $\eta \in S$ we have 
$a(X_{\eta}) \leq 1$.} \\
See \cite{FO-NP}, Section 2. A proof of this lemma uses Theorem \ref{cat}. \B

\vn
In Section \ref{7} we proved \ref{cat}, and obtained as corollary 
\ref{Defoa1} in the case of a simple $p$-divisible group. 
From the previous lemma we  
derive a proof for  Theorem \ref{Defoa1}.

\subsection{}{\bf Definition.} We say that $X_0$ over a field $K$ is a 
{\it specialization} of $X_{\eta}$ over a field $L$ if there exists 
an integral scheme $S \to \Spec(K)$, a $k$-rational point $0 \in S(K)$, 
and $\cX \to S$ such that $X_0 = \cX_0$, and for the generic point 
$\eta \in S$ we have $L = K(\eta)$ and $X_{\eta} = \cX_{\eta}$. \\
This can be used for $p$-divisible groups, for abelian schemes, etc.

\subsection{}{\bf Proposition.} {\it Let $X_0$ be a specialization of 
$X_{\eta} = Y_0$, and let $Y_0$  be a specialization of $Y_{\rho}$. 
Then $X_0$ is a specialization of $Y_{\rho}$.}

\vn
Using Theorem \ref{G} and    Theorem \ref{Defoa1} by the proposition above  
we derive a proof for the Grothendieck Conjecture Theorem \ref{GC}. 
\hfill \qed

\subsection{\bf Corollary} of Theorem \ref{GC}. 
{\it  Let $X_0$ be a $p$-divisible group, $\beta = \cN(X_0)$. 
Every component of the locus $\cW_{\beta}({\rm Def}(X_0))$ 
has dimension $\diam(\beta)$.}  \B

\subsection{\bf Definition.} Let $(X,\lambda)$ be a principally 
polarized $p$-divisible group over $S$. We say that 
$0 =  X^{(0)} \subset X^{(1)} \subset \cdots \subset X^{(s)} = X$ is 
{\it a maximal symplectic filtration} of $(X,\lambda)$ if:
\begin{itemize}
\item every geometric fiber of $Y^{(i)} := X^{(i)}/ X^{(i-1)}$ for 
$1 \leq i \leq s$ 
is simple of slope $\tau_i$,
\item $\tau_1 \geq  \tau_2 \geq  \cdots \geq \tau_s$, and
\item  $\lambda: X \to X^t$ induces an isomorphism
$$
\lambda_i: Y^{(i)} \to (Y^{(s+1-i)})^t \quad\mbox{for}
\quad 0 < i \leq (s+1)/2.
$$
\end{itemize}

\subsection{\bf Lemma.} {\it For every principally polarized  
$(X,\lambda)$ over 
$k$ there exists a maximal symplectic filtration.}\\
See \cite{FO-NP}, 3.5. \B

\subsection{} Using this definition, and this lemma 
we show the principally polarized analog \ref{DefAVtoa1} of  \ref{Defopa1}, 
see \cite{FO-NP}, Section 3. Hence Corollary \ref{DefAVtoa1} 
follows. Using \ref{DefAVtoa1} and Theorem \ref{pG} we derive a proof for:

\subsection{\bf Theorem} (an analog of the Grothendieck conjecture). 
{\it Let $K \supset \FF_p$. Let $(X_0,\lambda_0)$ be a} principally 
{\it quasi-polarized $p$-divisible group over $K$. 
We write $\cN(X_0) = \xi$ for its Newton polygon.  
Suppose given a Newton polygon 
$\zeta$ ``below" $\xi$, i.e.  $\xi \prec \zeta$. 
There exists a deformation 
$(X_{\eta},\lambda_{\eta})$ of $(X_0,\lambda_0)$ 
such that ${\cal N}({\cal X}_{\eta}) = \zeta$.}
\B

\subsection{\bf Corollary.}  {\it Let $K \supset \FF_p$. 
Let $(A_0,\lambda_0)$ be a} 
principally {\it polarized abelian variety over $K$. 
We write $\cN(A_0) = \xi$ 
for its Newton polygon.  
Suppose given a Newton polygon 
$\zeta$ ``below" $\xi$, i.e.  $\xi \prec \zeta$. 
There exists a deformation $(A_{\eta},\lambda_{\eta})$ of 
$(A_0,\lambda_0)$ 
such that ${\cal N}({\cal X}_{\eta}) = \zeta$.}
\B

\subsection{\bf Corollary.}\label{dimp} {\it Let $\xi$ be 
a symmetric Newton polygon. Every component of the stratum 
$W_{\xi} = \cW_{\xi}(\cA_{g,1})$ has dimension equal to $\triangle(\xi)$.}
\B


\section{Proof of the density of ordinary Hecke orbits}\label{9}

In this section we give a proof of Theorem \ref{dens}
on density of ordinary Hecke orbits, restated as
Theorem \ref{dens_siegel} below.
To establish Thm.\ \ref{dens}, we need the analogous statement for
a Hilbert modular variety; see \ref{dens_hilb} for the precise statement.
\smallbreak

Here is a list of tools we will use; many have been explained in 
previous sections.
\begin{itemize}
\item[(i)] Serre-Tate coordinates, see \S\ref{2}.

\item[(ii)] Local stabilizer principle, see \ref{lsp} and \ref{lsphilb}.

\item[(iii)] Local rigidity for group actions on formal tori, 
see \ref{localrig}.

\item[(iv)] Consequence of EO stratification, see \ref{EOconseq}.

\item[(iv)] Hilbert trick, see \ref{thm_hilb}. 
\end{itemize}
\smallbreak
The logical structure of the proof of Theorem \ref{dens} is as follows. 
We first prove the density of ordinary Hecke orbits on
Hilbert modular varieties. Then we use the Hilbert trick to show that
the Zariski closure of any prime-to-\(p\) Hecke orbit on 
\(\calA_{g,1,n}\) contains a hypersymmetric point; here 
we mean that the underlying abelian variety
is ordinary and hypersymmetric, i.e.\ it is isogenous to
\(B\times\cdots\times B\), where \(B\) is a simple ordinary hypersymmetric abelian
variety over \(\bbF\); see \cite{Chai.FO-Hypersymmetric}.
Finally we use the local stabilizer principle and the local rigidity
to conclude the proof of \ref{dens}
\medbreak

The Hilbert trick is based on the following observation.
Given an ordinary point 
\(x=[(A_x,\lambda_x,\eta_x)]\in\calA_{g,1,n}(\Fpbar)\),
the prime-to-\(p\) Hecke orbit of \(x\) contains, up to a
possibly inseparable isogeny correspondence, 
the (image of) the prime-to-\(p\) Hecke orbit of a point
\(h=[(A_y,\iota_y,\lambda_y,\eta_y)]\) of a Hilbert modular
variety \(\calM_{E,\calL,\calL^+,m}\) such that \(A_y\) is isogenous
to \(A_x\), because \(\End^0(A_x)\) contains a product
\(E=F_1\times\cdots\times F_r\) of totally real field with
\([E:\bbQ]=g\).
So if we can establish the density of the prime-to-\(p\) Hecke orbit
of \(y\) in \(\calM_{E,\calL,\calL^{+},m}\), then we know that
the Zariski closure of the prime-to-\(p\) Hecke orbit of \(x\)
contains the image of the Hilbert modular variety 
\(\calM_{E,\calL,\calL^{+},m}\) in \(\calA_{g,1,n}\) under
a finite isogeny correspondence, i.e.\
a scheme \(T\) over \(\Fp\) 
and finite \(\Fp\)-morphisms \(g: T\to \calM_{E,\calL,\calL^{+},m}\) and 
\(f: T\to \calA_{g,1,n}\) such that the pull-back by \(g\) of the 
universal abelian scheme over \(\calM_{E,\calL,\calL^{+},m}\)
is isogenous to the pull-back by \(f\) of the universal
abelian scheme over \(\calA_{g,1,n}\).
Since \(\calM_{E,\calL,\calL^+,m}\) contains ordinary hypersymmetric
points, \(\,\left(\calH^{(p)}_{\rm Sp}(x)\right)^{\rm Zar}\,\)
also contains an ordinary hypersymmetric point.
Then the linearization method afforded by the combination of the
\emph{local stabilizer principle} and
the \emph{local rigidity}
implies that the dimension of
\(\,\left(\calH^{(p)}_{\rm Sp}(x)\right)^{\rm Zar}\,\)
is equal to \(g(g+1)/2\), hence
\(\left(\calH^{(p)}_{\rm Sp}(x)\right)^{\rm Zar}=\calA_{g,1,n}\)
because \(\calA_{g,1,n}\) is geometrically irreducible, see \ref{irredA}.
\smallbreak

To prove the density of ordinary Hecke orbits on a Hilbert modular
variety, the linearization method is again crucial. 
Since a Hilbert modular variety \(\calM_{E,\calL,\calL^+,m}\)
is ``small'', there are only a
finite number of possibilities as to what (the formal completion of)
the Zariski closure of an ordinary Hecke orbit can be;
the possibilities are indexed by the set of all subsets of
prime ideals of \(\ringO_E\). To pin the number of possibilities
down to one, one can use either the consequence of 
EO-stratification that the Zariski closure
of any Hecke-invariant subvariety of a Hilbert modular variety 
contains a supersingular point, or de Jong's theorem on
extending homomorphisms between Barsotti-Tate groups.
We follow the first approach here, see 
\ref{rem:alt} and \cite[\S8]{Chai-Fam}
for the second approach.
\bigbreak

\begin{thm} \label{dens_siegel}
Let \(n\geq 3\) be an integer prime to \(p\).
Let \(x=[(A_x,\lambda_x,\eta_x)]\in \calA_{g,1,n}(\Fpbar)\)
such that \(A_x\) is ordinary.  
\begin{itemize}
\rmitem[(i)] The prime-to-\(p\) 
\(\Sp_{2g}(\bbA_f^{(p)})\)-Hecke orbit of \(x\) is
dense in the moduli space \(\calA_{g,1,n}\) over \(\,\Fpbar\,\), i.e.\
\[
\left(\cH^{(p)}_{\rm Sp}(x)\right)^{\rm Zar}=\calA_{g,1,n}\,.
\]

\rmitem[(ii)] The \(\Sp_{2g}(\Ql)\)-Hecke orbit of \(x\)
dense in the moduli space \(\calA_{g,1,n}\) over \(\,\Fpbar\,\), i.e.\
\[
\left(\cH_{\ell}^{\rm Sp}(x)\right)^{\rm Zar}=\calA_{g,1,n}\,.
\]
\end{itemize}
\end{thm}
\medbreak

\begin{thm} \label{dens_hilb}
Let \(n\geq 3\) be an integer prime to \(p\).
Let \(E=F_1\times\cdots\times F_r\), where \(F_1,\ldots,F_r\)
are totally real number fields.
Let \(\calL\) be an invertible \(\ringO_E\)-module,
and let \(\calL^+\) be a notion of positivity for \(\calL\).
Let \(y=[(A_y,\iota_y,\lambda_y,\eta_y)]
\in\calM_{E,\calL,\calL^{+},n}(\Fpbar)\) be a point of 
the Hilbert modular variety \(\calM_{E,\calL,\calL^+,n}\)
such that \(A_y\) is ordinary.
Then the \(\SL_2(E\otimes_{\bbQ}\bbA_f^{(p)})\)-Hecke orbit of
\(y\) on \(\calM_{E,\calL,\calL^+,n}\) is Zariski dense
in \(\calM_{E,\calL,\calL^+,n}\) over \(\Fpbar\).
\end{thm}

\begin{prop} \label{smooth}
Let \(n\geq 3\) be a integer prime to \(p\).
\begin{itemize}

\rmitem[(i)] Let \(x\in\calA_{g,1,n}(\Fpbar)\) be a closed point
of \(\calA_{g,1,n}\).  Let \(Z(x)\) be the Zariski closure
of the prime-to-\(p\) Hecke orbit \(\cH_{\Sp_{2g}}^{(p)}(x)\)
in \(\calA_{g,1,n}\) over \(\Fpbar\).  Then \(Z(x)\) is smooth at \(x\)
over \(\Fpbar\).

\rmitem[(ii)] Let \(y\in\calM_{E,\calL,\calL^+,n}(\Fpbar)\) be
a closed point of a Hilbert modular variety 
\(\calM_{F,\calL,\calL^+,n}\).  Let \(Z_F(y)\) be the
be the Zariski closure of the prime-to-\(p\) Hecke orbit 
\(\cH_{\SL_2}^{(p)}(y)\) on \(\calM_{F,\calL,\calL^+,n}\)
over \(\Fpbar\).  Then \(Z_F(y)\) is smooth at \(y\) over \(\Fpbar\).
\end{itemize}
\end{prop}

\proof We give the proof of (ii) here.  The proof of (i) is
similar and left to the reader.
\smallbreak

Because \(Z_F\) is reduced, there exists a dense open subset
\(U\subset Z_F\) which is smooth over \(\Fpbar\).
This open subset \(U\) must contain an element \(y'\) of the dense subset
\(\cH_{\SL_2}^{(p)}(y)\) of \(Z_F\), so \(Z_F\) is smooth over \(\Fpbar\)
at \(y'\). Since prime-to-\(p\) Hecke correspondences are defined
by schemes over 
\(\calM_{F,\calL,\calL^+,n}\times_{\Spec(\Fpbar)}\calM_{F,\calL,\calL^+,n}\)
such that both projections to \(\calM_{F,\calL,\calL^+,n}\) are \'etale,
\(Z_F\) is smooth at \(y\) as well.
\hfill\qed
\medbreak

\noindent{\bf Remark.}\enspace
(i) Prop.\ \ref{smooth} is an analog of the following well-known fact.
Let \(X\) be a reduced scheme over an algebraically closed field \(k\) on
which an algebraic group operates transitively.  
Then \(X\) is smooth over \(k\).
\smallbreak

(ii) The proof of Prop.\ \ref{smooth} also shows 
that all irreducible components
of \(Z(x)\) (resp.\ \(Z_F(y)\)) have the same dimension:
For any non-empty subset \(U_1\subset Z_F(y)\) and any
open subset \(W_1 \ni y\), there exist a
non-empty subset \(U_2\subset U_1\), an open subset \(W_2\ni y\)
and a non-empty \'etale correspondence between \(U_2\) and \(W_2\).
\medbreak

\begin{thm}\Bl\label{monconseq}
Let \(Z\) be a reduced closed subscheme of \(\calA_{g,1,n}\) over \(\Fpbar\)
such that no maximal point of \(Z\) is contained in the supersingular
locus of \(\calA_{g,1,n}\).
If \(Z\) is stable under all \(\Sp_{2g}(\Ql)\)-Hecke correspondences
on \(\calA_{g,1,n}\), then \(Z\) is stable under all
\(\Sp_{2g}(\bbA_f^{(p)})\)-Hecke correspondences.
\end{thm}
\smallbreak

\noindent{\bf Remark.} This is proved in \cite[Prop. 4.6]{Chai-Mon}. 
\bigbreak

\noindent{\scshape Local stabilizer principle}
\smallbreak

Let \(k\supset\Fp\) be an algebraically closed field.
Let \(Z\) be a reduced closed subscheme of \(\calA_{g,1,n}\) over
\(k\).
Let \(z=[(A_z,\lambda_z,\eta_z)]\in Z(k)\subset \calA_{g,1,n}(k)\) 
be a closed point of \(Z\).
Let \(*_z\) be the Rosati involution on \(\End^0(A_z)\).
Denote by \({\rm H}_z\) the unitary group attached to the
semisimple algebra with involution \((\End^0(A_z), *_z)\),
defined by
\[
{\rm H}_z(R)=\left\{x\in (\End^0(A_z)\otimes_{\bbQ}R)^{\times}\,
\mid\, x\cdot *_0(x)=*_0(x)\cdot x={\rm Id}_{A_z}
\right\}
\]
for any \(\bbQ\)-algebra \(R\).
Denote by \({\rm H}_z(\Zp)\) the subgroup of \({\rm H}_z(\Qp)\)
consisting of all elements \(x\in {\rm H}_z(\Qp)\) such that 
\(x\) induces an automorphism of 
\((A_z,\lambda_z)[p^{\infty}]\).
Denote by \({\rm H}_z(\bbZ_{(p)})\) the group 
\({\rm H}_z(\bbQ)\cap {\rm H}_z(\Zp)\), i.e.\ its elements
consists of all elements \(x\in {\rm H}_z(\bbQ)\) such that 
\(x\) induces an automorphism of 
\((A_z,\lambda_z)[p^{\infty}]\).
Note that the action of \({\rm H}_z(\Zp)\) on \(A_z[p^{\infty}]\)
makes \({\rm H}_z(\Zp)\) a subgroup of 
\({\rm Aut}((A_z,\lambda_z)[p^{\infty}])\).
Denote by \(\calA_{g,1,n}^{/z}\) (resp.\ \(Z^{/z}\)) the formal completion 
of \(\calA_{g,1,n}\) (resp.\ \(Z\)) at \(z\). 
The compact \(p\)-adic group
\({\rm Aut}((A_z,\lambda_z)[p^{\infty}])\)
operates naturally on the deformation space
\({\rm{Def}}\left((A_z,\lambda_z)[p^{\infty}]/k
\right)\).
So we have a natural action of
\({\rm Aut}((A_z,\lambda_z)[p^{\infty}])\)
on the formal scheme \(\calA_{g,1,n}^{/z}\) via
the canonical isomorphism
\[
\calA_{g,1,n}^{/z}
={\rm{Def}}\left((A_z,\lambda_z)/k\right)
\xrightarrow[\sim]{\mbox{\scriptsize{\rm Serre-Tate}}}
{\rm{Def}}\left((A_z,\lambda_z)[p^{\infty}])/k
\right)\,.
\]
\smallbreak

\begin{thm}[local stabilizer principle]\label{lsp}
Notation as above. Suppose that \(Z\) is stable under all 
\(\Sp_{2g}(\bbA_f^{(p)})\)-Hecke correspondences on \(\calA_{g,1,n}\).
Then the closed formal subscheme \(Z^{/z}\) in \(\calA_{g,1,n}^{/z}\)
is stable under the action of the subgroup
\({\rm H}_z(\Zp)\) of \({\rm Aut}((A_z,\lambda_z)[p^{\infty}])\).
\end{thm}

\proof
Consider the projective system 
\(\,\widetilde{\calA_{g,1}}=\varprojlim_m\calA_{g,1,m}\,\) over \(k\),
where \(m\) runs through all integers \(m\geq 1\) which are
prime to \(p\).
The pro-scheme \(\widetilde{\calA_{g,1}}\) classifies triples
\((A\to S, \lambda, \eta)\), where 
\(A\to S\) is an abelian scheme up to prime-to-\(p\) isogenies, 
\(\lambda\) is a principal polarization of \(A\to S\),
and 
\[
\eta: {\rm H}_1(A_z,\bbA_f^{(p)})\xrightarrow{\sim}
\underline{\rm H}_1(A/S, \bbA_f^{(p)})
\]
is a symplectic prime-to-\(p\) level structure.
Here we have used the first homology groups of \(A_z\) 
attached to the base point \(z\) to produce the 
standard representation representation of the symplectic group
\(\Sp_{2g}\).
Take \(S_z=\calA_{g,1,n}^{/z}\), let
\((\hat{A},\hat{\lambda})\to \calA_{g,1,n}^{/z}\) be the
restriction of the universal principally polarized abelian scheme
to \(\calA_{g,1,n}^{/z}\), and let
\(\hat{\eta}\) be the tautological prime-to-\(p\) level structure,
we get an \(S_z\)-point of the tower \(\widetilde{\calA_{g,1}}\)
that lifts \(S_z\hookrightarrow \calA_{g,1,n}\).
\smallbreak

Let \(\gamma\) be an element of \({\rm H}_z(\bbZ_{(p)})\).
Let \(\gamma_p\) (resp.\ \(\gamma^{(p)}\)) be the image of \(\gamma\)
in the local stabilizer subgroup 
\({\rm H}_z(\Zp)\subset {\rm Aut}((A_z,\lambda_z)[p^{\infty}])\)
(resp.\ in \({\rm H}_z(\bbA_f^{(p)})\).
From the definition of the action of 
\({\rm Aut}((A_z,\lambda_z)[p^{\infty}])\)
on \(\calA_{g,1,n}^{/z}\) we have a commutative diagram
\[
\xymatrix{
(\hat{A},\hat{\lambda})[p^{\infty}]
\ar[r]^{\quad f_{\gamma}[p^{\infty}]\quad}\ar[d]
&(\hat{A},\hat{\lambda})[p^{\infty}]\ar[d]
\\
\calA_{g,1,n}^{/z}\ar[r]^{u_{\gamma}}
&\calA_{g,1,n}^{/z}
}
\]
where \(u_{\gamma}\) is the action of \(\gamma_p\) on \(\calA_{g,1,n}^{/z}\)
and \(f_{\gamma}[p^{\infty}]\) is an isomorphism over \(u_{\gamma}\) whose
fiber over \(z\) is equal to \(\gamma_p\).
Since \(\gamma_p\) comes from a prime-to-\(p\) quasi-isogeny,
\(f_{\gamma}[p^{\infty}]\) extends to a prime-to-\(p\) quasi-isogeny 
\(f_{\gamma}\) over \(u_{\gamma}\), such that the diagram
\[
\xymatrix{
\hat{A}\ar[r]^{f_{\gamma}}\ar[d]
&\hat{A}\ar[d]
\\
\calA_{g,1,n}^{/z}\ar[r]^{u_{\gamma}}
&\calA_{g,1,n}^{/z}
}
\]
commutes and \(f_{\gamma}\) preserves the polarization \(\hat{\lambda}\).
Clearly the fiber of \(f_{\gamma}\) at \(z\) is equal to 
\(\gamma\) as a prime-to-\(p\) isogeny from \(A_z\) to itself.
From the definition of the action of the symplectic group
\(\Sp({\rm H}_1(A_z,\bbA_f^{(p)}), \langle\cdot,\cdot\rangle)\) one sees that
\(u_{\gamma}\) coincides with the action of 
\((\gamma^{(p)})^{-1}\) on 
\(\widetilde{\calA_{g,1}}\).
Since \(Z\) is stable under all 
\(\Sp_{2g}(\bbA_f^{(p)})\)-Hecke correspondences,
we conclude that \(Z^{/z}\) is stable under the action of \(u_{\gamma}\),
for every \(\gamma\in {\rm H}_z(\bbZ_{(p)})\).
\smallbreak

By the weak approximation theorem for linear algebraic groups
(see \cite{PR}, 7.3, Theorem 7.7 on page 415),
\({\rm H}_z(\bbZ_{(p)})\) is \(p\)-adically dense in 
\({\rm H}_z(\Zp)\). 
So \(Z^{/z}\) is stable under the action of \({\rm H}_z(\Zp)\)
by the continuity of the action of 
\({\rm Aut}((A_z,\lambda_z)[p^{\infty}])\).
\hfill\qed
\medbreak

\noindent{\bf Remark.}\enspace  
The group \({\rm H}_z(\bbZ_{(p)})\) can be thought of as
the ``stabilizer subgroup'' at \(z\) inside the family of
prime-to-\(p\) Hecke correspondences: Every element 
\(\gamma\in{\rm H}_z(\bbZ_{(p)})\) gives rise to a prime-to-\(p\)
Hecke correspondence with \(z\) as a fixed point.
\bigbreak

We set up notation for the local stabilizer principle for 
Hilbert modular varieties.
Let \(E=F_1\times\cdots\times F_r\), where \(F_1,\ldots,F_r\)
are totally real number fields.
Let \(\calL\) be an invertible \(\ringO_E\)-module, and 
let \(\calL^+\) be a notion of positivity for \(\calL\).
Let \(m\geq 3\) be a positive integer which is prime to \(p\).
Let \(Y\) be a reduced closed subscheme of \(\calM_{E,\calL,\calL^+,m}\) over 
\(\Fpbar\).
Let \(y=[(A_y,\iota_y,\lambda_y,\eta_y)]\in\calM_{E,\calL,\calL^+,m}(\Fpbar)\)
be a closed point in \(Y\subset \calM_{E,\calL,\calL^+,m}\).
Let \(*_y\) be the Rosati involution attached to \(\lambda\)
on the semisimple algebra
\(\End^0_{\ringO_E}(A_y)=\End_{\ringO_E}(A_y)\otimes_{\ringO_E}E\).
Denote by \({\rm H}_y\) the unitary group over \(\bbQ\) attached to
\((\End^0_{\ringO_E}(A_y),*_y)\), so 
\[
{\rm H}_y(R)=\left\{u\in \left(\End^0_{\ringO_E}(A_y)\otimes_{\bbQ}R
\right)^{\times}\,\mid\,u\cdot *_y(u)=*_y(u)\cdot u={\rm Id}_{A_y}
\right\}
\]
for every \(\bbQ\)-algebra \(R\).
Let \({\rm H}_y(\Zp)\) be the subgroup of \({\rm H}_y(\Qp)\) consisting
of all elements of \({\rm H}_y(\Qp)\) which induces an automorphism
of \((A_y[p^{\infty}],\iota_y[p^{\infty}],\lambda_y[p^{\infty}])\).
Denote by \({\rm H}_y(\bbZ_{(p)})\) the intersection of 
\({\rm H}_y(\bbQ)\) and \({\rm H}_y(\Zp)\) inside \({\rm H}_y(\Qp)\),
i.e.\ it consists of all elements \(u\in {\rm H}_y(\bbQ)\) such
that \(u\) induces an automorphism of 
\((A_y,\iota_y,\lambda_y)[p^{\infty}]\).
\smallbreak

The compact \(p\)-adic group
\({\rm Aut}((A_y,\iota_y,\lambda_y)[p^{\infty}])\)
operates naturally on the deformation space
\({\rm{Def}}\left((A_y,\iota_y,
\lambda_y)[p^{\infty}]/k
\right)\).
So we have a natural action of the compact \(p\)-adic group
\({\rm Aut}((A_y,\iota_y,\lambda_y)[p^{\infty}])\)
on the formal scheme \(\calM_{E,\calL,\calL^+,m}^{/y}\) via
the canonical isomorphism
\[
\calM_{E,\calL,\calL^+,m}^{/y}
={\rm{Def}}\left((A_y,\iota_y,\lambda_y)/k\right)
\xrightarrow[\sim]{\mbox{\scriptsize{\rm Serre-Tate}}}
{\rm{Def}}\left((A_y,\iota_y,
\lambda_y)[p^{\infty}]/k
\right)\,.
\]
\smallbreak

\begin{thm}\label{lsphilb}
Notation as above.  Assume that the closed subscheme 
\(Y\subset \calM_{E,\calL,\calL^+,m}\) over \(\Fpbar\) is stable under all 
\(\SL_2(E\otimes_{\bbQ}\bbA_f^{(p)})\)-Hecke correspondences
on the Hilbert modular variety \(\calM_{E,\calL,\calL^+,m}\).
Then the closed formal subscheme \(Y^{/y}\) of 
\(\calM_{E,\calL,\calL^+,m}^{/y}\) is stable under the
action by elements of the subgroup
\({\rm H}_y(\Zp)\) of 
\({\rm Aut}(A_y[p^{\infty}],\iota_y[p^{\infty}],\lambda_y[p^{\infty}])\).
\end{thm}

\proof The proof of Thm.\ \ref{lsphilb} is similar to 
that of Thm.\ \ref{lsp}, and is already contained in the proof of
Cor.\ \ref{cor:HMVlocalorbit}.
\hfill\qed
\medbreak

\begin{thm}\Bl \label{EOconseq}
Let \(n\geq 3\) be an integer relatively prime to \(p\).
Let \(\ell\) be a prime number, \(\ell\neq p\).
\begin{itemize}
\rmitem[(i)] Every closed subset of \(\calA_{g,n}\) over
\(\Fpbar\), which is stable 
under all Hecke correspondences on \(\calA_{g,n}\) coming from
\(\Sp_{2g}(\Ql)\),
contains a supersingular point.

\rmitem [(ii)] Similarly, every closed subset in a Hilbert modular
variety \(\calM_{E,\calL,\calL^+,n}\) over \(\Fpbar\),
which is stable under all 
\(\SL_2(E\otimes\Ql)\)-Hecke correspondences 
on \(\calM_{E,\calL,\calL^+,n}\), 
contains a supersingular point.
\end{itemize}
\end{thm}
\medbreak

\noindent{\bf Remark.}\enspace Thm.\ \ref{EOconseq} follows from 
the main theorem of \cite{FO-EO} and Prop.\ \ref{infho} below. 
See also \ref{Appl}.
\bigbreak

\begin{prop}\Bl \label{infho}
Let \(k\supset \Fp\) be an algebraically closed field.
Let \(\ell\) be a prime number, \(\ell\neq p\).
Let \(n\geq 3\) be an integer prime to \(p\).
\begin{itemize}
\rmitem[(i)] Let \(x=[(A_x,\lambda_x,\eta_x)]\in\calA_{g,1,n}(k)\)
be a closed point of \(\calA_{g,1,n}\).
If \(A_x\) is supersingular, then the prime-to-\(p\) Hecke
orbit \(\calH_{\Sp_{2g}}^{(p)}(x)\) is finite.
Conversely, if \(A_x\) is not supersingular, then
the \(\ell\)-adic Hecke orbit \(\calH_{\ell}^{\Sp_{2g}}(x)\)
is infinite for every prime number \(\ell\neq p\).

\rmitem[(ii)] Let \(y=[(A_y,\iota_y,\lambda_y,\eta_y)]
\in\calM_{E,\calL,\calL^+,n}(k)\) be a closed point of 
a Hilbert modular variety \(\calM_{E,\calL,\calL^+,n}\).
If \(A_y\) is  is supersingular, then the prime-to-\(p\) Hecke
orbit \(\calH_{\SL_{2,E}}^{(p)}(y)\) is finite.
Conversely,  if \(A_y\) is not supersingular, then
the \(v\)-adic Hecke orbit \(\calH_{v}^{\SL_{2,E}}(y)\)
is infinite for every prime ideal \(\wp_v\) of \(\ringO_E\)
which does not contain \(p\).
\end{itemize}
\end{prop}
\medbreak

\noindent{\bf Remark.} (1) The statement (i) is proved in Prop.\ 1,
p.\ 448 of \cite{Chai-Density}, see \ref{ssfinite}.  The proof of (ii) is similar.
The key to the proof of the second part of (i) is a bijection
\[
\calH_{\ell}^{\Sp_{2g}}(x) 
\xleftarrow{\sim} \left(
{\rm H}_x(\bbQ)\cap \prod_{\ell'\neq \ell}K_{\ell}\right)
\backslash {\rm Sp}_{2g}(\Ql)/K_{\ell}\,
\]
where \(\ell'\) runs through all prime numbers not equal to \(\ell\)
or \(p\), \({\rm H_x}\) is the unitary group attached to
\((\End^0(A_x),*_x)\) as in Thm.\ \ref{lsp}.
The compact groups \(K_{\ell'}\) and \(K_{\ell}\) are
defined as follows: for every prime number \(\ell'\neq p\),
\(K_{\ell'}=\Sp_{2g}(\bbZ_{\ell'}\) if \((\ell', n)=1\),
and \(K_{\ell'}\) consists of all elements \(u\in \Sp_{2g}(\bbZ_{\ell'}\)
such that \(u\equiv 1\pmod{n}\) if \(\ell'| n\).
We have an injection \({\rm H}_x(\bbA_f^{(p)})\to \Sp_{2g}(\bbA_f^{(p)}\)
as in Thm.\ \ref{lsp}, so that the
intersection \({\rm H}_x(\bbQ)\cap \prod_{\ell'\neq \ell}K_{\ell}\)
makes sense. 
The second part of (i) follows from the group-theoretic fact
that a double coset as above is finite if and only if
\({\rm H}_x\) is a form of \(\Sp_{2g}\).
\smallbreak

(2) When the abelian variety \(A_x\) in (i) (resp. \(A_y\) in (ii))
is ordinary, one can also use the canonical lifting to \(W(k)\) to 
show that  \(\calH_{\ell}^{\Sp_{2g}}(x)\) 
(resp.\ \(\calH_{v}^{\SL_{2,E}}(y)\)) is infinite.
\bigbreak

The following irreducibility statement is handy for the proof
of Thm.\ \ref{dens_hilb}, in that it shortens the argument
and simplifies the logical structure of the proof.
\smallbreak

\begin{thm} \Bl \label{thm:irredH}
Let \(W\) be a locally closed subscheme of \(\calM_{F,n}\) 
which is smooth over \(\Fpbar\) and stable under all 
\(\SL_2(F\otimes\bbA_f^{(p)})\)-Hecke correspondences.
Assume that the \(\SL_2(F\otimes\bbA_f^{(p)})\)-Hecke correspondences
operates transitively on the set \(\Pi_0(W)\) of irreducible
components of \(W\), and some (hence all) maximal point of \(W\)
corresponds to a non-supersingular abelian variety. 
Then \(W\) is irreducible.
\end{thm}

\noindent{\bf Remark.} The argument in \cite{Chai-Mon} works
in the situation of \ref{thm:irredH}.  The following observations
may be helpful.\\
(i) The linear algebraic group \({\rm Res}_{F/\bbQ}(\SL_2)\) over 
\(\bbQ\) is a semisimple, connected and simply connected.
Therefore every subgroup of finite index in 
\(\SL_2(F\otimes\Ql)\) is equal to \(\SL_2(F\otimes\Ql)\),
for every prime number \(\ell\).
Consequently \(\SL_2(F\otimes \bbA_f^{(p)})\) has no proper subgroup
of finite index.\\
(ii) The only part of the argument in \cite{Chai-Mon} that needs
to be supplemented is the end of (4.1), where the fact that
\(\Sp_{2g}\) is simple over \(\Ql\) is used.
Let \(G_{\ell}\) be the image group of the \(\ell\)-adic monodromy
\(\rho_Z\) attached to \(Z\).
By definition, \(G_{\ell}\) is a closed subgroup of
\(\SL_2(F\otimes \Ql)=\prod_{v|\ell} \SL_2(F_v)\),
where \(v\) runs through all places of \(F\) above \(\ell\).
In the present situation of a Hilbert modular variety \(\calM_F\),
we need to know the fact that the projection of 
\(G_{\ell}\) to the factor \(\SL_2(F_v)\) is non-trivial
for every place \(v\) of \(F\) above \(\ell\) and
for every \(\ell\neq p\).
\bigbreak

\begin{thm}[Hilbert trick] \label{thm_hilb}
Given \(x_0\in\calA_{g,1,n}(\Fpbar)\), 
then there exist 
\begin{itemize}
\rmitem[(a)] totally real number fields \(F_1,\ldots,F_r\)
such that \(\sum_{i=1}^r [E_i:\bbQ]=g\),

\rmitem[(b)] an invertible \(\ringO_E\)-module \(\calL\)
with a notion of positivity \(\calL^+\), i.e.\
\(\calL^+\) is a union of connected components 
of \(\calL\otimes_{\bbQ}\bbR\) such that
such that \(\calL\otimes\bbR\) is the 
disjoint union of \(\calL^{+}\)
with \(-\calL^{+}\),

\rmitem[(c)] a positive integer \(a\) and 
a positive integer \(m\) such that 
\((m,p)=1\) and \(m\equiv 0\pmod{n}\),

\rmitem[(d)] a finite flat morphism
\(
g:\calM_{E,\calL,\calL^+m;a}^{\rm ord}\to 
\calM_{E,\calL, \calL^{+}, m}^{\rm ord}\,,
\)

\rmitem[(e)] a finite morphism
\(
f:\calM_{E,\calL,\calL^+,m;a}^{\rm ord}\to\calA_{g,n}^{\rm ord}\,,
\)

\rmitem[(f)] a point \(y_0\in\calM_{E,\calL,\calL^+,m;a}^{\rm ord}(\Fpbar)\) 
\end{itemize}
such that the following properties are satisfied.
\begin{itemize}
\rmitem[(i)] There is a projective system 
\(\widetilde\calM_{E,\calL,\calL^+;a}^{\rm ord}\)
of finite \'etale coverings of 
\(\calM_{E,\calL,\calL^+,m;a}\) on which the group 
\(\SL_2(E\otimes\bbA_{f}^{(p)})\) operates.
This \(\SL_2(E\otimes\bbA_{f}^{(p)})\)-action induces 
Hecke correspondences on \(\calM_{E,\calL,\calL^+,m;a}^{\rm ord}\)

\rmitem[(ii)] The morphism \(g\) is equivariant w.r.t.\ Hecke correspondences
coming from the group \(\SL_2(E\otimes\bbA_{f}^{(p)})\).
In other words, there is a \(\SL_2(E\otimes\bbA_f^{(p)})\)-equivariant 
morphism \(\tilde{g}\) from
the projective system \(\widetilde\calM_{E,\calL,\calL^+;a}^{\rm ord}\)
to the projective system 
\(\left(\calM_{E,\calL,\calL^+,md}^{\rm ord}\right)_{d\in \bbN-p\bbN}\)
which lifts \(g\).

\rmitem[(iii)] There exists an injective homomorphism
\(j_{E}:\SL_2(E\otimes_{\bbQ}\bbA_f^{(p)})\to \Sp_{2g}(\bbA_f^{(p)})\)
such that the finite morphism \(f\) is Hecke equivariant 
w.r.t.\ \(j_E\).

\rmitem[(iv)] We have \(f(y_0)=x_0\). 

\rmitem[(v)]  For every geometric point \(z\in\calM_{E,m;a}^{\rm ord}\),
the abelian variety underlying the fiber over \(g(z)\in\calM_{E,m}^{\rm ord}\)
of the universal abelian scheme over \(\calM_{E,m}^{\rm ord}\)
is isogenous to the abelian variety underlying the fiber
over \(f(z)\in\calA_{g,n}^{\rm ord}(\Fpbar)\) of
the universal abelian scheme over \(\calA_{g,n}^{\rm ord}(\Fpbar)\).
\end{itemize}
\end{thm}
\medbreak

\noindent{\bf Lemma.}\enspace 
Let \(A\) be an ordinary abelian variety over \(\Fpbar\) which is
simple.  Then 
\begin{itemize}
\rmitem[(i)] \(K:=\End^0(A)\) is a totally imaginary quadratic
extension of a totally real number field \(F\).

\rmitem[(ii)] \([F:\bbQ]=\dim(A)\).

\rmitem[(iii)] \(F\) is the fixed by the Rosati involution attached
to any polarization of \(A\).

\rmitem[(iii)]  Every place \(\wp\) of \(F\) above \(p\) splits in \(K\).
\end{itemize}

\proof The statements (i)--(iv) are immediate consequences of Tate's theorem
for abelian varieties over finite fields; see \cite{Tate-Bourb}.
\hfill\qed
\smallbreak

\noindent{\bf Lemma.}\enspace
Let \(K\) be a CM field, let \(E:={\rm M}_d(K)\), and let
\(*\) be a positive involution on \(E\) which induces the 
complex conjugation on \(K\).
Then there exists a CM field \(L\) which contains \(K\)
and a \(K\)-linear ring homomorphism \(h:L\to E\)
such that \([L:K]=d\) and \(h(L)\) is stable under the
involution \(*\).

\proof 
This is an exercise in algebra.  A proof using Hilbert irreducibility
can be found on p.\ 458 of \cite{Chai-Density}.
\qed
\bigbreak

\noindent{{\bf\scshape Proof of Thm.\ \ref{thm_hilb} (Hilbert trick).}}
\enspace
\smallbreak

\noindent{\scshape Step 1.}\enspace
Consider the abelian variety \(A_{0}\) attached to the given
point \(x_0=[(A_0,\lambda_0,\eta_0)]\in\calA_{g,1,n}^{\rm ord}(\Fpbar)\). 
By the two lemmas above
 there exist totally real number fields
\(F_1,\ldots, F_r\) and an embedding
\(\iota_{0}:E:=F_1\times\cdots\times F_r\hookrightarrow \End^0(A_{0})\)
such that 
\(E\) is fixed under the Rosati involution on \(\End^0(A_0)\)
attached to the principal polarization \(\lambda_0\),
and 
\([E:\bbQ]=g=\dim(A_0)\).
\smallbreak

The intersection of \(E\) with \(\End(A_0)\) is an order \(\ringO_1\)
of \(E\), so we can regard \(A_0\) as an abelian variety with
action by \(\ringO_1\). We claim that there exists an 
\(\ringO_E\)-linear abelian variety \(B\) and an
\(\ringO_1\)-linear isogeny \(\alpha:B \to A_0\).
This claim follows from a standard ``saturation construction'' as 
follows. 
Let \(d\) be the  order of the finite
abelian group \(\ringO_E/\ringO_1\). 
Since \(A_0\) is ordinary, one sees by Tate's theorem 
(the case when \(K\) is a finite field in Thm.\ \ref{lmon}) 
that \((d,p)=1\).
For every prime divisor \(\ell\neq p\) of \(d\),
consider the \(\ell\)-adic Tate module \({\rm T}_{\ell}(A_0)\) as a lattice
inside the free rank two \(E\)-module \({\rm V}_{\ell}(A_0)\).
Then the lattice \(\Lambda_{\ell}\) 
generated by \(\ringO_E\cdot {\rm T}_{\ell}(A_0)\)
is stable under the action of \(\ringO_E\) by construction.
The finite set of lattices \(\,\{\Lambda_{\ell}: \ell|d\}\,\)
defines an \(\ringO_E\)-linear abelian variety \(B\)
and an \(\ringO_1\)-linear isogeny \(\beta_0:A_0\to B\)
which is killed by a power \(d^i\) of \(d\).
Let \(\alpha:B\to A_0\) be the isogeny such that 
\(\alpha\circ\beta_0=[d^i]_{A_0}\).  The claim is proved.
\medbreak

\noindent{\scshape Step 2.}\enspace
The construction in Step 1 gives us a triple \((B, \alpha, \iota_{x_0})\),
where \(B\) is
an abelian variety \(B\) over \(\Fpbar\,\),
\(\alpha:B \to A_x\) is an isogeny over \(\Fpbar\), and 
\(\iota_B:\ringO_E\to\End(B)\) is an an injective ring homomorphism
such that \(\alpha^{-1}\circ \iota_x(u)\circ \alpha =\iota_B(u)\)
for every \(u\in \ringO_E\).
Let \(\calL_B:={\rm Hom}_{\ringO_E}^{\rm sym}(B,B^t)\) be set of 
all \(\ringO_E\)-linear symmetric homomorphisms from 
\(B\) to the dual \(B^t\) of \(B\).
The set \(\calL_B\) has a natural structure as an \(\ringO_E\)-module.
By Tate's theorem (the case when \(K\) is a finite field in
Thm.\ \ref{lmon}, see \ref{struct2}) one sees that 
\(\calL_B\) is an invertible \(\ringO_E\)-module, and
the natural map
\[
\lambda_B:B\otimes_{\ringO_E}\calL_B\to B^t
\]
is an \(\ringO_E\)-linear isomorphism.
The subset of elements in \(\calL\) which are polarizations
defines a notion of positivity \(\calL^+\) on \(\calL\)
such that \(\calL_B\cap \calL_B^{+}\) is the subset of
\(\ringO_E\)-linear polarizations on \((B,\iota_B)\).
\medbreak

\noindent{\scshape Step 3.}\enspace
Recall that the Hilbert modular variety \(\calM_{E,\calL,\calL_+,n}\)
classifies (the isomorphism class of) all quadruples 
\((A\to S,\iota_A,\lambda_A,\eta_A)\), where
\((A\to S,\iota_A)\) is an 
\(\ringO_E\)-linear abelian schemes,
\(\lambda_A: \calL\to {\rm Hom}^{\rm sym}_{\ringO_E}(A,A^t)\)
is an injective \(\ringO_E\)-linear map such that
the resulting morphism
\(\calL\otimes A\xrightarrow{\sim} A^t\)
is an isomorphism of abelian schemes and
every element of \(\calL\cap \calL^+\) gives rise to 
an \(\ringO_E\)-linear polarization,
and \(\eta_A\) is
an \(\ringO_E\)-linear level structure on \((A,\iota_A)\).
In the preceding paragraph, if we choose an \(\ringO_E\)-linear
level-\(n\) structure \(\eta_B\) on \((B,\iota_B)\),
then \(y_1:=[(B,\iota_B,\lambda_B,\eta_B)]\) is an 
\(\Fpbar\)-point of the Hilbert modular variety
\(\calM_{E,\calL_B,\calL_B^+,n}\).
The element \(\alpha^*(\lambda_0)\) is 
an \(\ringO_E\)-linear polarization on 
\(B\), hence it is equal to \(\lambda_B(\mu_0)\)
for a unique element \(\mu_0\in\calL\cap \calL^+\).
\smallbreak

Choose a positive integer \(m_1\) with \({\rm gcd}(m_1,p)=1\)\
and \(a\in \bbN\) such that
\({\rm Ker}(\alpha)\) is killed by \(m_1 p^a\).
Let \(m=m_1 n\).
Let \((A,\iota_A,\lambda_A,\eta_A)\to 
\calM_{E,\calL,\calL^+,m}^{\rm ord}\) be the universal polarized
\(\ringO_E\)-linear abelian scheme over the ordinary
locus \(\calM_{E,\calL,\calL^+,m}^{\rm ord}\) of 
\(\calM_{E,\calL,\calL^+,m}\).
Define a scheme \(\calM_{E,\calL,\calL^+,m;a}^{\rm ord}\)
over \(\calM_{E,\calL,\calL^+,m}^{\rm ord}\) by
\[
\calM_{E,\calL,\calL^+,m;a}^{\rm ord}:=
\underline{\rm Isom}_{\calM_{E,\calL,\calL^+,m}^{\rm ord}}^{\ringO_E}
\left((B,\iota_B,\lambda_B)[p^a]\times_{\Spec(\Fpbar)}
\calM_{E,\calL,\calL^+,m}^{\rm ord},
(A,\iota_A,\lambda_A)[p^a]
\right)\,.
\]
In other words \(\calM_{E,\calL,\calL^+,m;a}^{\rm ord}\) 
is the moduli
space of \(\ringO_E\)-linear ordinary abelian varieties with
level-\(m p^a\) structure,
where we have used the \(\ringO_E\)-linear polarized
truncated Barsotti-Tate group 
\((B,\iota_B,\lambda_B)[p^m]\) as the ``model'' for
the \(mp^a\)-torsion subgroup scheme of the universal abelian
scheme over \(\calM_{E,\calL,\calL^+,m}^{\rm ord}\).
Let 
\[
g:\calM_{E,\calL,\calL^+,m;a}^{\rm ord}
\to \calM_{E,\calL,\calL^+,m}^{\rm ord}
\]
be the structural morphism of \(\calM_{E,\calL,\calL^+,m;a}^{\rm ord}\),
the source of \(g\) being an fppf sheaf of sets
on the target of \(g\).
Notice that the structural morphism
\(g:\calM_{E,\calL,\calL^+,m;a}^{\rm ord}
\to \calM_{E,\calL,\calL^+,m}^{\rm ord}\) has a natural
structure as a torsor over the constant finite flat group scheme
\[\underline{\rm Aut}\left((B,\iota_B,\lambda_B)[p^a]\right)
\times_{\Spec(\Fpbar)}\calM_{E,\calL,\calL^+,m}^{\rm ord}\,.
\]
We have constructed the finite flat morphism \(g\) as promised in
Thm.\ \ref{thm_hilb} (d). We record some properties of this morphism.
\medbreak

\n
The group
\(\underline{\rm Aut}\left(B,\iota_B,\lambda_B)[p^a]\right)\)
sits in the middle of a short exact sequence
\[
0\to \underline{\rm Hom}_{\ringO_E}\left(B[p^a]_{\rm et},B[p^a]_{\rm mult}
\right)
\to \underline{\rm Aut}\left((B,\iota_B,\lambda_B)[p^a]\right)
\to  \underline{\rm Aut}\ringO_E\left(B[p^a]_{\rm et}\right)
\to 0\,.
\]
The morphism  \(g:\calM_{E,\calL,\calL^+,m;a}^{\rm ord}
\to \calM_{E,\calL,\calL^+,m}^{\rm ord}\)
factors as
\[
\calM_{E,\calL,\calL^+,m;a}^{\rm ord}
\xrightarrow{g_1}
{\calM_{E,\calL,\calL^+,m;a}^{\rm ord,et}}
\xrightarrow{g_2}
\calM_{E,\calL,\calL^+,m}^{\rm ord}\,,
\]
where \(g_1\) 
is defined as the
push-forward by the surjection 
\[\underline{\rm Aut}\left(B,\iota_B,\lambda_B)[p^a]\right)
\twoheadrightarrow
\underline{\rm Aut}\ringO_E\left(B[p^a]_{\rm et}\right)
\]
of the \(\underline{\rm Aut}\left(B,\iota_B,\lambda_B)[p^a]\right)\)-torsor
\(\calM_{E,\calL,\calL^+,m;a}^{\rm ord}
\).
Notice that the morphism \(g_1\) is finite flat and purely inseparable,
and \(\calM_{E,\calL,\calL^+,m;a}^{\rm ord,et}\) is integral.
Moreover \({\calM_{E,\calL,\calL^+,m;a}^{\rm ord}}\)
and \(\calM_{E,\calL,\calL^+,m;a}^{\rm ord}\) are irreducible
by \cite{Ribet}, \cite{Deligne.Ribet}, \cite{Rap} and \cite{D.P}.
\medbreak

\noindent{\scshape Step 4.}\enspace
Let \(\pi_{n,m}:\calM_{E,\calL,\calL^+,m}\to \calM_{E,\calL,\calL^+,n}\)
be the natural projection.
Denote by \[A[mp^a]\to \calM_{E,\calL,\calL^+,m}^{\rm ord}\] the
kernel of \([mp^a]\) on
\(A\to \calM_{E,\calL,\calL^+,m}^{\rm ord}\),
and let \(g^*A[mp^a]\to \calM_{E,\calL,\calL^+,m;a}^{\rm ord}\)
be the pull-back of \(A[mp^a]\to \calM_{E,\calL,\calL^+,m}^{\rm ord}\) 
by \(g\).
By construction the \(\ringO_E\)-linear finite flat group scheme
\(\,g^*A[mp^a]\to
\calM_{E,\calL,\calL^+,m;a}^{\rm ord}\) is constant
via a tautological trivialization 
\[
\psi:
\underline{\rm Aut}\left(B,\iota_B,\lambda_B)[p^a]\right)
\times_{\Spec(\Fpbar)} \calM_{E,\calL,\calL^+,m}^{\rm ord}
\xrightarrow{\sim}
\calM_{E,\calL,\calL^+,m;a}^{\rm ord}
\]
Choose a point \(y_0\in\calM_{E,\calL,\calL^+,m;a}^{\rm ord}(\Fpbar)\)
such that \((\pi_{n,m}\circ g)(y_0)=y_1\).
The fiber over \(y_0\) of \(g^*A[mp^a]\to 
\calM_{E,\calL,\calL^+,m;a}^{\rm ord}\)
is naturally identified with \(B[mp^a]\).
Let \(K_0:={\rm Ker}(\alpha:B\to A_0)\),
and let 
\[K:=\psi\left(K_0\times_{\Spec(\Fpbar)}
\calM_{E,\calL,\calL^+,m;a}^{\rm ord}
\to \calM_{E,\calL,\calL^+,m;a}^{\rm ord}\right)\,,
\]
the subgroup scheme of \(g^*A[mp^a]\) which corresponds to
the constant group \(K_0\) under the trivialization \(\psi\).
The element \(\mu_0\in\calL\cap\calL^+\) defines a polarization
on the abelian scheme \(g^*A\to \calM_{E,\calL,\calL^+,m;a}^{\rm ord}\), the 
pull-back by \(g\) of the universal polarized \(\ringO_E\)-linear
abelian scheme over
\(A\to \calM_{E,\calL,\calL^+,m}^{\rm ord}\).
The group \(K\) is a maximal totally isotropic subgroup scheme
of \(g^*{\rm Ker}(\lambda_A(\mu_0))\to \calM_{E,\calL,\calL^+,m;a}^{\rm ord}
\), because 
\(g^*{\rm Ker}(\lambda_A(\mu_0))\) is constant
and \(K_0\) is a maximal totally isotropic subgroup scheme
of \({\rm Ker}(\lambda_B(\mu_0))\).
\medbreak

Consider the quotient abelian scheme 
\(A'\to \calM_{E,\calL,\calL^+,m;a}^{\rm ord}\) of
\(g^*A\to \calM_{E,\calL,\calL^+,m;a}^{\rm ord}\)
by \(K\). Recall that we have defined an element
\(\mu_0\in\calL\cap \calL^+\) in Step 3.
The polarization \(g^*(\lambda_A(\mu_0))\) on the abelian scheme
\(g^* A\to \calM_{E,\calL,\calL^+,m;a}^{\rm ord}\) descends to
the quotient abelian scheme
\(A'\to \calM_{E,\calL,\calL^+,m;a}^{\rm ord}\),
giving it a principal polarization \(\lambda_{A'}\).
Moreover the \(n\)-torsion subgroup scheme
\(A'[n]\to \calM_{E,\calL,\calL^+,m;a}^{\rm ord}\)
is constant, as can be checked easily.
Choose a level-\(n\) structure \(\eta_{A'}\) for \(A'\).
The triple \((A',\lambda_{A'},\eta_{A'})\) over 
\(\calM_{E,\calL,\calL^+,m;a}^{\rm ord}\) defines a 
morphism \(f:\calM_{E,\calL,\calL^+,m;a}^{\rm ord}
\to \calA_{g,1,n}^{\rm ord}\)
by the modular definition of \(\calA_{g,1,n}^{\rm ord}\),
since every fiber of \(A'\to \calM_{E,\calL,\calL^+,m;a}^{\rm ord}\)
is ordinary by construction.
We have constructed the morphism \(f\) as required in
\ref{thm_hilb} (e), and also the point \(y_0\) as required
in \ref{thm_hilb} (f).
\medbreak

\noindent{\scshape Step 5.}\enspace
So far we have constructed the morphisms
\(g\) and \(f\) as required in Thm.\ \ref{thm_hilb}.
To construct the homomorphism \(j_E\) as required in (iii),
one uses the first homology group
\(V:={\rm H}_1(B,\bbA_f^{(p)})\), and the symplectic pairing
\(\langle\cdot,\cdot\rangle\) induced by the polarization
\(\alpha^*(\lambda_0)=\lambda_B(\mu_0)\) constructed
in Step 3.
Notice that \(V\) has a natural structure
as a free \(E\otimes_{\bbQ}\bbA_f^{(p)}\)-module of rank two. 
Also, \(V\) is a free \(\bbA_f^{(p)}\)-module of rank \(2g\).
So we get an embedding
\(j_E:\SL_{E\otimes\bbA_f^{(p)}}(V)\hookrightarrow
\Sp_{\bbA_f^{(p)}}(V,\langle\cdot,\cdot\rangle)\).
We have finished the construction of \(j_E\).
\smallbreak 

We define 
\(\widetilde{\calM}_{E,\calL,\calL^+,a}^{\rm ord}\)
to be the projective system 
\(\,\varprojlim_{md} \calM_{E,\calL,\calL^+,md;a}\),
where \(d\) runs through all positive integers which are prime to \(p\).
This finishes the last construction needed for Thm.\ \ref{thm_hilb}.
\medbreak

By construction we have \(f(y_0)=x_0\),
which is statement (iv). The statement (v) is clear by construction.
The statements (i)--(iii) can be verified without difficulty from
the construction. 
\hfill\qed

\noindent{{\bf\scshape Proof of Theorem \ref{dens_hilb}.} (Density of
ordinary Hecke orbits in \(\calM_{E,\calL,\calL^+,n}\))}
\smallbreak

\noindent{\scshape Reduction step.}\enspace
\smallbreak
From the product decomposition
\[
\calM_{E,\calL,\calL^+,n}=\calM_{F_1,\calL_1,\calL_1^+,n}
\times_{\Spec(\Fpbar)}\cdots\times_{\Spec(\Fpbar)}
\calM_{F_r,\calL_r,\calL_r^+,n}
\]
of the Hilbert modular variety \(\calM_{E,\calL,\calL^+,n}\),
we see that it suffices to prove Thm.\ \ref{dens_hilb} when
\(r=1\), i.e.\ \(E=F_1=:F\) is a totally real number field.
Assume this is the case from now on.
\medbreak

The rest of the proof is divided into four steps.
\smallbreak

\noindent{{\bf\scshape Step 1} (Serre-Tate coordinates for Hilbert modular
varieties)}.\enspace
\smallbreak

\noindent{\scshape Claim.}\enspace
The Serre-Tate local coordinates at a closed ordinary point of
\(z\in\calM_{F,\calL,\calL^+,n}^{\rm ord}\) of 
a Hilbert modular variety \(\calM_{F,\calL,\calL^+,n}\)
admits a canonical decomposition
\[
\calM_{F,\calL,\calL^+,n}^{/z}\cong 
\prod_{\wp\in\Sigma_{F,p}} \calM^{z}_{\wp}\,,
\qquad
\calM^{z}_{\wp}=\underline{\Hom}_{\ringO_{F,\wp}}\left(
{\rm T}_p(A_z[\wp^{\infty}]_{\rm et}),\ e_{\wp}\cdot\widehat{A_z}
\right)\,,
\]
where 
\begin{itemize}
\item the indexing set \(\Sigma_{F,p}\) is the finite set consisting of 
of all prime ideals of \(\ringO_F\) above \(p\),

\item the \((\ringO_F\otimes\Zp)\)-linear formal torus \(\widehat{A_z}\) is 
the formal completion of the ordinary abelian variety \(A_z\), 

\item \(e_{\wp}\) is the irreducible idempotent in 
\(\ringO_F\otimes_{\bbZ}\Zp\)
so that \(e_{\wp}\cdot (\ringO_F\otimes_{\bbZ}\Zp)\) is
equal to the factor \(\ringO_{F_{\wp}}\) of \(\ringO_F\otimes_{\bbZ}\Zp\).
\end{itemize}
Notice that \(e_{\wp}\widehat{A_z}\) is the formal torus
attached to the multiplicative Barsotti-Tate group
\(A_z[\wp^{\infty}]_{\rm mult}\) over \(\Fpbar\).
\medbreak

\noindent{Proof of Claim.}\enspace 
The decomposition 
\(\ringO_F\otimes_{\bbZ}\Zp=\prod_{\wp\in \Sigma_{F,p}}\ringO_{F_{\wp}}\)
induces a decomposition of the formal scheme
\(\,\calM_{F,\calL,\calL^+,n}^{/z}\,\) into a product
\(\,\calM_{F,\calL,\calL^+,n}^{/z}=\prod_{\wp\in\Sigma_{F,p}} \calM_{\wp}^z\,\)
for every closed point \(z\) of \(\calM_{F,\calL,\calL^+,n}\):
Let \((A/R, \iota)\) be an \(\ringO_F\)-linear abelian
scheme over an Artinian local ring \(R\).
Then we have a decomposition 
\(A[p^{\infty}]=\prod_{\wp\in\Sigma_{F,p}}A[\wp^{\infty}]\)
of the Barsotti-Tate group attached to \(A\),
and each \(A[\wp^{\infty}]\) is a deformation of 
\(A\times_{\Spec(R)}\Spec(R/\grm)\) over \(\Spec(R)\).
\smallbreak

If \(z\) corresponds to an ordinary abelian variety \(A_z\),
then 
\(\calM_{\wp}^z\) is canonically isomorphic to the 
\(\ringO_{F,\wp}\)-linear formal torus 
\(\underline{\Hom}_{\ringO_{F,\wp}}(A_z[\wp^{\infty}]_{\rm et}, 
e_{\wp}\cdot\widehat{A_z})\),
which is the factor ``cut out'' in the 
\((\ringO_F\otimes_{\bbZ}\Zp)\)-linear formal torus
\[
\calM_{F,\calL,\calL^+,n}^{/z}=\underline{\rm Hom}_{\ringO_F\otimes\Zp}
\left({\rm T}_p(A_z[p^{\infty}],\
\widehat{A_z}
\right)
\]
by the idempotent \(e_{\wp}\) in \(\ringO_F\otimes\Zp\).
Each factor \(\calM^{z}_{\wp}\) in the above decomposition
is a formal torus of dimension \([F_{\wp}:\Qp]\), with a
natural action by \(\ringO_{F,\wp}^{\times}\);
it is non-canonically isomorphic to the 
\(\ringO_{\wp}\)-linear formal torus
\(\widehat{A_z}\).
\hfill\qed
\medbreak

\noindent{\bf\scshape Step 2.}\enspace (Linearization)
\smallbreak

\noindent{\scshape Claim.}\enspace 
For every closed point \(z\in Z_F^{\rm ord}(\Fpbar)\) in the ordinary locus
of \(Z_F\), there exists a non-empty subset \(S_z\subset \Sigma_{F,p}\)
such that 
\(Z_F^{/z}=\prod_{\wp \in S_z} \calM^{z}_{\wp}\),
where \(\calM^{z}_{\wp}\) is the factor of the Serre-Tate
formal torus \(\calM_{F,\calL,\calL^+,n}^{/z}\)
corresponding to \(\wp\).
\small

\proof The \(\ringO_F\)-linear abelian variety \(A_z\) is an
ordinary abelian variety defined over \(\Fpbar\).
Therefore \(\End^0_{\ringO_F}(A_z)\) is a totally imaginary
quadratic extension field \(K\) of \(F\) which is split
over every prime ideal \(\wp\) of \(\ringO_F\) above \(p\),
by Tate's theorem (the case when \(K\) is a finite field in
Thm.\ \ref{lmon}).
By the local stabilizer principle, 
\(Z_F^{/z}\) is stable under the norm-one subgroup \(U\)
of \((\ringO_K\otimes_{\bbZ}\Zp)^{\times}\).
Since every prime \(\wp\) of \(\ringO_F\) above \(p\) splits in
\(\ringO_K\), \(U\) is isomorphic to 
\(\prod_{\wp\in\Sigma_{F,p}}\ringO_{F,\wp}^{\times}\)
through its action on the \((\ringO_F\otimes\Zp)\)-linear 
formal torus \(\widehat{A_z}\).
The factor \(\ringO_{F,\wp}^{\times}\) of \(U\) operates
on the \(\ringO_{F,\wp}\)-linear formal torus
\(\calM_{\wp}^z\) through the character \(t\mapsto t^2\),
i.e.\ a typical element 
\(t\in U=\prod_{\wp\in\Sigma_{F,p}}\ringO_{F,\wp}^{\times}\) 
operates on the \((\ringO_F\otimes\Zp)\)-linear formal torus
\(\calM_{F,\calL,\calL^+,n}^{/z}\) through the element 
\(t^2\in U=(\ringO_F\otimes\Zp)^{\times}\).
The last assertion can be seen through the formula
\[
\calM_{F,\calL,\calL^+,n}^{/z}=\underline{\Hom}_{\ringO_{F}\otimes\Zp}\left(
{\rm T}_p(A_z[p^{\infty}]_{\rm et}),\ \widehat{A_z}
\right)
\,,
\]
because any element \(t\) of \(U\xrightarrow{\sim}\ringO_F\otimes_{\bbZ}\Zp\)
operates via \(t\) (resp.\ \(t^{-1}\)) on the \(\ringO_{F}\otimes\Zp\)-linear
formal torus \(e_{\wp}\widehat{A_z}\) 
(resp.\ the  \(\ringO_{F}\otimes\Zp\)-linear Barsotti-Tate group
\(A_z[p^{\infty}]_{\rm et}\)).
\smallbreak

The local rigidity
theorem \ref{localrig} implies that 
\(Z^{/z}_F\) is a formal subtorus of the Serre-Tate formal torus
\(\calM_F^{/z}\).  
For every \(\wp \in \Sigma_{F,p}\), 
let \(X_{\wp,*}\) be the cocharacter group of the 
\(\ringO_{F_{\wp}}\)-linear formal torus \(\calM_{\wp}^z\),
so that \(\prod_{\wp\in\Sigma_{F,p}} X_{\wp,*}\) is the cocharacter group
of the Serre-Tate formal torus \(\calM_{F,\calL,\calL^+,n}^{/z}\).
Let \(Y_*\) be the cocharacter group of the formal torus
\(Z^{/z}_F\).
We know that \(Y_*\) is a co-torsion free \(\Zp\)-submodule
of the rank-one free 
\(\,\left(\prod_{\wp\in\Sigma_{F,p}}\ringO_{F,\wp}\right)\)-module
\(\,\prod_{\wp\in\Sigma_{F,p}} X_{\wp,*}\,\),
and  \(Y_*\) is stable under  
multiplication by elements of the subgroup
\(\prod_{\wp\in\Sigma_{F,p}} (\ringO_{F,\wp}^{\times})^2\)
of \(\prod_{\wp\in\Sigma_{F,p}}\ringO_{F,\wp}^{\times}\).
It is easy to see that the additive subgroup generated
by \(\,\prod_{\wp\Sigma_{F,p}} (\ringO_{F,\wp}^{\times})^2\,\)
is equal to \(\,\prod_{\wp\in\Sigma_{F,p}}\ringO_{F,\wp}\,\),
i.e.\ \(Y_*\) is a 
\(\,\left(\prod_{\wp\in\Sigma_{F,p}}\ringO_{F,\wp}\right)\)-submodule
of \(\prod_{\wp} X_{\wp,*}\).
Hence there exists a subset 
\(S_z\subseteq \Sigma_{F,p}\) such that 
\(Y_*=\prod_{\wp\in S_z} X_{\wp,*}\).
Since the prime-to-\(p\) Hecke orbit 
\(\calH_{\SL_{2,F}}^{(p)}(x)\) is infinite by \ref{infho}, we have
\(0<\dim(Z_F)=\sum_{\wp\in S_z} [F_{\wp}:\Qp]\),
hence \(S_z\neq \emptyset\) for every ordinary point
\(z\in Z_F(x)(\Fpbar)\).
\hfill\qed
\medbreak

\noindent{{\bf\scshape Step 3} (globalization)}\enspace
\smallbreak

\noindent{\scshape Claim.} The finite set \(S_z\) is
independent of the point \(z\), i.e.\ there exists a 
subset \(S\subset\Sigma_{F,p}\) such that 
\(S_z=\Sigma_{F,p}\) for all \(z\in Z_F^{\rm ord}(\Fpbar)\).
\medbreak

\noindent{\scshape Proof of Claim.}\enspace
Consider the diagonal embedding
\(\Delta_Z:Z_F\to Z_F\times_{\Spec(\Fpbar)} Z_F\),
the diagonal embedding 
\(\Delta_{\calM}:\calM_{F,n}\to \calM_{F,n}\times_{\Spec(\Fpbar)}
\calM_{F,n}\),
and the map \(\Delta_{Z,\calM}\)
from \(\Delta_Z\) to \(\Delta_{\calM}\)
induced by the inclusion \(Z_F\hookrightarrow \Delta_{\calM}\).
Let \(\calP_Z\) be the formal completion of 
\(Z_F\times_{\Spec(\Fpbar)} Z_F\) along \(\Delta_{Z}(Z_F)\),
and let \(\calP_{\calM}\) be the formal completion of
\(\calM_{F,n}\times_{\Spec(\Fpbar)}\calM_{F,n}\)
along \(\Delta_{\calM}(\calM_{F,n})\).
The map \(\Delta_{Z,\calM}\) induces a closed embedding
\(i_{Z,\calM}: \calP_{Z}\hookrightarrow \calP_{\calM}\).
We regard \(\calP_Z\) (resp.\ \(\calP_{\calM}\))
as a formal scheme over \(Z_F\) (resp.\ \(\calM_{F,n}\))
via the first projection.
\smallbreak

The decomposition \(\ringO_F\otimes_{\bbZ}\Zp
=\prod_{\wp\in\Sigma_{F,p}}\ringO_{F,\wp}\) 
induces a fiber product decomposition
\[
\calP_{\calM}=\prod_{\wp\in\Sigma_{F,p}} \left(
\calP_{\wp}\to \calM_{F,n}\right)
\]
over the base scheme \(\calM_{F,n}\), where 
\(\calM_{\wp}\to \calM_{F,n}\) is a smooth formal scheme
of relative dimension \([F_{\wp}:\Qp]\)
with a natural section \(\delta_{\wp}\),
for every \(\wp\in\Sigma_{F,p}\),
and the formal completion of the fiber of \((\calM_{\wp},\delta_{\wp})\)
over any closed point \(z\) of the base scheme 
\(\calM_{F,n}\) is canonically isomorphic to the
formal torus \(\calM_{\wp}^z\).
In fact one can show that 
\(\calM_{\wp}\to \calM_{F,n}\) has a natural structure 
as a formal torus of relative dimension
\([F_{\wp}:\Qp]\), with \(\delta_{\wp}\) as the zero section;
we will not need this fact here.
Notice that \(\calP_Z\to Z_F\) is a closed formal subscheme
of \(\calP_{\calM}\times_{\calM_{F,n}}Z_F\to Z_F\).
The above consideration globalizes the ``pointwise'' construction
of formal completions at closed points.
\smallbreak

By Prop.\ \ref{thm:irredH}, \(Z_F\) is irreducible. 
We conclude from the
irreducibility of \(Z_F\) that there is a non-empty subset
\(S\subset \Sigma_{F,p}\) such that 
the restriction of \(\calP_Z\to Z_F\) to the ordinary locus
\(Z_F^{\rm ord}\) is equal to the fiber product over
\(Z_F^{\rm ord}\) of formal schemes
\(\calP_{\wp}\times_{\calM_{F,n}}Z_F^{\rm ord}\to Z_F^{\rm ord}\)
over \(Z_F^{\rm ord}\),
where \(\wp\) runs through the subset \(S\subseteq \Sigma_{F,p}\).
\hfill\qed
\medbreak

\noindent{\bf Remark.}\enspace (i) Without using Prop.\ \ref{thm:irredH},
the above argument shows that for each irreducible component 
\(Z_1\) of \(Z_F^{\rm ord}\), there exists a subset \(S\subset \Sigma_{F,p}\)
such that \(S_z=S\) for every closed point 
\(z\in Z_1(\Fpbar)\).
\smallbreak

(ii) There is an alternative proof of the claim 
that \(S_z\) is independent of \(z\): By Step 2, \(Z_F^{\rm ord}\)
is smooth over \(\Fpbar\).  Consider the relative cotangent sheaf
\(\Omega_{Z_F^{\rm ord}/\Fpbar}^1\), which is a locally free
\(\ringO_{Z_F^{\rm ord}}\)-module.
We recall that \(\Omega_{\calM_{F,\calL,\calL^+,n}^{\rm ord}/\Fpbar}^1\)
has a natural structure as a \(\ringO_F\otimes_{\bbZ}\Fp\)-module,
from the  Serre-Tate coordinates explained in Step 1.
By Step 2, we have
\[
\Omega_{Z_F^{\rm ord}/\Fpbar}^1\otimes_{\ringO_{Z_F,z}}\widehat{\ringO_{Z_F,z}}
= \sum_{\wp\in S_z} e_{\wp}\cdot 
\Omega_{\calM_{F,\calL,\calL^+,n}^{\rm ord}/\Fpbar}^1
\otimes_{\ringO_{\calM_{F,\calL,\calL^+,n}^{\rm ord}}}
\widehat{\ringO_{Z_F,z}}
\]
for every \(z\in Z_F^{\rm ord}(\Fpbar)\),
where \(\widehat{\ringO_{Z_F,z}}\) is the formal completion of the
local ring of \(Z_F\) at \(z\).
Therefore for each irreducible component \(Z_1\) of \(Z_F^{\rm ord}\)
there exists a subset \(S\subset \Sigma_{F,p}\) such
that 
\[
\Omega_{Z_1/\Fpbar}^1
=
\sum_{\wp\in S} e_{\wp}\cdot
\Omega_{\calM_{F,\calL,\calL^+,n}^{\rm ord}/\Fpbar}^1
\otimes_{\ringO_{\calM_{F,\calL,\calL^+,n}^{\rm ord}}}
\ringO_{Z_1/\Fpbar}\,.
\]
Hence \(S_z=S\) for every \(z\in Z_1(\Fpbar)\).
This argument was used in \cite{Chai-Density}; see Prop.\ 5 on p.\ 473
in \emph{loc.\ cit.}
\smallbreak

(iii) One advantage of the globalization argument in Step 3 is that
it makes the final Step 4 of the proof of Thm.\ \ref{dens_hilb} easier,
as compared with the two-page proof of Prop.\ 7 on p.\ 474 of
\cite{Chai-Density}.
\bigbreak

\noindent{\bf\scshape Step 4.}\enspace 
We have \(S=\Sigma_{F,p}\).
Therefore \(Z_F=\calM_{F,n}\).

\medbreak

\noindent{\scshape Proof of Step 4.}\enspace
\smallbreak

Notation as in Step 3 above. 
For every closed point \(s\) of \(Z_F\), the formal completion
\(Z_F^{/s}\) contains the product
\(\prod_{\wp\in S}\calM_{\wp}^{s}\).
By Thm.\ \ref{EOconseq}, \(Z_F\) contains a supersingular point 
\(s_0\).  
Consider the formal completion \(\hat{Z}:=Z_F^{/s_0}\),
which contains \(\hat{W}:=\prod_{\wp\in S}\calM_{\wp}^{s_0}\),
and the generic point \(\eta_{\hat{W}}\) of 
\(\Spec\left({\rm H}^0(\hat{W},\ringO_{\hat{W}})\right)\)
is a maximal point of \(\Spec
\left({\rm H}^0(\hat{Z},\ringO_{\hat{Z}})\right)\).
The restriction of the universal abelian scheme to \(\eta_{\hat{W}}\)
is an ordinary abelian variety.
Hence \(S=\Sigma_{F,p}\), otherwise \(A_{\eta_{\hat{W}}}\)
has slope \(1/2\) with multiplicity at least 
\(2 \sum_{\wp\notin S} [F_{\wp}:\Qp]\).
Theorem \ref{dens_hilb} is proved.
\hfill\qed
\medbreak

\noindent{\bf Remark.} The proof of Thm.\ \ref{dens_hilb} can be
finished without using Prop.\ \ref{thm:irredH} as follows.
We saw in the Remark after Step 3 that \(S_z\) depends only
on the irreducible component of \(Z_F^{\rm ord}\) which contains \(z\).
The argument in Step 4 shows that at least the subset 
\(S\subset \Sigma_{F,p}\) attached to one irreducible
component \(Z_1\) of \(Z_F^{\rm ord}\) is equal to \(\Sigma_{F,p}\).
So \(\dim(Z_1)=\dim(\calM_{F,\calL,\calL^+,n})=[F:\bbQ]\).
Since \(\calM_{F,\calL,\calL^+,n}\) is irreducible, we
conclude that \(Z_F=\calM_{F,\calL,\calL^+,n}\).
\bigbreak

\noindent{{\bf\scshape Proof of Theorem \ref{dens_siegel}.} (Density
of ordinary Hecke orbits in \(\calA_{g,1,n}\))}
\smallbreak

\noindent{\scshape Reduction step.}\enspace  By Thm.\ \ref{monconseq},
the weaker statement \ref{dens_siegel} (i) implies \ref{dens_siegel} (ii).
So it suffices to prove \ref{dens_siegel} (i).
\medbreak

\noindent{\bf Remark.} Our argument can be used to prove (ii)
directly without appealing to Thm.\ \ref{monconseq}.
But some statements, including the local stabilizer
principal for Hilbert modular varieties, need to be modified.
\medbreak

\noindent{{\bf\scshape Step 1} (Hilbert trick)}\enspace
Given \(x\in\calA_{g,n}(\Fpbar)\), 
Apply Thm.\ \ref{thm_hilb} to produce
a finite flat flat morphism
\[
g:\calM_{E,\calL,\calL^+,m;a}^{\rm ord}\to 
\calM_{E,\calL,\calL^+,m}^{\rm ord}\,,
\]
where \(E\) is a product of totally real number fields,
a finite morphism
\[
f:\calM_{E,\calL,\calL^+, m;a}^{\rm ord}\to\calA_{g,n}\,,
\]
and a point \(y_0\in\calM_{E,\calL,\calL^+,m;a}^{\rm ord}(\Fpbar)\) 
such that the following properties are satisfied.
\begin{itemize}
\item[(i)] There is a projective system 
\(\widetilde\calM_{E,\calL,\calL^+;a}^{\rm ord}\)
of finite \'etale coverings of 
\(\calM_{E,\calL,\calL^+,m;a}\) on which the group 
\(\SL_2(E\otimes\bbA_{f}^{(p)})\) operates.
This \(\SL_2(E\otimes\bbA_{f}^{(p)})\)-action induces 
Hecke correspondences on \(\calM_{E,\calL,\calL^+,m;a}^{\rm ord}\)

\item[(ii)] The morphism \(g\) is equivariant w.r.t.\ Hecke correspondences
coming from the group \(\SL_2(E\otimes\bbA_{f}^{(p)})\).
In other words, there is a \(\SL_2(E\otimes\bbA_f^{(p)})\)-equivariant 
morphism \(\tilde{g}\) from
the projective system 
\(\widetilde\calM_{E,\calL,\calL^+,a}^{\rm ord}\)
to the projective system 
\(\left(\calM_{E,\calL,\calL^+,mn}^{\rm ord}\right)_{n\in \bbN-p\bbN}\)
which lifts \(g\).

\item[(iii)] The finite morphism \(f\) is 
Hecke equivariant w.r.t.\ an injective homomorphism
\[
j_E:\SL_2(E\otimes_{\bbQ}\bbA_f^{(p)})\to \Sp_{2g}(\bbA_f^{(p)})
\,.
\]

\item[(iv)] For every geometric point 
\(z\in\calM_{E,\calL,\calL^+,m;a}^{\rm ord}\),
the abelian variety underlying the fiber over 
\(g(z)\in\calM_{E,\calL,\calL^+,m}^{\rm ord}\)
of the universal abelian scheme over \(\calM_{E,\calL,\calL^+,m}^{\rm ord}\)
is isogenous to the abelian variety underlying the fiber
over \(f(z)\in\calA_{g,n}^{\rm ord}(\Fpbar)\) of
the universal abelian scheme over \(\calA_{g,n}^{\rm ord}(\Fpbar)\).

\item[(v)] We have \(f(y_0)=x_0\). 
\end{itemize}
Let \(y:=g(y_0)\in \calM_{E,\calL,\calL^+,m}^{\rm ord}\).
\bigbreak

\noindent{\bf\scshape Step 2.} \enspace
Let \(Z_{y_0}\) be the Zariski closure of 
the \(\SL_2(E\otimes \bbA_f^{(p)})\)-Hecke orbit of 
\(y_0\) on \(\calM_{E,\calL,\calL^+,m;a}^{\rm ord}\),
and let \(Z_{y}\) be the Zariski closure of the 
\(\SL_2(E\otimes \bbA_f^{(p)})\)-Hecke orbit on
\(\calM_{E,\calL,\calL^+,m}\).
By Thm.\ \ref{dens_hilb}
we know that
\(Z_{y}=\calM_{E,\calL,\calL^+,m}\).
Since \(g\) is finite flat, 
we conclude that \(g(Z_{y_0})=Z_{y}\cap \calM_{E,\calL,\calL^+,m}^{\rm ord}
=\calM_{E,\calL,\calL^+,m}^{\rm ord}\).
We know that \(f(Z_{y_0})\subset Z_x\) because \(f\) is Hecke-equivariant.
\bigbreak

\noindent{\bf\scshape Step 3.}\enspace 
Let \(E_1\) be an ordinary elliptic curve over \(\Fpbar\).
Let \(y_1\) be an \(\Fpbar\)-point of \(\calM_{E,\calL,\calL^+,m}\)
such that \(A_{y_1}\) is isogenous to \(E_1\otimes_{\bbZ}\ringO_E\) 
and 
\(\calL_{y_1}\) contains the 
an \(\ringO_E\)-submodule of finite index in
\(\lambda_{E_1}\otimes\ringO_E\),
where \(\lambda_{E_1}\) denotes the canonical principal polarization  
on \(E_1\). 
In the above the tensor product \(E_1\otimes_{\bbZ}\ringO_E\)
is taken in the category of fppf sheaves over \(\Fpbar\);
the tensor product is represented by an abelian variety
isomorphic to the product of \(g\) copies of \(E_1\), 
with an action by \(\ringO_E\).
It is not difficult to check that such a point \(y_1\) exists.
\smallbreak

Let \(z_1\) be a point of \(Z_{y_0}\) such that
\(g(z_1)=y_1\). Such a point \(y_1\) exists because
\(g(Z_{y_0})=\calM_{E,\calL,\calL^+,m}^{\rm ord}\)
The point 
\(x_1=f(z_1)\) is contained in the Zariski closure
\(Z(x)\) of the prime-to-\(p\) Hecke orbit of \(x\) on 
\(\calA_{g,n}\).
Moreover \(A_{x_1}\) is isogenous to the product of 
\(g\) copies of \(E_1\) by property (iv) in Step 1.
So \(\End^0(A_{x_1})\cong {\rm M}_g(K)\), where
\(K=\End^0(E)\) is an imaginary quadratic extension field of
\(\bbQ\) which is split above \(p\).
The local stabilizer principle says that 
\(Z(x)^{/x_1}\) is stable under the natural action
of an open subgroup of 
\({\rm SU}(\End^0(A_{x_1}),\lambda_{x_1})(\Qp)\cong \GL_g(\Qp)\).
\bigbreak

\noindent{\bf\scshape Step 4.}\enspace
%
%
%
%
We know that \(Z(x)\) is smooth  at the ordinary point \(x\)
over \(k\),
so \(Z(x)^{/x}\) is reduced and irreducible.
By the local stabilizer principle \ref{lsp}, 
\(Z(x)^{/x}\) is stable under the natural action of the
open subgroup \({\rm H}_x\) of 
\({\rm SU}(\End^0(A_{x_1}),\lambda_{x_1})\) consisting of all
elements \(\gamma\in {\rm SU}(\End^0(A_{x_1}),\lambda_{x_1})(\Qp)\)
such that \(\gamma(A_{x_1}[p^{\infty}])=A_{x_1}[p^{\infty}]\).
By Thm.\ \ref{localrig}, \(Z(x)^{/x_1}\) is a formal subtorus of the 
formal torus \(\calA_{g,n}^{/x_1}\), which is
stable under the action of an open subgroup of
\({\rm SU}(\End^0(A_{x_1},\lambda_{x_1}))(\Qp)\cong \GL_g(\Qp)\).
\smallbreak

Let \(X_*\) be the cocharacter group of the
Serre-Tate formal torus \(\calA_{g,n}^{/x_1}\),
and let \(Y_*\) be the cocharacter group of the formal subtorus
\(Z(x)^{/x_1}\). Both \(X_*\) and \(X_*/Y_*\) are free \(\Zp\)-modules.
It is easy to see that the restriction to \(\SL_g(\Qp)\) of the
linear action of 
\({\rm SU}(\End^0(A_{x_1}),\lambda_{x_1})(\Qp)\cong \GL_g(\Qp)\) on
\(X_*\otimes_{\bbZ}\Qp\) is isomorphic to the second
symmetric product of the standard representation of 
\(\SL_g(\Qp)\).
It is well-known that the latter is an absolutely irreducible 
representation of \(\SL_g(\Qp)\).
Since the prime-to-\(p\) Hecke orbit of \(x\) is infinite,
\(Y_*\neq (0)\), hence \(Y_*=X_*\). 
In other words \(Z(x)^{/x_1}=\calA_{g,n}^{/x_1}\). Hence
\(Z(x)=\calA_{g,n}\) because \(\calA_{g,n}\) is irreducible.
\hfill\qed
\bigbreak

\begin{rem}\label{rem:alt}
We mentioned at the beginning of this section that 
there is an alternative argument for Step 4 of the
proof of Thm.\ \ref{dens_hilb},
which uses \cite{J} instead of Thm.\ \ref{EOconseq}, 
therefore independent of depends on \cite{FO-EO}.
We sketch the idea here; see \S8 of \cite{Chai-Fam} for more details.
\medbreak

We keep the notation in Step 3 of the proof of \ref{dens_hilb}.
Assume that \(S\neq \Sigma_{F,p}\).
Consider the universal \(\ringO_F\)-linear
abelian scheme \((A\to Z_F^{\rm ord}, \iota)\)
and the \((\ringO_F\otimes\Zp)\)-linear Barsotti-Tate group 
\((A\to Z_F^{\rm ord},\iota)[p^{\infty}]\) over the 
base scheme \(Z_F^{\rm ord}\), which is smooth over \(\bbF\).
We have a canonical decomposition of 
\(X_{\wp}:=A[p^{\infty}]\to Z_F^{\rm ord}\) as the fiber product over
\(Z_F^{\rm ord}\) of 
\(\ringO_{\wp}\)-linear Barsotti-Tate groups
\(A[\wp^{\infty}]\to Z_F^{\rm ord}\), where 
\(\wp\) runs through the finite set \(\Sigma_{F,p}\) of
all places of \(F\) above \(p\).
Let \(X_1\to Z_{F}^{\rm ord}\) (resp.\ \(X_2\to Z_F^{\rm ord}\))
be the fiber product over \(Z_{F}^{\rm ord}\)
of those \(X_{\wp}\)'s with \(\wp\in S\)
(resp.\ with \(\wp\notin S\)), so that we have
\(A[p^{\infty}]=X_1\times_{Z_F^{\rm ord}} X_2\).
\smallbreak

We know that for every closed point \(s\) of \(Z_F^{\rm ord}\),
the restriction to the formal completion \(Z_F^{/s}\) of the 
\(\left(\prod_{\wp\notin S}\ringO_{\wp}\right)\)-linear
Barsotti-Tate group \(X_2\to Z_F^{\rm ord}\) is constant.
This means that \(X_2\to Z_F^{\rm ord}\) is the twist of 
a constant \(\left(\prod_{\wp\notin S}\ringO_{\wp}\right)\)-linear
Barsotti-Tate group by a character 
\[
\chi:\pi_1^{\rm et}(Z_F^{\rm ord})\to 
\prod_{\wp\notin S}\ringO_{\wp}^{\times}\,.
\]
More precisely, one twists the \'etale part and toric part 
of the constant Barsotti-Tate group by \(\chi\) and \(\chi^{-1}\) 
respectively.
Consequently
\(\End_{\prod_{\wp\notin S}\ringO_{\wp}}(X_2)\supseteq
\prod_{\wp\notin S}(\ringO_{\wp}\times \ringO_{\wp})\).
\smallbreak

By the main results in \cite{J}, we have an isomorphism
\[
\End_{\ringO_F}(A/Z_F^{\rm ord})\otimes_{\bbZ}\Zp
\xrightarrow{\sim}
\End_{\ringO_F\otimes_{\bbZ}\Zp}(A[p^{\infty}]/Z_F^{\rm ord})
=\End_{\prod_{\wp\in S}\ringO_{\wp}}(X_1)
\times \End_{\prod_{\wp\notin S}\ringO_{\wp}}(X_2)\,.
\]
Since 
\(\End_{\prod_{\wp\notin S}\ringO_{\wp}}(X_2)\supseteq
\prod_{\wp\notin S}(\ringO_{\wp}^{\times}\times \ringO_{\wp}^{\times})\), 
we conclude that 
\(\End_{\ringO_F}(A/Z_F^{\rm ord})\otimes_{\bbZ}\bbQ\)
is a either a totally imaginary quadratic extension field of \(F\) or
a central quaternion algebra over \(F\).
This implies that the abelian scheme 
\(A\to Z_F^{\rm ord}\) as smCM, therefore it is isotrivial.
We have arrived at a contradiction because \(\dim(Z_F^{\rm ord})>0\)
by \ref{infho}.  Therefore \(S=\Sigma_{F,p}\).
\end{rem}


\section{Notations and some results we will use}\label{Notat}

\subsection{} {\bf Warning.} In most recent papers there is a distinction 
between an abelian variety defined over a field $K$ on the one hand, 
and $A \otimes_K K'$ over $K' \supsetneqq K$ on the other hand.  
The notation $\End(A)$ stands for the ring of endomorphisms of $A$ over $K$. 
This is the way Grothendieck taught us to choose our notation.

In pre-Grothendieck literature and in some modern papers there is a 
confusion between on the one hand $A/K$ and ``the same'' abelian variety over 
any extension field. In such papers there is a  confusion. 
Often it is not clear what is meant by ``a point on $A$'', 
the notation $\End_K(A)$ can stand for the 
``endomorphisms defined over $K$'', but then sometimes    
$\End(A)$ can stand for the ``endomorphisms defined over $\overline{K}$''.

Please adopt the Grothendieck convention that a scheme $T \to S$ is 
what it is, and any scheme obtained by base extension $S' \to S$ is 
denoted by $T \times_S S' = T_{S'}$, etc. For an abelian scheme $X \to S$ 
write $\End(X)$ for the endomorphism ring of $X \to S$ 
(old terminology ``endomorphisms defined over $S$''). 
Do not write  $\End_T(X)$ but $\End(X \times_S T)$.

\subsection{}\label{Sendo}
We write $\End(A)$ for the endomorphism ring of $A$ and 
$\End^0(A) = \End(A) \otimes_{\ZZ} \QQ$ for the endomorphism algebra of $A$. 
By Wedderburn's theorem every central simple algebra is a matrix algebra 
over a division algebra. If $A$ is $K$-simple the algebra $\End^0(A)$ is 
a division algebra; in that case we write: 
$$
\QQ \quad\subset\quad L_0 \quad\subset\quad L 
:= {\rm Centre}(D) \quad\subset\quad D = \End^0(A);
$$
here $L_0$ is a totally real field, and either $L = L_0$ or $[L:L_0] = 2$ 
and in that case $L$ is a CM-field. In case $A$ is simple  
$\End^0(A)$ is one of the four types in the Albert classification.
We write: 
$$
[L_0:\QQ] =e_0, \quad [L:\QQ] = e, \quad [D:L] = d^2.
$$

\subsection{}
The Rosati involution $\dag: D \to D$ is positive definite. \\
{\bf Definition.} A \emph{simple division algebra of finite degree 
over $\QQ$ with a positive definite involution}, i.e.\
an anti-isomorphism of order two which is positive 
definite, is called an \emph{Albert algebra}. \\
Applications to abelian varieties and the classification have been 
described by Albert,  \cite{Albert1},  \cite{Albert2}, \cite{Albert3}.

\subsection{}{\bf Albert's classification.}\label{Albert} 
Any Albert algebra belongs to one of the following types.

\vspace{2mm}\n
Type I($e_0$) \quad Here $L_0 = L = D$ is a totally real field. 

\vspace{2mm}\n
Type II($e_0$)\quad Here $d=2$, $e=e_0$, $\inv_v(D) = 0$ for 
all infinite $v$, and $D$ is an indefinite quaternion algebra over 
the totally real field  $L_0 = L$.

\vspace{2mm}\n
Type III($e_0$)\quad Here  $d=2$, $e=e_0$, $\inv_v(D) \not= 0$ for 
all infinite $v$, and $D$ is an definite quaternion algebra over the 
totally real field  $L_0 = L$.

\vspace{2mm}\n
Type IV($e_0,d$)\quad Here $L$ is a CM-field, $[F:\QQ] = e = 2e_0$, 
and $[D:L] = d^2$.

\subsection{}{\bf smCM}\label{smCM} We say that an abelian variety $X$ over 
a field $K$ {\it admits sufficiently many 
complex multiplications over} $K$, 
abbreviated by ``smCM over $K$", if $\End^0(X)$ contains a commutative 
semi-simple subalgebra of rank $2{\cdot}\dm(X)$ over $\QQ$. 
Equivalently: for every simple abelian variety $Y$ over $K$ which 
admits a non-zero homomorphism to $X$ the algebra $\End^0(Y)$ 
contains a field of degree $2{\cdot}\dm(Y)$ over $\QQ$. 
For other characterizations see \cite{Deligne-900}, page 63, 
see \cite{Mumford-Discgroups}, page 347.  

Note that if a simple abelian variety $X$ of dimension $g$ over a field 
{\it of characteristic zero} admits smCM then its endomorphism algebra 
$L = \End^0(X)$ is a  CM-field of degree $2g$ over $\QQ$.
We will use he terminology ``CM-type" in the case of an abelian variety $X$  
over $\CC$ which admits smCM, and where the type is given, i.e. 
the action of the endomorphism algebra on the tangent space 
$T_{X,0} \cong \CC^g$ is part of the data. 

Note however that there exist (many) abelian varieties $A$ admitting smCM 
(defined over a field of positive characteristic), such that $\End^0(A)$ 
is not a field. 

\vn
By Tate we know that an abelian variety over a finite field admits smCM, 
see \ref{struct2}. By Grothendieck we know that an abelian variety which 
admits smCM can be defined up to isogeny over a finite field, 
see \ref{GrothsmCM}. 

\vn
{\bf Terminology.} Let $\va \in \End^0(A)$. Then $d\va$ is a 
$K$-linear endomorphism of the tangent space. 
If the base field is $K = \CC$, this is just multiplication by 
a complex matrix $x$, and every  multiplication by a complex matrix  
$x$ leaving invariant the lattice $\Lambda$, where 
$A(\CC) \cong \CC^g/\Lambda$, gives rise to an endomorphism of $A$. 
If $g=1$, i.e. $A$ is an elliptic curve, and $\va \not\in \ZZ$ 
then $x\in\CC$ and $x \not\in \RR$. Therefore
an endomorphism of an elliptic curve over $\CC$ which is not in $\ZZ$ 
is sometimes called ``a complex multiplication''. 
Later this terminology was extended to all abelian varieties.\\
{\bf Warning.} Sometimes the terminology ``an abelian variety with CM'' 
is used, when one wants to say ``admitting smCM''. An elliptic curve 
$E$ has $\End(E) \supsetneqq \ZZ$ if and only if it admits smCM. 
However it is easy to give an abelian variety $A$ which ``admits CM'', 
meaning that $\End(A) \supsetneqq \ZZ$, such that $A$ does not admit smCM. 
However we will use the terminology ``a CM-abelian variety'' 
for an abelian variety which admits smCM.

\subsection{}{\bf Exercise.} {\it Show there exists an abelian variety 
$A$ over a field $k$ such that $\ZZ \subsetneqq \End(A)$ and such 
that $A$ does not admit} smCM.

\subsection{}{\bf Theorem}\label{struct2} (Tate). {\it  Let $A$ be 
an abelian variety over a finite field.} \\
{\bf (1)} {\it The algebra $\End^0(A)$ is semi-simple. 
Suppose $A$ is simple; the center of $\End^0(A)$  
equals $L := \QQ(\pi_A)$.}\\
{\bf (2)}  {\it Suppose $A$ is} simple; {\it then 
$$2g  \quad=\quad [L : \QQ]{\cdot}\sqrt{[D:L]},$$
where $g$ is the dimension of $A$. Hence: every abelian 
variety over a finite field admits {\rm smCM}.} See \ref{smCM}.
{\it We have:}
$$
f_A \quad=\quad \left({\rm Irr}_{\pi_A}\right)^{\sqrt{[D:L]}}.
$$
Here \(f_Z\) is the characteristic polynomial of the Frobenius 
morphism \({\rm Fr}_{A,\bbF_q}:A\to A\),
and \({\rm Irr}_{\pi_A}\) is the irreducible polynomial 
over \(\bbQ\) of the element \(\pi_A\) in the finite extension
\(L/\bbQ\).\\
{\bf (3)} {\it Suppose $A$ is} simple,
$$
\QQ \quad\subset\quad L := \QQ(\pi_A) \quad\subset\quad D = \End^0(A).
$$
{\it The central simple algebra $D/L$ 
\begin{itemize}
\item does not split at every real place of $L$, 
\item does split at every finite place not above $p$, 
\item and for $v \mid p$ the invariant of  $D/L$ is given by
\end{itemize}
$$
\inv_v(D/L) = \frac{v(\pi_A)}{v(q)}{\cdot}[L_v:\QQ_p] \quad{\rm mod} \ \ 1,
$$ 
where $L_v$ is the local field obtained from $L$ by completing at $v$.
}\\
See \cite{Tate-Endo}, \cite{Tate-Bourb}. 

\subsection{}{\bf Remark.} An abelian variety over a field of characteristic 
zero which admits smCM is defined over a number field; e.g. 
see \cite{Shim.Tan}, Proposition 26 on page 109.

\subsection{}{\bf Remark.}\label{GrothsmCM}  The converse of Tate's 
result \ref{struct2} (2) is almost true. Grothendieck showed: 
{\it Let $A$ be an abelian variety over a field $K$ which admits} 
smCM; {\it then $A_k$ is isogenous with an abelian variety 
defined over a finite extension of the prime field, 
where $k = \overline{K}$}; see \cite{FO-isog}.

It is easy to give an example of an abelian variety  (over a field 
of characteristic $p$), with smCM which is not defined over a finite field.  

\subsection{}{\bf Exercise.} {\it Give an example of a simple abelian 
variety $A$ over a field $K$ such that $A \otimes \overline{K}$ 
is not simple.}

\subsection{}\label{grsch}
We fix a prime number $p$. Base schemes and base fields will be 
of characteristic $p$, unless otherwise stated. 
We write $k$ or $\Om$ for an algebraically closed field 
of characteristic $p$.

\vn
We write $\cA_g$ for the moduli space of polarized abelian varieties 
of dimension $g$ in characteristic $p$ (this we write instead of 
$\cA_g \otimes \FF_p$). If we write ``work over $K$'' 
we mean that we consider $\cA_g \otimes K$. 
We write $\cA_{g,d}$ in case only polarizations of degree $d^2$ 
are considered. We write $\cA_{g,d,n}$ by 
considering polarized abelian varieties with a symplectic 
level-$n$-structure; in this case it is assumed that 
$n$ is not divisible by $p$.

\vn
The dimension of an abelian variety usually we will denote by $g$. 
If $m \in \ZZ_{>1}$ and $A$ is an abelian variety we write 
$A[m]$ for the the group scheme of $m$-torsion. 
Note that if $m$ is not divisible by $p$, 
then $A[m]$ is a group scheme \'etale over $K$; in this case 
it is uniquely determined by the Galois module $A[m](k)$. 
If $p$ divides $m$, then $A[m]$ is a group scheme which is not reduced. 

\vn
Group schemes considered will be assumed to be {\it commutative}. 
If $G$ is a finite abelian group, and $S$ is a scheme, 
we write $\underline{G}_S$ for the constant group scheme over $S$ 
with fiber equal to $G$.

\vn
Let $N \to S$ be a finite, flat group scheme. 
We write $N^D \to S$ for its Cartier dual, see \cite{FO-CGS}, I.2.

\subsection{} For the definition of an abelian variety, an abelian scheme,
see \cite{Mumford-AV},  II.4, \cite{GIT}, 6.1. 
The dual of an abelian scheme $A \to S$ will be denoted by $A^t \to S$, 
denoted by $\hat{A}$ in \cite{GIT}, 6.8.  

An isogeny $\va: A \to B$ of abelian schemes is a finite, 
surjective homomorphism. It follows that $\Ker(\va)$ is finite and 
flat over the base,  \cite{GIT}, Lemma 6.12. 
This defines a dual isogeny $\va^t: B^t \to A^t$. And see \ref{dualth}. 

\vn
A divisor $D$ on an abelian scheme $A/S$ defines a morphism 
$\va_D: A \to A^t$, see \cite{Mumford-AV}, 
theorem on page 125, see see \cite{GIT}, 6.2. 
A {\it polarization} on an abelian scheme $\mu: A \to A^t$ 
is an isogeny such that for every geometric point $s \in S(\Om)$ 
there exists an ample divisor $D$ on $A_s$ 
such that $\lambda_s = \va_D$, 
see \cite{Mumford-AV}, Application 1 on page 60, 
and \cite{GIT}, Definition 6.3. 
Note that a polarization is {\it symmetric} in the sense that
$$
\left(\lambda: A \to A^t \right) \quad=\quad 
\left(A \stackrel{\kappa}{\longrightarrow} 
A^{tt}\stackrel{\lambda^t}{\longrightarrow} A^t \right),
$$  
where $\kappa: A \to A^{tt}$ is the canonical isomorphism.

\vn
Writing $\va: (B,\mu) \to (A,\lambda)$ 
we mean that $\va: A \to A$         
and $\va^{\ast}(\lambda) = \mu$, i.e.
$$
\mu = \left(B \stackrel{\va}{\longrightarrow}  A   
\stackrel{\lambda}{\longrightarrow}  A^t   
\stackrel{\va^t}{\longrightarrow} B^t \right).
$$

\subsection{\bf The Frobenius morphism.} For a scheme $S$ over $\FF_p$ 
(i.e. $p{\cdot}1=0$ in all fibers of $\cO_S$), we define the 
absolute Frobenius morphism $\fr: S \to S$; if $S=\Spec(R)$ 
this is given by $x \mapsto x^p$ in $R$.

For a scheme $A \to S$ we define $A^{(p)}$ as the fiber 
product of $A \to S
\stackrel{\fr}{\longleftarrow} S$. The morphism 
$\fr: A \to A$ factors through $A^{(p)}$. 
This defines $F_A: A \to A^{(p)}$, a morphism over $S$; 
this is called {\it the relative Frobenius morphism}. 
If $A$ is a group scheme over $S$, the morphism 
$F_A: A \to A^{(p)}$ is a homomorphism of group schemes. 
For more details see \cite{SGA3}, Exp. VII$_A$.4. 
The notation $A^{(p/S)}$ is (maybe) more correct.

\vn
{\bf Example.} Suppose $A \subset \AA^n_R$ is given 
as the zero set of a polynomial $\sum_I a_IX^I$ 
(multi-index notation). Then $A^{(p)}$ is given by  
$\sum_I a^p_IX^I$, and $A \to A^p$ is given, 
on coordinates, by raising these to the power $p$. 
Note that if a point $(x_1, \cdots , x_n) \in A$ then 
indeed $(x^p_1, \cdots , x^p_n) \in A^{(p)}$, and 
$x_i \mapsto x^p_i$ describes  $F_A: A \to A^{(p)}$ 
on points.

\vn Let $S = \Spec(\FF_p)$; for any 
$T \to S$ we have a canonical isomorphism 
$T \cong T^{(p)}$. In this case $F_T = \fr: T \to T$.

\subsection{}{\bf Verschiebung.} Let $A$ be a 
{\it commutative} group scheme over 
a characteristic $p$ base scheme. 
In \cite{SGA3}, Exp. VII$_A$.4 we find the definition of 
the ``relative Verschiebung''
$$V_A:  A^{(p)} \to A; \quad\mbox{\rm we have:} 
\quad F_A{\cdot}V_A = [p]_{A^{(p)}}, \ \ V_A{\cdot}F_A = [p]_A.$$
In case $A$ is an abelian variety we see that $F_A$ is surjective, 
and $\Ker(F_A) \subset A[p]$. In this case we do not need 
the somewhat tricky construction of  \cite{SGA3}, Exp. VII$_A$.4, 
but we can define $V_A$ by $V_A{\cdot}F_A = [p]_A$ and check that 
$F_A{\cdot}V_A = [p]_{A^{(p)}}$.

\subsection{}\label{FV}
{\bf Remark.}  We use covariant Dieudonn\'e module theory. 
The Frobenius on a group scheme $G$ defines the Verschiebung on $\DD(G)$; 
this we denote by $\cV$, in order to avoid possible confusion. 
In the same way as  ``$\DD(F) = \cV$'' we have ``$\DD(V) = \cF$''. 
See \cite{FO-EO}, 15.3.

\subsection{\bf Algebraization.}\label{algebr}  
(1)  Suppose given a formal $p$-divisible group $X_0$ over $k$ with 
$\cN(X_0)=\gamma$ ending at $(h,c)$. We write $\cD = {\rm Def}(X_0)$ for 
the universal deformation space in equal characteristic $p$.
By this we mean the following. 
Formal deformation theory of $X_0$ is prorepresentable; 
we obtain a formal scheme $\Spf(R)$ and a prorepresenting family 
$\cX' \to \Spf(A)$.  
However ``a finite group scheme over a formal scheme actually 
is already defined over an actual scheme''. Indeed,
by \cite{AJdJ-Crys}, Lemma 2.4.4 on page 23, 
we know that there is an equivalence of categories 
of $p$-divisible groups 
over $\Spf(R)$ respectively over \Spec(R). 
We will say that $\cX \to \Spec(R) = \cD = {\rm Def}(X_0)$ 
is the universal deformation of $X_0$ 
if the corresponding $\cX' \to \Spf(R) = \cD^{\wedge}$ 
prorepresents the deformation functor.\\
Note that for a formal $p$-divisible group $\cX \to \Spf(R)$, where $R$ 
is moreover an integral domain, it makes sense to consider 
``the generic fiber'' of $X/\Spec(R)$. 

\vn
(2) Let $A_0$ be an abelian variety. The deformation functor 
${\rm Def}(A_0)$ is prorepresentable. We obtain the universal family 
$A \to \Spf(R)$, which is a formal abelian scheme. 
If $\dim(A_0) > 1$ this family is {\it not algebraizable}, 
i.e. it does not come from an actual scheme over $\Spec(R)$. 

\vn
(3) Let $(A_0,\mu_0)$ be a polarized abelian variety. 
The deformation functor ${\rm Def}(A_0,\mu_0)$ is prorepresentable. 
We can use the   Chow-Grothendieck theorem, 
see \cite{EGA}, III$^1$.5.4 (this is also called a theorem of 
``GAGA-type"):  
the  formal polarized abelian scheme obtained is algebraizable, 
and we obtain the universal deformation as a polarized abelian scheme 
over $\Spec(R)$.

\vn
The subtle differences between (1), (2) and (3) will be used without 
further mention.

\subsection{}{\bf Theorem}\label{irredA} (Irreducibility of moduli spaces.) 
\Bl \ {\it Let $K$ be a field, and consider $\cA_{g,1,n} \otimes K$ 
the moduli space of principally polarized abelian varieties over $K$-schemes, 
where $n \in \ZZ_{>0}$ is prime to the characteristic of $K$. 
This moduli scheme is geometrically irreducible.}\\
For fields of characteristic zero this follows by complex uniformization. 
For fields of positive characteristic  this was proved by Faltings in 1984, 
see \cite{Falt}, at the same time by Chai in his Harvard PhD-thesis; 
also see \cite{F.C}, IV.5.10. for a pure characteristic-$p$-proof see 
\cite{FO-EO}, 1.4.

\subsection\label{et} {\bf\'etale finite group schemes as 
Galois modules.} 
(Any characteristic.)
Let $K$ be a field, and let $G = \Gal(K^{\rm sep}/K)$. The main theorem 
of Galois theory says that there is an equivalence between the category 
of algebras\'etale and finite over $K$, and the category of finite sets 
with a continuous $G$-action. Taking group-objects on both sides 
we arrive at:\\
{\bf Theorem.} {\it There is an equivalence between the 
category of\'etale 
finite group schemes over $K$ and the category of 
finite continuous $G$-modules.} \\
See \cite{Waterhouse-Introd}, 6.4. 
Note that this equivalence also holds in the case of 
not necessarily commutative group schemes.

\vn
Naturally this can be generalized to: let $S$ be a connected scheme, 
and let $s \in S$ be a base point; let $\pi = \pi_1(S,s)$.  
{\it There is an equivalence between the category of
\'etale finite group schemes} 
(not necessarily commutative) 
{\it over $S$ and the category of finite continuous $\pi$-sets.}

\section{A remark and a question}\label{ques}
\subsection{} In \ref{HdenseNP} we have seen that the closure of the 
{\it full} Hecke orbit equals the related Newton polygon stratum. 
That result finds its origin in the construction of two {\it foliations}, 
as in \cite{FO-Fol}: Hecke-prime-to-$p$ actions ``move'' 
a point in a central leave, and Hecke actions only involving compositions 
of isogenies with kernel isomorphic with $\alpha_p$  ``move'' a point in 
an isogeny leaf, called $\cH_{\alpha}$-actions; as an open Newton polygon 
stratum, up to a finite map, is equal to the product of a central leaf 
and an isogeny leaf the result \ref{HdenseNP} for an irreducible component 
of a Newton polygon stratum follows if we show that $\cH_{\ell}(x)$ is dense 
in the central leaf passing through $x$. 

\vn
In case of ordinary abelian varieties the central leaf is the whole 
open Newton polygon stratum. As the Newton polygon goes up central leaves 
get smaller. Finally for supersingular points, a central leaf is 
finite and see \ref{ssfinite}, and  an isogeny leaf of a supersingular 
point is the whole supersingular locus. 

\vn
In order to finish a proof of   \ref{HdenseNP} one  shows that Hecke-$\alpha$ 
actions act transitively on the set of geometric components of the 
supersingular locus, and that any Newton polygon stratum in $\cA_{g,1}$ 
which is not supersingular is geometrically irreducible, 
see \cite{Chai.FO-Mon}. 

\subsection{} In Section \ref{5}, in particular see the proofs of \ref{G} 
and \ref{diagram}, we have seen a natural way of introducing coordinates 
in the formal completion at a point $x$ where $a \leq 1$ 
on an (open) Newton polygon stratum:
$$
\left(\cW_{\xi}(\cA_{g,1,n})\right)^{\wedge x} = \Spf(B_{\xi}),
$$ 
see the proof of \ref{pG}.
It would be nice to have a better understanding and interpretation of 
these ``coordinates''.

\vn
As in \cite{FO-Dimleaves} we write  
$$\triangle(\xi;\xi^{\ast}) := \{(x,y) \in \ZZ \mid (x,y) \prec \xi,   
\quad  (x,y) \succneqq \xi^{\ast}, \quad  x \leq g\}.$$
 We write 
$$B_{(\xi;\xi^{\ast})} = k[[Z_{x,y} \mid (x,y) 
\in \triangle(\xi;\xi^{\ast})]].$$
The inclusion $\triangle(\xi;\xi^{\ast}) \subset \triangle(\xi)$ defines 
$B_{\xi} \twoheadrightarrow B_{(\xi;\xi^{\ast})}$ by equating elements not in 
$\tr(\xi;\xi^{\ast})$ to zero.
Hence $\Spf(B_{(\xi;\xi^{\ast})}) \subset \Spf(B_{\xi})$.  
We also have the inclusion $\cC(x) \subset \cW_{\xi}(\cA_{g,1,n})$. \\
{\bf Question.} Does the inclusion 
$\triangle(\xi,\xi^{\ast})]]) \subset \triangle(\xi)$ 
define the inclusion 
$(\cC(x))^{\wedge x} \subset (\cW_{\xi}(\cA_{g,1,n}))^{\wedge x}$?
\\
A positive answer would give more insight in these coordinates, 
also along a central leaf, and perhaps a new proof of results in 
\cite{FO-Dimleaves}.


\hfill
\begin{tabbing}

        \= Department of MathematicsXXXXXX\= Mathematisch InstituutXXXX\= \kill 
        \=Ching-Li Chai \>Frans Oort \> \\ 
        \>Department of Mathematics \>Mathematisch Instituut \> \\
        \>University of Pennsylvania \>Budapestlaan 6 \> Postbus 80010\\   
       \>Philadelphia, PA 19104-6395 \>NL - 3584 CD TA Utrecht\> NL - 3508 TA Utrecht\\ 
       \>USA \>The Netherlands\> The Netherlands\\ 
        \>email:  chai@math.upenn.edu\>email:oort@math.uu.nl  \\
\end{tabbing}


\end{document}